\documentclass[11pt]{article}

\def\sqr#1#2{{\vcenter{\vbox{\hrule height.#2pt
              \hbox{\vrule width.#2pt height#1pt \kern#1pt \vrule width.#2pt}
              \hrule height.#2pt}}}}
\def\signed #1{{\unskip\nobreak\hfil\penalty50
              \hskip2em\hbox{}\nobreak\hfil#1
              \parfillskip=0pt \finalhyphendemerits=0 \par}}
\def\endpf{\signed {$\sqr69$}}
\def\dbE{{\mathop{\rm l\negthinspace E}}}
\def\dbR{{\mathop{\rm l\negthinspace R}}}

\def\3n{\negthinspace \negthinspace \negthinspace }
\def\1n{\negthinspace }

\def\dbE{{\mathop{\rm l\negthinspace E}}}

\def\dbN{{\mathop{\rm l\negthinspace N}}}

\def\dbR{{\mathop{\rm l\negthinspace R}}}
\def\={\buildrel \triangle \over =}

\def\ds{\displaystyle}

\def\ns{\noalign{\ss}}

\def\a{\alpha}
\def\b{\beta}

\def\th{\theta}
\def\o{\omega}

\def\O{\Omega}

\def\ss{\smallskip}
\def\ms{\medskip}

\def\q{\quad}
\def\qq{\qquad}

\def\cd{\cdot}

\def\supp{\hbox{\rm supp$\,$}}

\def\Re{{\mathop{\rm Re}\,}}
\def\Im{{\mathop{\rm Im}\,}}

\def\|{\Big |}
\def\({\Big (}
\def\){\Big )}
\def\[{\Big[}
\def\]{\Big]}
\def\bde{\begin{definition}}
\def\ede{\end{definition}}
\def\be{\begin{equation}}
\def\bel{\begin{equation}\label}
\def\ee{\end{equation}}
\def\bt{\begin{theorem}}
\def\et{\end{theorem}}
\def\bc{\begin{corollary}}
\def\ec{\end{corollary}}
\def\bl{\begin{lemma}}
\def\el{\end{lemma}}
\def\bp{\begin{proposition}}
\def\ep{\end{proposition}}
\def\bas{\begin{assumption}}
\def\eas{\end{assumption}}
\def\br{\begin{remark}}
\def\er{\end{remark}}
\def\ba{\begin{array}}
\def\ea{\end{array}}
\def\ed{\end{document}}

\def\square#1{\vbox{\hrule\hbox{\vrule height#1%
     \kern#1\vrule}\hrule}}
\def\rectangle#1#2{\vbox{\hrule\hbox{\vrule height#1%
     \kern#2\vrule}\hrule}}

\font\tenbb=msbm10 \font\sevenbb=msbm7 \font\fivebb=msbm5

\newfam\bbfam
\scriptscriptfont\bbfam=\fivebb \textfont\bbfam=\tenbb
\scriptfont\bbfam=\sevenbb

\newtheorem{lemma}{Lemma}[section]
\newtheorem{remark}{Remark}[section]

\newtheorem{theorem}{Theorem}[section]
\newtheorem{corollary}{Corollary}[section]

\newtheorem{definition}{Definition}[section]
\newtheorem{proposition}{Proposition}[section]
\newtheorem{assumption}{Assumption}[section]

\makeatletter
   
   \@addtoreset{equation}{section}
\makeatother

\begin{document}

\title{\bf The Theory of Stochastic Pseudo-differential Operators
and Its Applications, I\thanks{This work is partially supported by
the NSF of China under grants 10901032 and 10831007, and by the
National Basic Research Program of China (973 Program) under grant
2011CB808002.}}

\date{}

\author{Xu Liu\thanks{School of Mathematics and Statistics, Northeast Normal
University, Changchun 130024, China. E-mail address:
liuxu@amss.ac.cn. \ms } \q and \q Xu Zhang
\thanks{Key Laboratory of Systems and Control, Academy of
Mathematics and System Science, Chinese Academy of Sciences, Beijing
100190, China; Yangtze Center of Mathematics, Sichuan University,
Chengdu 610064, China. E-mail address: xuzhang@amss.ac.cn.}}

\date{}

\maketitle

\begin{abstract}
The purpose of this paper is to establish the theory of stochastic
pseudo-differential operators and give its applications in
stochastic partial differential equations. First, we introduce some
concepts on stochastic pseudo-differential operators and prove their
fundamental properties. Also, we present the boundedness theory,
invertibility of stochastic elliptic operators and the G{\aa}rding
inequality. Moreover, as an application of the theory of stochastic
pseudo-differential operators, we give a Calder\'on-type uniqueness
theorem on the Cauchy problem of stochastic partial differential
equations. The proof of the uniqueness theorem is based on a new
Carleman-type estimate, which is adapted to the stochastic setting.
\end{abstract}

\section{Introduction}

During the past seventy years, more and more studies of stochastic
phenomena have appeared. Research in this area is stimulated by the
need to take account of random effects in the engineering and
physical systems. For such systems, stochastic processes give a
natural replacement for deterministic functions as mathematical
descriptions. Since K. It\^{o} introduced the stochastic integral,
one of the main topics in Probability Theory and Stochastic Process
has been the issue on stochastic differential equations. Compared to
 deterministic differential equations, stochastic differential
equations are much more complicated. Indeed,
 one has to distinguish forward stochastic
differential equations, backward stochastic differential equations
and forward-backward stochastic differential equations in most of
the interesting cases.  Now, stochastic ordinary differential
equations have been well-developed. However, stochastic partial
differential equations (SPDEs for short) make slow progress. By now,
the main tool in this field is of functional analysis nature. As far
as we know, the real PDE-based approach is not well-developed.
Meanwhile, we notice that the theory of pseudo-differential
operators has been a powerful tool in the study of general partial
differential equations, since it was established in 1960's. It plays
a crucial role in the studies of existence, uniqueness and
propagation of singularities for the solutions of partial
differential equations. Therefore, we would like to introduce
 such a theory to the stochastic setting and regard it as a tool to solve
 the problems related to SPDEs.

For this purpose, we establish the theory of stochastic
pseudo-differential operators (SPDOs for short) and give some
applications to SPDEs. First of all, we introduce some basic
notions, which are adapted to the stochastic setting, including
symbol, amplitude, SPDO, kernel and uniformly properly supported
SPDO. Since a stochastic process has the variables on time and
sample point, we add these two variables to all notions and endow
them with suitable integrability. Also, in order for the symbol
calculus, we study asymptotic expansions of a symbol, and eventually
establish an algebra and generalized module of SPDOs. On the other
hand, we establish the $L^p$- boundedness theory, invertibility of
stochastic elliptic operators and the G{\aa}rding inequality, which
are some fundamental results related to the energy estimates.
Moreover, as an application of theory of SPDOs, we present a
Calder\'on-type uniqueness theorem on the Cauchy problem of SPDEs.
In his remarkable paper \cite{1}, A.-P. Calder\'on established a
fundamental result on the uniqueness of the non-characteristic
Cauchy problem for general partial differential equations. One of
the main tools introduced in \cite{1} is a preliminary version of
the symbol calculation technique. Later, Calder\'on's uniqueness
theorem was extended to the operators with characteristics of high
multiplicity. We refer to \cite{9} and the references cited therein
for some deep results in this topic. However, as far as we know,
there is no work addressing the uniqueness on the Cauchy problem for
general SPDEs. In this paper, we give a Calder\'on-type uniqueness
result in the stochastic setting, by virtue of the theory of SPDOs.
In order to present the key idea in the simplest way, we do not
pursue the full technical generality. More precisely, we focus
mainly on the Cauchy problem for SPDEs in the case of at most double
characteristics.

It is a little surprising that the theory of SPDOs was not available
in the previous literatures although a related but clearly different
theory for random pseudo-differential operators was introduced in
\cite{DS}. It deserves to point out that the study of SPDOs seems to
be of independent interest. We divide our results on SPDOs into two
parts. In this paper, we present the first part and its
applications. We will establish the stochastic micro-local analysis
and develop singularity propagation theory for stochastic hyperbolic
equations of second order in the forthcoming paper.

The rest of this paper is organized as follows. Section 2 is devoted
to the basic concepts and properties of SPDOs. In Section 3 we
establish the boundedness theory. In Section 4 we give invertibility
of stochastic elliptic operators and the G{\aa}rding inequality. As
an application of the theory of SPDOs, a Calder\'{o}n-type
uniqueness theorem is presented in Section 5.

\section{Calculus of stochastic pseudo-differential
operators}

Pseudo-differential operators developed from the theory of singular
integral operators, which were essentially pseudo-differential
operators with homogeneous symbol of order 0. The appearance of both
Calder\'on's uniqueness theorem (\cite{1}) and the index theorem for
elliptic operators by M. F. Atiyah and I. M. Singer (\cite{as})
showed the importance of the theory of singular integral operators.
Shortly afterwards, J. J. Kohn and L. Nirenberg (\cite{kn}) removed
the restriction to order 0 and generalized the notions of
pseudo-differential operators to the case of general polyhomogeneous
symbols. Later, L. H\"{o}rmander (\cite{h}) introduced
pseudo-differential operators with the symbols of type $(\rho,
\delta)$, by the need to incorporate fundamental solutions of
hypoelliptic operators of constant strength. Since the theory of
pseudo-differential operators was established in 1960's, it has been
an important mathematical branch (\cite{5}, \cite{7}). It plays an
important role in many fields, such as partial differential
equations, harmonic analysis and differential geometry etc.

In order to introduce the theory of pseudo-differential operators to
the stochastic setting, in this section, we shall present some basic
concepts and properties of SPDOs. First, we introduce some locally
convex topological vector spaces, which will be used later. Then, we
give the notions of symbol, amplitude, SPDO, kernel and uniformly
properly supported SPDO in sequence. Moreover, we give asymptotic
expansions of a symbol and establish an algebra and generalized
module of SPDOs.

To begin with, we give some usual notations. Throughout this paper,
$(\Omega, \mathcal{F}, \{\mathcal{F}_t\}_{t\geq 0}, P)$ is a
complete filtered probability space, on which a one dimensional
standard Brownian motion $\{w(t)\}_{t\ge 0}$ is defined.  Fix $T>0$,
$n\in \dbN\backslash\{0\}$, $m\in\dbN$, $\ell\in\dbR$, $p,\ q\in[1,
\infty]$ and a domain $G$ of $\dbR^n$. $i$ denotes the imaginary
unit. Let $H$ be a Banach space. We denote by $L^p_{\mathcal{F}}(0,
T; H)$ the set of all $H$-valued $\{\mathcal{F}_t\}_{t\geq
0}$-adapted process $X(\cdot)$ such that $\dbE\displaystyle\int^T_0
|X(t)|^p_H dt<\infty$; by $L_{\mathcal{F}}^\infty(0,T;H)$ the Banach
space consisting of all $H$-valued $\{\mathcal{F}_t\}_{t\ge
0}$-adapted bounded processes; and  by
$L_{\mathcal{F}}^p(\O;C^m([0,T];H))$ the Banach space consisting of
all $H$-valued $\{\mathcal{F}_t\}_{t\ge 0}$-adapted $m$-th order
continuously differential  processes $X(\cd)$ such that
$\dbE(|X(\cd)|_{C^m([0,T];H)}^p)<\infty$. Moreover, we simply write
$L^p_{\mathcal{F}}(0, T; \dbR)$ as $L^p_{\mathcal{F}}(0, T)$, and
have the similar notations for $L^\infty_\mathcal{F}(0, T)$ and
$L^p_\mathcal{F}(\Omega; C^m([0, T]))$. Furthermore, we denote by
$C(\cdot)$ a generic constant, which may be different from one place
to another.

\subsection{Basic function spaces}

In this subsection, we introduce some locally convex topological
vector spaces, which will be used as the domain or range of a SPDO
later. To begin with, we denote by $\mathcal{D}(G)$ the topological
space (with the usual inductive topology, see Page 54 in \cite{sch})
of infinitely differentiable functions supported by $G$; by
$\mathcal{E}(G)$ the topological space of infinitely differentiable
functions defined on $G$; and by $\mathcal{S}$ the topological space
of rapidly decreasing functions. Let $\{K_j\}_{j\in\dbN}$ stand for
a sequence of compact sets satisfying $K_0\subseteq
K_1\subseteq\cdots$ and $\displaystyle\bigcup_{j\in\dbN}K_j=G$.
Then, we write
\begin{eqnarray*}
&&|v|_{j, k, 1}=\sup\limits_{x\in K_j,\ |\alpha|\leq
k}|\partial_x^\alpha v(x)|,\ \ (j,\ k\in\dbN);\\[4mm]
&&|u|_{p, j, k, 1}=\left|\sup\limits_{x\in K_j,\ |\alpha|\leq
k}|\partial_x^\alpha u(\cdot, \cdot,
x)|\right|_{L^p_{\mathcal{F}}(0, T)},\ \ (j,\ k\in\dbN);\\[1mm]
&&|u|_{p, j, k, 2}=\left|\sup\limits_{x\in \dbR^n,\atop |\alpha|\leq
j, |\beta|\leq k}|x^\alpha\partial_x^\beta u(\cdot, \cdot,
x)|\right|_{L^p_{\mathcal{F}}(0, T)},\ \ (j,\ k\in\dbN),
\end{eqnarray*}
where $\alpha$ and $\beta$ are multi-indices. On the other hand, by
Theorem 1.36 and Theorem 1.35 in \cite{ru}, if $\mathcal{B}$ is a
0-neighborhood base for the inductive topology on $\mathcal{D}(G)$,
then for any $\gamma\in\mathcal{B}$,
$$
\gamma=\{\ v\in\mathcal{D}(G);\ \mu_\gamma(v)<1\ \},
$$
where $\mu_\gamma$ is the Minkowski functional of $\gamma$. Also,
$\{\mu_\gamma\}_{\gamma\in\mathcal{B}}$ is a family of generating
semi-norms on $\mathcal{D}(G)$. Set
$$
|u|_{p, \gamma}=|(\mu_\gamma(u))(\cdot,
\cdot)|_{L^p_{\mathcal{F}}(0, T)},\ \ (\gamma\in\mathcal{B}).
$$
Next, we define the following locally convex spaces:
\begin{eqnarray*}
&&L^p_{\mathcal{F}}(0, T; \mathcal{E}(G))=\left\{u\ |\ u(t, \omega,
\cdot)\in \mathcal{E}(G),\ \mbox{a.e. }(t, \omega)\in (0,
T)\times\Omega;\
u(\cdot, \cdot, x)\mbox{ is }\{\mathcal{F}_t\}_{t\geq 0}\mbox{-adapted},\right.\\
&&\left.\quad\quad\quad\quad\quad\quad\quad\quad\quad\quad
\mbox{for any } x\in G;\ \mbox{and }|u|_{p, j, k, 1}<\infty,\ j, k=0, 1, \cdots \right\},\\[1mm]
&&\mbox{which is a Fr\'{e}chet space, generated by a  sequence of
semi-norms
}\{|\cdot|_{p, j, k, 1}\}_{j, k\in\dbN};\\[1mm]
&&L^p_{\mathcal{F}}(0, T; \mathcal{S})=\left\{u\ |\ u(t, \omega,
\cdot)\in \mathcal{S},\ \mbox{a.e. }(t, \omega)\in (0,
T)\times\Omega;\ u(\cdot, \cdot, x)\mbox{ is
}\{\mathcal{F}_t\}_{t\geq
0}\mbox{-adapted},\right.\\
&&\left.\quad\quad\quad\quad\quad\quad\quad\quad \mbox{for any }x\in
\dbR^n; \mbox{ and }|u|_{p, j, k, 2}<\infty,\ j, k=0, 1, \cdots\right\},\\[1mm]
&&\mbox{which is a Fr\'{e}chet space, generated by a  sequence of
semi-norms
}\{|\cdot|_{p, j, k, 2}\}_{j, k\in\dbN};\\[1mm]
&&L^p_{\mathcal{F}}(0, T; \mathcal{D}(G))=\left\{u\ |\ u(t, \omega,
\cdot)\in \mathcal{D}(G),\ \mbox{a.e. }(t, \omega)\in (0,
T)\times\Omega;\ u(\cdot, \cdot, x)\mbox{
is }\{\mathcal{F}_t\}_{t\geq 0}\mbox{-adapted},\right.\\
&&\left.\quad\quad\quad\quad\quad\quad\quad\quad\quad\quad \mbox{for
any }x\in
G; \mbox{ and }|u|_{p, \gamma}<\infty,\ \mbox{for each }\gamma\in\mathcal{B}\right\},\\[1mm]
&&\mbox{which is generated by a  family of semi-norms }\{|\cdot|_{p,
\gamma}\}_{\gamma\in\mathcal{B}}.
\end{eqnarray*}

In addition, for any compact set $K\subseteq G$ and $j\in\dbN$, we
write
$$
|v|_{K, j}=\sup\limits_{x\in K,\ |\alpha|\leq j}|\partial^\alpha_x
v(x)|,\quad\quad\quad |u|_{p, K, j}=\left| |u(\cdot, \cdot,
\cdot)|_{K, j} \right|_{L^p_{\mathcal{F}}(0, T)}.
$$
Then, $\mathcal{D}_K=\{\ v\in C^\infty(G)\ |\ \supp v\subseteq K\}$
is a locally convex topological vector space, endowed with a
sequence of semi-norms $\{|\cdot|_{K, j}\}_{j\in\dbN}$. Also, we
introduce the following locally convex spaces:
\begin{eqnarray*}
&&L^p_{\mathcal{F}}(0, T; \mathcal{D}_K)=\left\{u\ |\ u(t, \omega,
\cdot)\in \mathcal{D}_K,\ \mbox{a.e. }(t, \omega)\in (0,
T)\times\Omega;\ u(\cdot, \cdot, x)\mbox{
is }\{\mathcal{F}_t\}_{t\geq 0}\mbox{-adapted},\right.\\
&&\left.\quad\quad\quad\quad\quad\quad\quad\quad\quad \mbox{for any
}x\in
G; \mbox{ and }|u|_{p, K, j}<\infty,\ j=0, 1, \cdots \right\},\\[1mm]
&&\mbox{which is generated by a sequence of semi-norms
}\{|\cdot|_{p, K, j}\}_{j\in\dbN};\\[3mm]
&&L^p_{\mathcal{F}}(0, T; \mathcal{D}_G)=\bigcup_{K\subseteq G
\mbox{ is compact}}L^p_{\mathcal{F}}(0, T; \mathcal{D}_K),\mbox{
which is endowed with the inductive topology.}
\end{eqnarray*}

In order to characterize the topology and convergence in
$L^p_{\mathcal{F}}(0, T; \mathcal{D}_G)$, we first recall the
following known result.
\begin{lemma}\label{103}{\rm (\cite[Page 54]{sch})}
Let $\mathcal{Z}$ and $\mathcal{Z}_K$ $(K\in \Theta,\ \Theta\mbox{
is an index set})$ be vector spaces, let $g_K$ be a linear mapping
of $\mathcal{Z}_K$ into $\mathcal{Z}$, and let $\Gamma_K$ be a
locally convex topology on $\mathcal{Z}_K$. If we denote by $\Gamma$
the inductive topology on $\mathcal{Z}$ with respect to the family
$\{(\mathcal{Z}_K, \Gamma_K, g_K); K\in\Theta\}$, then a
0-neighborhood base for $\Gamma$ is given by the family $\{U\}$ of
all balanced, convex, absorbing subsets of $\mathcal{Z}$, such that
for each $K\in\Theta$, $g^{-1}_K(U)$ is a 0-neighborhood in
$(\mathcal{Z}_K, \Gamma_K)$.
\end{lemma}
For our problem, for any compact set $K\subseteq G$, we take
$\mathcal{Z}=L^p_{\mathcal{F}}(0, T; \mathcal{D}_G)$ and
$\mathcal{Z}_K=L^p_{\mathcal{F}}(0, T; \mathcal{D}_K)$. $g_K:
\mathcal{Z}_K\rightarrow \mathcal{Z}$ is canonical embedding and the
topology $\Gamma_K$ is generated by a sequence of semi-norms
$\{|\cdot|_{p, K, j}\}_{j\in\dbN}$. Then for the inductive topology
$\Gamma$ on $L^p_{\mathcal{F}}(0, T; \mathcal{D}_G)$, a
0-neighborhood base $\mathcal{B}$ is given by the family $\{U\}$ of
all balanced, convex, absorbing subsets of $L^p_{\mathcal{F}}(0, T;
\mathcal{D}_G)$, such that for any compact set $K\subseteq G$,
$g_K^{-1}(U)$ is a 0-neighborhood in $(L^p_{\mathcal{F}}(0, T;
\mathcal{D}_K), \Gamma_K)$.

\medskip

Moreover, we recall another known result on locally convex spaces.
\begin{lemma}\label{106}{\rm (\cite[Theorem 1.37]{ru})}
Suppose that $\mathcal{P}$ is a separating family of semi-norms on a
vector space $\mathcal{Z}$. For each $\phi\in\mathcal{P}$ and
positive integer $k$, set
$$
V(\phi, k)=\left\{\ u\in \mathcal{Z};\ \phi(u)<\frac{1}{k}\right\}.
$$
Let $\mathcal{B}$ be the collection of all finite intersections of
the  sets $V(\phi, k)$. Then $\mathcal{B}$ is a balanced convex
absorbing local base, which turns $\mathcal{Z}$ into a locally
convex space.
\end{lemma}
For any compact set $K$, by Lemma \ref{106} and the definition of
semi-norms $\{|\cdot|_{p, K, j}\}_{j\in\dbN}$ on
$L^p_{\mathcal{F}}(0, T; \mathcal{D}_K)$, it is easy to show that
for any $\{u_k\}_{k\in\dbN}\subseteq L^p_{\mathcal{F}}(0, T;
\mathcal{D}_K)$, $\lim\limits_{k\rightarrow \infty}u_k=0$ in
$L^p_{\mathcal{F}}(0, T; \mathcal{D}_{K})$ if and only if for any
$j\in\dbN$,  $\lim\limits_{k\rightarrow
\infty}\dbE\displaystyle\int^T_0 \sup\limits_{x\in K,\ |\alpha|\leq
j}|\partial^\alpha_x u_k(t, \omega, x)|^pdt= 0$.

\medskip

In the remainder of this subsection, we give the following result on
the convergence in $L^p_{\mathcal{F}}(0, T; \mathcal{D}_G)$.
\begin{proposition}\label{104}
For any $\{u_j\}_{j\in\dbN}\subseteq L^p_{\mathcal{F}}(0, T;
\mathcal{D}_G)$, $\lim\limits_{j\rightarrow \infty}u_j=0$ in
$L^p_{\mathcal{F}}(0, T; \mathcal{D}_G)$ if and only if the
following two conditions hold:

\noindent $(1)$ there exists a compact set $K_*$ such that $\supp
u_j(t, \omega, \cdot)\subseteq K_*$ for a.e. $(t, \omega)\in(0,
T)\times\Omega$ and any $j\in\dbN$;

\noindent $(2)$ for any $k\in\dbN$, $\lim\limits_{j\rightarrow
\infty}\dbE\displaystyle\int^T_0 \sup\limits_{x\in K_*,\
|\alpha|\leq k}|\partial^\alpha_x u_j(t, \omega, x)|^pdt= 0$.
\end{proposition}
{\bf Proof. }By Lemma \ref{103} and Lemma \ref{106}, we have only to
 prove (1). Assume the contrary. Then there exist a subsequence of
 $\{u_{j_k}\}_{k\in\dbN}$ of $\{u_j\}_{j\in\dbN}$ and a sequence of
compact sets $\{K_{k,
*}\}_{k\in\dbN}$ satisfying $K_{0,
*}=\emptyset$, $K_{1, *}\subseteq K_{2, *}\subseteq\cdots$ and
$\bigcup\limits_{k\in\dbN}K_{k, *}=G$, such that
\begin{equation}\label{105}
|u_{j_k}(t, \omega, x_k)|\geq \varepsilon_k,\quad\quad k=1, 2,
\cdots,
\end{equation}
for $(t, \omega)\in T_k\times\Omega_k$ and a positive constant
$\varepsilon_k$, where $x_k\in K_{k, *}\backslash K_{k-1, *}$,
$\{T_k\}_{k=1, 2, \cdots}$ and $\{\Omega_k\}_{k=1, 2, \cdots}$ are
two sequences of measurable sets (with positive measures) of $(0,
T)$ and $\Omega$, respectively. Then, we define a semi-norm
$|\cdot|_{p,
*}$ on $L^p_{\mathcal{F}}(0, T; \mathcal{D}_G)$ $(p\geq 1)$ as follows:
\begin{eqnarray*}
&&|u|_{p,
*}^p=\sum_{k=1}^{\infty}\frac{1}{P(\Omega_k)|T_k|}\int_{\Omega_k}
\int_{T_k}\sup\limits_{x\in K_{k, *}\backslash
 K_{k-1, *}} \left|\frac{u(t, \omega, x)}{u_{j_k}(t, \omega,
 x_k)}\right|^pdtdP,
\end{eqnarray*}
where $|T_k|$ denotes the Lebesgue measure of $T_k$. Notice that for
any $u\in L^p_{\mathcal{F}}(0, T; \mathcal{D}_G)$, the right side of
the above equality is indeed a finite sum. Since $x_k\in K_{k,
*}\backslash K_{k-1, *}$, it is easy to see that $|u_{j_k}|_{p,*}\geq 1$
for any $k=1, 2, \cdots$. Therefore, if we write $U_*=\{\ u\in
L^p_{\mathcal{F}}(0, T; \mathcal{D}_G)\ |\ |u|_{p,
*}<1\}$, then any $u_{j_k}$ $(k=1, 2, \cdots)$ does not belong to $U_*$.

On the other hand, $U_*$ is a 0-neighborhood in
$L^p_{\mathcal{F}}(0, T; \mathcal{D}_G)$. In fact, by Lemma
\ref{103}, it remains to prove that for any compact set $K$,
$V_K=U_*\cap L^p_{\mathcal{F}}(0, T; \mathcal{D}_K)$ is a
0-neighborhood in $L^p_{\mathcal{F}}(0, T; \mathcal{D}_K)$. By the
definition of the semi-norm $|\cdot|_{p, *}$, it follows that
\begin{eqnarray*}
&&V_K=\left\{\ u\in L^p_{\mathcal{F}}(0, T; \mathcal{D}_K)\ \left|\
\sum_{k=1}^{\infty}\frac{1}{P(\Omega_k)|T_k|}\int_{\Omega_k}
\int_{T_k}\sup\limits_{x\in (K_{k, *}\backslash
 K_{k-1, *})\cap K} \left|\frac{u(t, \omega, x)}{u_{j_k}(t, \omega,
 x_k)}\right|^pdtdP<1\right\}.\right.
\end{eqnarray*}
We suppose that $K$ has nonempty intersections with the sets $K_{k,
*}\backslash
 K_{k-1, *}$ $(k=1, 2, \cdots, i_*)$, where $i_*=i_*(K)$ is a positive integer.
Then, for any $u\in L^p_{\mathcal{F}}(0, T; \mathcal{D}_K)$,
\begin{eqnarray}\label{108}
\begin{array}{rl}
&\displaystyle\sum\limits_{k=1}^{\infty}\frac{1}{P(\Omega_k)|T_k|}\displaystyle\int_{\Omega_k}
\int_{T_k}\sup\limits_{x\in (K_{k, *}\backslash
 K_{k-1, *})\cap K} \left|\frac{u(t, \omega, x)}{u_{j_k}(t, \omega,
 x_k)}\right|^pdtdP\\[2mm]
&=\displaystyle\sum\limits_{k=1}^{i_*}\displaystyle\frac{1}{P(\Omega_k)|T_k|}\int_{\Omega_k}
\int_{T_k}\sup\limits_{x\in (K_{k, *}\backslash
 K_{k-1, *})\cap K} \left|\frac{u(t, \omega, x)}{u_{j_k}(t, \omega,
 x_k)}\right|^pdtdP\\[2mm]
 &\leq \displaystyle\sum\limits_{k=1}^{i_*}\displaystyle\frac{1}{P(\Omega_k)|T_k|\varepsilon^p_k}\dbE
\int_0^{T}\sup\limits_{x\in K} \left|u(t, \omega, x)\right|^pdt.
\end{array}
\end{eqnarray}
If we take $N_*=1+\left[\displaystyle\sum_{k=1}^{i_*}\displaystyle
\frac{1}{P(\Omega_k)|T_k|\varepsilon^p_k}\right]$, then by
(\ref{108}),
\begin{eqnarray*}
&&\left\{\ u\in L^p_{\mathcal{F}}(0, T; \mathcal{D}_K)\ \left|\
\dbE\int^T_0\sup\limits_{x\in K}|u(t, \omega,
x)|^pdt<\frac{1}{N_*}\right\}\subseteq V_K,\right.
\end{eqnarray*}
where $[\ell]$ denotes the integral part of a real number $\ell$.
Therefore, $V_K$ is a 0-neighborhood of $L^p_{\mathcal{F}}(0, T;
\mathcal{D}_K)$. This implies that $\{u_j\}_{j\in\dbN}$ cannot be a
sequence converging to 0 in $L^p_{\mathcal{F}}(0, T;
\mathcal{D}_G)$. This contradiction proves that (1) must be true.
Similarly, we can get the desired result for $p=\infty$. The proof
is completed.\endpf

\subsection{Symbol and stochastic pseudo-differential operators}

In this subsection, we use the Fourier integral representation to
define SPDOs. For this purpose,
 we first introduce the notion of symbols. Compared to the
 classical one in the deterministic case, we add two variables $t$ and $\omega$, and endow them
 with the integrability.
\begin{definition}\label{d1}
A complex-valued function $a$ is called a symbol of order $(\ell,
p)$ if $a$ satisfies the following conditions:\\[2mm]
\noindent $(1)$ $a(t, \omega, \cdot, \cdot)\in C^\infty(G\times
\dbR^n)$,
a.e. $(t, \omega)\in (0, T)\times \Omega$;\\[1mm]
\noindent $(2)$ $a(\cdot, \cdot, x, \xi)$ is
$\{\mathcal{F}_t\}_{t\geq
0}$-adapted, $\forall\ (x, \xi)\in G\times \dbR^n$;\\[1mm]
\noindent $(3)$ for any two multi-indices $\alpha$ and $\beta$, and
any compact set $K\subseteq G$, there exists a nonnegative function
$M_{\alpha, \beta, K}(\cdot, \cdot)\in L^p_{\mathcal{F}}(0, T)$ such
that $$\left|
\partial^\alpha_\xi\partial^\beta_x a(t, \omega, x, \xi)\right|
\leq M_{\alpha, \beta, K}(t, \omega)(1+|\xi|)^{\ell-|\alpha|},$$ for
a.e. $(t, \omega)\in (0, T)\times\Omega$ and any $(x, \xi)\in
K\times \dbR^n$. We write $a\in S^{\ell}_p(G\times \dbR^n)$ for
short.
\end{definition}
\begin{remark}\label{d2}
In Definition \ref{d1}, if $G=\dbR^n$ and $M_{\alpha, \beta,
K}(\cdot, \cdot)\equiv M_{\alpha, \beta}(\cdot, \cdot)$, which is
independent of the compact set $K$, we set $a\in S^{\ell}_p$.
\end{remark}

By Definition \ref{d1}, it is easy to show that for any $\ell_1,\
\ell_2\in \dbR$, if $\ell_1\leq \ell_2$, then $S^{\ell_1}_p(G\times
\dbR^n)\subseteq S^{\ell_2}_p(G\times \dbR^n)$. Therefore, we write
$$S^{\infty}_p(G\times \dbR^n)=\bigcup_{\ell\in\dbR}S^{\ell}_p(G\times \dbR^n),\ \
S^{-\infty}_p(G\times
\dbR^n)=\bigcap_{\ell\in\dbR}S^{\ell}_p(G\times \dbR^n).$$ Moreover,
if $a_1\in S^{\ell_1}_p(G\times \dbR^n)$ and $a_2\in
S^{\ell_2}_q(G\times \dbR^n)$, then for any multi-indices $\alpha$
and $\beta$,
$$
\partial^\alpha_\xi\partial^\beta_x a_1\in S^{\ell_1-|\alpha|}_p
(G\times \dbR^n),\  a_1+a_2\in S^{\max{\{\ell_1,
\ell_2\}}}_{\min{\{p, q\}}}(G\times \dbR^n),\ a_1 a_2\in
S^{\ell_1+\ell_2}_{q^*}(G\times \dbR^n), $$ here and hereafter $q^*$
denotes a constant defined as follows: $q^*=\frac{pq}{p+q},\mbox{
for }p, q\geq 1, pq\geq p+q;\ q^*=p,\mbox{ for }p\geq 1, q=\infty;\
q^*=q,\mbox{ for }q\geq 1, p=\infty;\ q^*=\infty,\mbox{ for
}p=q=\infty$.

\medskip

Now, we introduce the notion of SPDOs.

\begin{definition}\label{d3}
A linear operator $A$ is called a SPDO of order $(\ell, p)$ if $a\in
S^\ell_p(G\times \dbR^n)$ and for any $u\in L^q_{\mathcal{F}}(0, T;
\mathcal{D}(G))$,
\begin{eqnarray*}
&&(Au)(t, \omega, x) =(2\pi)^{-n}\displaystyle\int_{\dbR^n}e^{i
x\cdot\xi}a(t, \omega, x, \xi)\hat{u}(t, \omega,
\xi)d\xi,\end{eqnarray*} where $\hat{u}(t, \omega,
\xi)=\displaystyle\int_G e^{-i x\cdot\xi}u(t, \omega, x)dx$. We
write $A\in \mathcal{L}^\ell_p(G)$. Moreover, if $a\in S^\ell_p$, we
set $A\in \mathcal{L}^\ell_p$.
\end{definition}
Notice that for a.e. $(t, \omega)\in(0, T)\times\Omega$, the SPDO
$A$ in Definition \ref{d3} is indeed a usual pseudo-differential
operator of order $\ell$ in the deterministic case. For simplicity
of notation, we write it simply $A$ when no confusion can arise in
this paper.

\medskip

In the following, we will investigate the domain and range of SPDOs.
Before that, we present a useful lemma. Write
$|v|_{\alpha}=\sup\limits_{x\in\dbR^n,\atop |\beta|\leq
2(|\alpha|+\ell+1+n)}
\left|(1+|x|)^{1+n}\partial_{x}^{\beta}v(x)\right|$, for any $v\in
C^\infty_0(G)$. Then, we have the following result.
\begin{lemma}\label{109}
$|\cdot|_\alpha$ is a generating semi-norm on $\mathcal{D}(G)$.
\end{lemma}
\noindent {\bf Proof. }By Lemma 1.34 in \cite{ru}, the set $U^*=\{\
v\in C_0^\infty(G);\ |v|_\alpha<1\}$ is balanced, convex, absorbing,
and $|\cdot|_\alpha$ is the Minkowski functional of $U^*$.
Therefore, it remains to prove that $U^*$ is the member of a
0-neighborhood base for the inductive topology on $\mathcal{D}(G)$.
By Lemma \ref{103}, it suffices to show that for any compact set
$K\subseteq G$, $V^*=\{\ v\in \mathcal{D}_K;\ |v|_\alpha<1\}$ is a
0-neighborhood in $\mathcal{D}_K$. By the definition of
$|\cdot|_\alpha$, for any $v\in \mathcal{D}_K$,
$$
\sup\limits_{x\in\dbR^n,\atop |\beta|\leq 2(|\alpha|+\ell+1+n)}
\left|(1+|x|)^{1+n}\partial_{x}^{\beta}v(x)\right|\leq N^*
\sup\limits_{x\in K,\atop |\beta|\leq 2(|\alpha|+\ell+1+n)}
|\partial_{x}^{\beta}v(x)|,
$$
where $N^*=1+[(1+\sup\limits_{x\in K}|x|)^{1+n}]$ and $[\ell]$
denotes the integral part of a real number $\ell$. This implies that
$$
\left\{\ v\in\mathcal{D}_K\ \left|\ |v|_{K, 2(|\alpha|+\ell+1+n)}=
\sup\limits_{x\in K,\atop |\beta|\leq
2(|\alpha|+\ell+1+n)}|\partial^\beta_x
v(x)|<\frac{1}{N^*}\right\}\subseteq V^*\right..
$$
Therefore, $V^*$ is a 0-neighborhood in $\mathcal{D}_K$. This
finishes the proof.\endpf

\medskip

Based on Lemma \ref{109}, we get the following result.
\begin{theorem}\label{t1}
Suppose that $A$ is a SPDO determined by a symbol $a$.

\noindent$(1)$ If $a\in S^\ell_p(G\times \dbR^n)$, $A:
L^q_{\mathcal{F}}(0, T; \mathcal{D}(G))\rightarrow
L^{q^*}_{\mathcal{F}}(0, T;
\mathcal{E}(G))$ is continuous;\\[1mm]
\noindent$(2)$ If $a\in S^\ell_p$, $A: L^q_{\mathcal{F}}(0, T;
\mathcal{S})\rightarrow L^{q^*}_{\mathcal{F}}(0, T; \mathcal{S})$ is
continuous.
\end{theorem}

\noindent{\bf Proof. }\ For any $u\in L^q_{\mathcal{F}}(0, T;
\mathcal{D}(G))$, multi-index $\alpha$ and any compact set
$K\subseteq G$, by Definition \ref{d1}, we have that
\begin{eqnarray*}
&&|\partial^\alpha_x[e^{i x\cdot\xi}a(t, \omega, x, \xi) \hat{u}(t,
\omega, \xi)]|
=\left|\sum\limits_{\alpha_1+\alpha_2=\alpha}\displaystyle\frac{\alpha!}
{\alpha_1!\alpha_2!}i^{|\alpha_1|}\xi^{\alpha_1} e^{i
x\cdot\xi}\partial^{\alpha_2}_x
a(t, \omega, x, \xi)\hat{u}(t, \omega, \xi)\right|\\[1mm]
&&=\left|\sum\limits_{\alpha_1+\alpha_2=\alpha}\displaystyle\frac{\alpha!}
{\alpha_1!\alpha_2!}i^{|\alpha_1|}\xi^{\alpha_1} e^{i
x\cdot\xi}\partial^{\alpha_2}_x a(t, \omega, x,
\xi)\displaystyle\int_G
e^{-ix\cdot\xi}u(t, \omega, x)dx\right|\\[1mm]
&&\leq  C(n, \alpha)M_{\alpha, K}(t,
\omega)(1+|\xi|)^{|\alpha|+\ell}\left|\displaystyle\int_G
e^{-ix\cdot\xi}u(t, \omega, x)dx\right|\\[1mm]
&&\leq  C(n, \alpha, \ell)M_{\alpha, K}(t,
\omega)(1+|\xi|)^{-n-1}\left|\displaystyle\int_G
(1+|\xi|^2)^{|\alpha|+\ell+1+n}e^{-ix\cdot\xi}u(t, \omega, x)dx\right|\\[1mm]
&&=C(n, \alpha, \ell)M_{\alpha, K}(t,
\omega)(1+|\xi|)^{-n-1}\left|\displaystyle\int_G
(1-\Delta_x)^{|\alpha|+\ell+1+n}e^{-ix\cdot\xi}u(t, \omega, x)dx\right|\\[1mm]
&&\leq C(n, \alpha, \ell)M_{\alpha, K}(t,
\omega)(1+|\xi|)^{-n-1}\sup\limits_{x\in \dbR^n}\left|(1+|x|)^{1+n}
(1-\Delta_x)^{|\alpha|+\ell+1+n}u(t, \omega, x)\right|,
\end{eqnarray*}
for a.e. $(t, \omega)\in(0, T)\times\Omega$ and any $(x, \xi)\in
K\times \dbR^n$, where $\Delta_x$ denotes the Laplacian operator
with respect to $x$. It follows that for a.e. $(t, \omega)\in(0,
T)\times\Omega$, $(Au)(t, \omega, \cdot)\in C^\infty(G)$. Moreover,
$$
\sup\limits_{x\in K}\left|(\partial^\alpha_x (Au))(t, \omega,
x)\right|\leq C(n, \alpha, \ell)M_{\alpha, K}(t, \omega) |u(t,
\omega, \cdot)|_\alpha.
$$
Hence,   $$\left|\sup\limits_{x\in K}(\partial^\alpha_x (Au))(\cdot,
\cdot, x)\right|_{L^{q^*}_{\mathcal{F}}(0, T)}\leq C(n, \alpha,
\ell)|M_{\alpha, K}(\cdot, \cdot)|_{L^p_{\mathcal{F}}(0,
T)}\left||u(\cdot, \cdot,
\cdot)|_\alpha\right|_{L^q_{\mathcal{F}}(0, T)}.
$$
By Lemma \ref{109},  this implies the desired continuity of $A$.
Also, notice that both the limit of a sequence of measurable
functions and the sum of a finite number of measurable functions are
measurable. Therefore, for any $x\in G$, $(Au)(\cdot, \cdot, x)$ is
$\{\mathcal{F}_t\}$-adapted.

(2) in Theorem \ref{t1} can be derived in the same way. The proof is
completed.\endpf

\begin{remark}\label{111}
It is easy to check that if $a\in S^\ell_p(G\times \dbR^n)$, then
the associated SPDO $A$ satisfies the following conditions:

\noindent $(1)$ for a.e. $(t, \omega)\in(0, T)\times \Omega$, $A:
\mathcal{D}(G)\rightarrow \mathcal{E}(G)$ is continuous;\\[1mm]
\noindent $(2)$ $A: L^q_{\mathcal{F}}(0, T;
\mathcal{D}_G)\rightarrow L^{q^*}_{\mathcal{F}}(0, T;
\mathcal{E}(G))$ is continuous.
\end{remark}
\begin{remark}
We can endow $S^{\ell}_p(G\times\dbR^n)$ and
$\mathcal{L}^{\ell}_p(G)$ with suitable topological structures,
respectively, such that any SPDO of order $(\ell, p)$ has structural
stability with respect to its amplitude. Indeed, let
$\{K_j\}_{j\in\dbN}$ be a sequence of compact sets satisfying
$K_0\subseteq K_1\subseteq\cdots$ and
$\displaystyle\bigcup_{j\in\dbN}K_j=G$. For any $k, N, j\in\dbN$,
set
$$
|a|_{k, N, j, \ell, p}=\left|\sup\limits_{x\in K_j, \xi\in\dbR^n,
\atop |\alpha|\leq k, |\beta|\leq N
}\displaystyle\frac{|\partial^\alpha_\xi\partial^\beta_x a(\cdot,
\cdot, x,
\xi)|}{(1+|\xi|)^{\ell-|\alpha|}}\right|_{L^p_{\mathcal{F}}(0, T)},
\quad\forall\ a\in S^{\ell}_p(G\times\dbR^n).
$$
Then $S^{\ell}_p(G\times\dbR^n)$ is a Fr\'{e}chet space, generated
by a sequence of semi-norms $\{|\cdot|_{k, N, j, \ell, p}\}_{k, N,
j\in\dbN}$. On the other hand, by the proof of Theorem \ref{t1}, we
see that for any SPDO $A$ determined by a symbol $a$, any
multi-index $\alpha$, $j\in\dbN$ and $u\in L^q_{\mathcal{F}}(0, T;
\mathcal{D}(G))$,
$$|Au|_{q^*,
 j, |\alpha|, 1}\leq C(n, \alpha, \ell)|a|_{0, |\alpha|, j, \ell,
p}\left||u(\cdot, \cdot, \cdot)|_\alpha\right|_{L^q_{\mathcal{F}}(0,
T)},$$ where $|\cdot|_{q^*,
 j, |\alpha|, 1}$ and $\left||\cdot|_\alpha\right|_{L^q_{\mathcal{F}}(0,
T)}$ are one of the generating semi-norms of
$L^{q^*}_{\mathcal{F}}(0, T; \mathcal{E}(G))$ and
$L^q_{\mathcal{F}}(0, T; \mathcal{D}(G))$, respectively. Therefore,
if we write
$$
|A|_{\alpha, j, \ell, p}=\sup\limits_{u\in L^q_{\mathcal{F}}(0, T;
\mathcal{D}(G)), \atop u\neq 0}\displaystyle\frac{|Au|_{q^*, j,
|\alpha|, 1}}{\left||u(\cdot, \cdot,
\cdot)|_\alpha\right|_{L^q_{\mathcal{F}}(0, T)}}, \quad\forall\ A\in
\mathcal{L}^\ell_p(G),
$$
then $\mathcal{L}^{\ell}_p(G)$ is a Fr\'{e}chet space, generated by
a sequence of semi-norms $\{|\cdot|_{\alpha, j, \ell,
p}\}_{j\in\dbN, |\alpha|\geq 0}$. Moreover, $|A|_{\alpha, j, \ell,
p}\leq C(n, \alpha, \ell)|a|_{0, |\alpha|, j, \ell, p}$. This
implies that the mapping $S^{\ell}_p(G\times\dbR^n)\rightarrow
\mathcal{L}^{\ell}_p(G),\quad a\mapsto A_a$ is continuous, where
$A_a$ denotes the SPDO determined by a symbol $a$ of order $(\ell,
p)$.
\end{remark}

\subsection{Amplitude and stochastic pseudo-differential operator}

In this subsection, we introduce  a class of SPDOs, which is
apparently more general than the class of SPDOs defined in the last
subsection. However, we shall show that these two classes in fact
coincide under some circumstances later. First of all, we give the
notion of amplitudes.
\begin{definition}\label{d21}
A complex-valued function $a$ is called an amplitude of order
$(\ell, p)$ if $a$ satisfies the following conditions:\\[2mm]
\noindent $(1)$ $a(t, \omega, \cdot, \cdot, \cdot)\in
C^\infty(G\times G\times \dbR^n)$,
a.e. $(t, \omega)\in (0, T)\times \Omega$;\\[1mm]
\noindent $(2)$ $a(\cdot, \cdot, x, y, \xi)$ is
$\{\mathcal{F}_t\}_{t\geq
0}$-adapted, $\forall\ (x, y, \xi)\in G\times G\times \dbR^n$;\\[1mm]
\noindent $(3)$ for any two multi-indices $\alpha$ and $\beta$, and
any compact set $K\subseteq G\times G$, there exists a nonnegative
function $M_{\alpha, \beta, K}(\cdot, \cdot)\in L^p_{\mathcal{F}}(0,
T)$ such that $$\left|
\partial^\alpha_\xi\partial^\beta_x a(t, \omega, x, y, \xi)\right|
\leq M_{\alpha, \beta, K}(t, \omega)(1+|\xi|)^{\ell-|\alpha|},$$ for
a.e. $(t, \omega)\in(0, T)\times\Omega$ and any $(x, y, \xi)\in
K\times \dbR^n$. We write $a\in S^\ell_p(G\times G\times\dbR^n)$ for
short.
\end{definition}

Based on Definition \ref{d21}, we have more freedom to construct
SPDOs.
\begin{definition}\label{d22}
The linear operator $A$ is called a SPDO of order $(\ell, p)$ if
$a\in S^\ell_p(G\times G\times \dbR^n)$ and for any $u\in
L^q_{\mathcal{F}}(0, T; \mathcal{D}(G))$,
$$(Au)(t, \omega, x)=(2\pi)^{-n}\displaystyle\int_{\dbR^n}\int_G e^{i
(x-y)\cdot \xi} a(t, \omega, x, y, \xi) u(t, \omega, y)dy d\xi.$$ We
write $A\in \mathcal{L}^\ell_p(G\times G)$.
\end{definition}

\begin{remark}\label{r21} For a.e. $(t, \omega)\in (0, T)\times\Omega$ and any $x\in  G$,
the integral in Definition \ref{d22} is understood as follows:
$$(Au)(t, \omega, x)=(2\pi)^{-n}\lim\limits_{\varepsilon\rightarrow 0}\displaystyle\int_{\dbR^n}
\displaystyle\int_{G}\chi(\varepsilon \xi)e^{i (x-y)\cdot\xi} a(t,
\omega, x, y, \xi)u(t, \omega, y)dy d\xi,$$ where $\chi\in
C^\infty_0(\dbR^n)$ and $\chi=1$ in a neighborhood of the origin.
Similar to   the deterministic case, it is easy to show that
Definition \ref{d22} is well posed.
\end{remark}

Also, by the same method as that used in the proof of Theorem
\ref{t1}, we obtain the following result for the SPDOs determined by
amplitudes.
\begin{theorem}\label{t21}
If $A\in \mathcal{L}^\ell_p(G\times G)$, then  the following
assertions hold:

\noindent $(1)$ $A:\ L^q_{\mathcal{F}}(0, T;
\mathcal{D}(G))\rightarrow L^{q^*}_{\mathcal{F}}(0, T;
\mathcal{E}(G))$ is continuous;

\noindent $(2)$ $A:\ L^q_{\mathcal{F}}(0, T;
\mathcal{D}_G)\rightarrow L^{q^*}_{\mathcal{F}}(0, T;
\mathcal{E}(G))$ is continuous.
\end{theorem}

In the following, we extend the domain of SPDOs to a space of
distributions. For this purpose,  we first introduce the definition
of  transpose operators.
\begin{definition}\label{d23}
Suppose that  $A$ is a SPDO of order $(\ell, p)$ and $a$ is its
amplitude. Then an operator $^tA$ is called  the transpose operator
of $A$ if for any $u\in L^q_{\mathcal{F}}(0, T; \mathcal{D}(G))$,
$$(^tAu)(t, \omega, x)=(2\pi)^{-n}
\displaystyle\int_{\dbR^n}\displaystyle\int_{G} e^{i
(x-y)\cdot\xi}a(t, \omega, y, x, -\xi)u(t, \omega, y)dy d\xi.$$
\end{definition}
\begin{remark}\label{126}
It is easy to check that for any $A\in \mathcal{L}^\ell_p(G\times
G)$, $^tA\in \mathcal{L}^\ell_p(G\times G)$. Therefore, Theorem
\ref{t21} holds for $^tA$.
\end{remark}

Denote by $(L^p_{\mathcal{F}}(0, T; \mathcal{D}(G)))'$ the dual
space of the locally convex space $L^p_{\mathcal{F}}(0, T;
\mathcal{D}(G))$; and by $(L^p_{\mathcal{F}}(0, T;
\mathcal{E}(G)))'$ the dual space of the locally convex space
$L^p_{\mathcal{F}}(0, T; \mathcal{E}(G))$. Let $q_*$ be a constant
defined as follows: $q_*=\frac{pq}{p-q}$, for $p, q\geq 1$, $p> q$;
$q_*=q$, for $q\geq 1$, $p=\infty$; $q_*=\infty$, for $p=q$. Next,
we present the following result on an extension of the domain of
SPDOs.
\begin{theorem}\label{t22}
Suppose that $A\in\mathcal{L}^\ell_p(G\times G)$. Then $A:
(L^q_{\mathcal{F}}(0, T; \mathcal{E}(G)))'\rightarrow
(L^{q_*}_{\mathcal{F}}(0, T; \mathcal{D}(G)))'$ is continuous.
\end{theorem}

\noindent{\bf Proof. }\ For any $u \in (L^q_{\mathcal{F}}(0, T;
\mathcal{E}(G)))'$ and  $v\in L^{q_*}_{\mathcal{F}}(0, T;
\mathcal{D}(G))$, define
\begin{equation}\label{21}\langle Au,
v\rangle _{(L^{q_*}_{\mathcal{F}}(0, T; \mathcal{D}(G)))',
L^{q_*}_{\mathcal{F}}(0, T; \mathcal{D}(G))}= \langle u,
^tAv\rangle_{(L^q_{\mathcal{F}}(0, T;
\mathcal{E}(G)))',L^q_{\mathcal{F}}(0, T;
\mathcal{E}(G))}.\end{equation} Since $u\in (L^q_{\mathcal{F}}(0, T;
\mathcal{E}(G)))'$, there exists a semi-norm $|\cdot|_{q, j_0, k_0,
1}$ $(j_0, k_0\in\dbN)$ defined on $L^q_{\mathcal{F}}(0, T;
\mathcal{E}(G))$ such that
$$
|\langle u, ^tAv\rangle_{(L^q_{\mathcal{F}}(0, T;
\mathcal{E}(G)))',L^q_{\mathcal{F}}(0, T; \mathcal{E}(G))}|\leq
C(u)\ |^tA v|_{q, j_0, k_0, 1}.
$$
By the above result and Remark \ref{126}, we can find a semi-norm
$|\cdot|_{q_*, \alpha_0}$  defined on $L^{q_*}_{\mathcal{F}}(0, T;
\mathcal{D}(G))$ for a multi-index $\alpha_0$, such that
$$
|\langle u, ^tAv\rangle_{(L^q_{\mathcal{F}}(0, T;
\mathcal{E}(G)))',L^q_{\mathcal{F}}(0, T; \mathcal{E}(G))}|\leq C(u,
a)\ |v|_{q_*, \alpha_0},
$$
which implies that $Au \in (L^{q_*}_{\mathcal{F}}(0, T;
\mathcal{D}(G)))'$. Moreover,  the continuity of $A$ is clear from
(\ref{21}). The proof of Theorem \ref{t22} is completed.
\endpf
\begin{remark}
It is regrettable that we fail to give a characterization of
$(L^q_{\mathcal{F}}(0, T;
 \mathcal{E}(G)))'$ clearly. Here we only present two classes of function spaces, which are contained
 in $(L^q_{\mathcal{F}}(0, T;
 \mathcal{E}(G)))'$. Let $\mathcal{E}'(G)$ and $\mathcal{D}'(G)$
 denote the dual spaces of $\mathcal{E}(G)$
 and $\mathcal{D}(G)$, respectively.
 Then we see that $\mathcal{E}'(G)\subseteq(L^q_{\mathcal{F}}(0, T;
 \mathcal{E}(G)))'$. Indeed, for any $u\in \mathcal{E}'(G)$
 and $v\in L^q_{\mathcal{F}}(0, T;
 \mathcal{E}(G))$,
 define
$$
\langle u, v\rangle_{(L^q_{\mathcal{F}}(0, T;
 \mathcal{E}(G)))', L^q_{\mathcal{F}}(0, T;
 \mathcal{E}(G))}=\dbE\displaystyle\int^T_0
\langle u(\cdot), v(t, \omega,
\cdot)\rangle_{\mathcal{E}'(G),\mathcal{E}(G)}dt.$$ Then there exist
two nonnegative integers $j_1$ and $k_1$, such that
\begin{eqnarray*}
&&|\langle u, v\rangle_{(L^q_{\mathcal{F}}(0, T;
 \mathcal{E}(G)))', L^q_{\mathcal{F}}(0, T;
 \mathcal{E}(G))}|\leq
 \dbE\int^T_0 C(u) |v(t, \omega, \cdot)|_{j_1, k_1,
1}dt\\[2mm]
&&\leq C(u, T, q)\left| |v(\cdot, \cdot, \cdot)|_{j_1, k_1,
1}\right|_{L^q_{\mathcal{F}}(0, T)}=C(u, T, q)|v|_{q, j_1, k_1, 1}.
\end{eqnarray*}
This implies that $u\in(L^q_{\mathcal{F}}(0, T;
 \mathcal{E}(G)))'$.

On the other hand, we define a constant $q'$ as follows:
$q'=\frac{q}{q-1}$ for $q>1$; $q'=1$ for $q=\infty$; $q'=\infty$ for
$q=1$. Then, for any given compact set $K$, it is also easy to prove
that
$$\{\ u\in L^{q'}_{\mathcal{F}}(0, T; L^1(G))\ |\ \supp u(t, \omega,
\cdot)\subseteq K, \mbox{ for a.e. }(t, \omega)\in (0,
T)\times\Omega\}\subseteq (L^{q}_{\mathcal{F}}(0, T;
 \mathcal{E}(G)))'.$$

\end{remark}
\begin{remark}
For any $A\in \mathcal{L}^\ell_p(G\times G)$, it is easy to show
that for a.e. $(t, \omega)\in(0, T)\times \Omega$, $A:
\mathcal{E}'(G)\rightarrow \mathcal{D}'(G)$ is continuous.
 Moreover, if we let
\begin{eqnarray*}
&&\mathcal{D}_{s}(G)=\left\{u\ |\ \mbox{for a.e. }(t, \omega)\in (0,
T)\times\Omega,\ u(t, \omega, \cdot)\in
\mathcal{D}(G);\mbox{ and for any }x\in G, \right.\\[1mm]
&&\quad\quad\quad\quad\quad\quad\left. u(\cdot, \cdot, x)\mbox{ is
}\{\mathcal{F}_t\}_{t\geq 0}
\mbox{-adapted}\right\};\\[1mm]
&&\mathcal{E}_{s}(G)=\left\{u\ |\ \mbox{for a.e. }(t, \omega)\in (0,
T)\times\Omega,\ u(t, \omega, \cdot)\in
\mathcal{E}(G);\mbox{ and for any }x\in G, \right.\\[1mm]
&&\quad\quad\quad\quad\quad\quad\left. u(\cdot, \cdot, x)\mbox{ is
}\{\mathcal{F}_t\}_{t\geq 0} \mbox{-adapted}\right\};\\[1mm]
&&\mathcal{D}'_{s}(G)=\left\{u\ |\ \mbox{for a.e. }(t, \omega)\in
(0, T)\times\Omega,\ u(t, \omega, \cdot)\in \mathcal{D}^\prime(G);\
\mbox{and for any }
v\in \mathcal{D}_{s}(G), \right.\\[1mm]
&&\quad\quad\quad\quad\quad\quad\quad\left. \langle u(\cdot, \cdot,
\cdot), v(\cdot, \cdot, \cdot) \rangle_{\mathcal{D}(G),
\mathcal{D}^\prime(G)} \mbox{ is }\{\mathcal{F}_t\}_{t\geq 0}
\mbox{-adapted}\right\};\\[2mm]
&&\mathcal{E}'_{s}(G)=\left\{u\ |\ \mbox{for a.e. }(t, \omega)\in
(0, T)\times\Omega,\ u(t, \omega, \cdot)\in
\mathcal{E}^\prime(G);\ \mbox{and for any }v\in \mathcal{E}_{s}(G), \right.\\[1mm]
&&\quad\quad\quad\quad\quad\quad\quad\left. \langle u(\cdot, \cdot,
\cdot), v(\cdot, \cdot, \cdot) \rangle_{\mathcal{E}(G),
\mathcal{E}^\prime(G)} \mbox{ is }\{\mathcal{F}_t\}_{t\geq 0}
\mbox{-adapted}\right\},
\end{eqnarray*}
then for any $u\in \mathcal{E}'_{s}(G)$, $Au\in\mathcal{D}'_s(G)$.
\end{remark}
\subsection{Kernel and pseudo-local property}

In this subsection, we introduce the definition of kernels, and then
prove the pseudo-local property of SPDOs (see Theorem \ref{t32}).
First of all, we give the notion of kernels.

\begin{definition}\label{d31}
Suppose that $a\in S^\ell_p(G\times G\times \dbR^n)$ and $A$ is the
associated SPDO. Then $K_A\in \mathcal{D}'_{s}(G\times G)$ is called
a kernel of $A$ if for a.e. $(t, \omega)\in (0, T)\times \Omega$,
$$\langle K_A(t, \omega, \cdot, \cdot), \ v\rangle_{\mathcal{D}^\prime(G\times
G), \mathcal{D}(G\times G)}
=(2\pi)^{-n}\displaystyle\int_{\dbR^n}\displaystyle\int_{G}
\displaystyle\int_{G} e^{i (x-y)\cdot\xi}a(t, \omega, x, y,
\xi)v(x,y)dxdyd\xi,$$ for any $v\in C^\infty_0(G\times G)$.
\end{definition}
Notice that for a.e. $(t, \omega)\in(0, T)\times\Omega$, the kernel
in Definition \ref{d31} is indeed the one in the deterministic case.
Therefore, the following assertions hold:
\begin{proposition}\label{t31}
Suppose that $A$ is a SPDO and $a$ is its amplitude. Then,\\[1mm]
$(1)$ if $a\in S^\ell_p(G\times G\times \dbR^n)$, $K_A$ is
$C^\infty$ off the diagonal in $G\times G$ for a.e. $(t, \omega)\in
(0, T)\times \Omega$. Moreover, for any $(x, y)\in G\times G$ with
$x\neq y$, $K_A(\cdot, \cdot, x, y)$
 is $\{\mathcal{F}_t\}_{t\geq 0}$-adapted, and for any compact set $K\subseteq (G\times G)\setminus \{x=y\}$,
 $\sup\limits_{(x, y)\in K}|K_A(\cdot, \cdot, x, y)|\in L^p_{\mathcal{F}}(0, T)$;\\[1mm]
$(2)$ if $a\in S^{-\infty}_p(G\times G\times \dbR^n)$, $K_A\in
C^\infty(G\times G)$ for a.e. $(t, \omega)\in (0, T)\times \Omega$.
Moreover, for any $(x, y)\in G\times G$, $K_A(\cdot, \cdot, x, y)$
 is $\{\mathcal{F}_t\}_{t\geq 0}$-adapted, and for any compact set $K\subseteq G\times G$,
 $\sup\limits_{(x, y)\in K}|K_A(\cdot, \cdot, x, y)|\in L^p_{\mathcal{F}}(0, T)$.
\end{proposition}
{\bf Sketch of the proof. } First, for any $a\in S^\ell_p(G\times
G\times \dbR^n)$,  it is easy to show that for a.e. $(t, \omega)\in
(0, T)\times \Omega$ and any $(x, y)\in G\times G $ with $x\neq y$,
the integral $\displaystyle\int_{\dbR^n} e^{i (x-y)\cdot\xi}a(t,
\omega, x, y, \xi)d\xi$ can be understood in the following two
senses equivalently:
\begin{eqnarray*}
&&\displaystyle\int_{\dbR^n} e^{i (x-y)\cdot\xi}a(t, \omega, x, y,
\xi)d\xi=\lim\limits_{\varepsilon\rightarrow
0}\displaystyle\int_{\dbR^n}\chi(\varepsilon\xi)e^{i
(x-y)\cdot\xi}a(t, \omega, x, y, \xi)d\xi,\\
&&=\displaystyle\int_{\dbR^n} e^{i
(x-y)\cdot\xi}(-1)^k|x-y|^{-2k}\Delta^k_\xi a(t, \omega, x, y,
\xi)d\xi,
\end{eqnarray*}
where $\chi\in C^\infty_0(\dbR^n)$, $\chi=1$ in a neighborhood of
the origin and $k\in \dbN$ with $\ell-2k<-n$.

Also, by the meaning of the above integral, for a.e. $(t, \omega)\in
(0, T)\times \Omega$ and any compact set $K\subseteq(G\times
G)\setminus\{x=y\}$, we have that $\displaystyle\int_{\dbR^n} e^{i
(x-y)\cdot\xi}a(t, \omega, x, y, \xi)d\xi\in C^\infty(K)$.

Moreover, by the Lebesgue dominated convergence theorem, it is easy
to check that $K_A(t, \omega, x, y)$ $
=(2\pi)^{-n}\displaystyle\int_{\dbR^n} e^{i (x-y)\cdot\xi}a(t,
\omega, x, y, \xi)d\xi$, for a.e. $(t, \omega)\in (0, T)\times
\Omega$ and any $(x, y)\in G\times G $ with $x\neq y$. Therefore,
 $K_A(t, \omega, \cdot, \cdot)\in C^\infty((G\times G)\backslash\{x=y\})$
 for a.e. $(t, \omega)\in(0, T)\times\Omega$.

Furthermore, by the properties of measurable functions, we see that
$K_A(\cdot, \cdot, x, y)$
 is $\{\mathcal{F}_t\}_{t\geq 0}$-adapted,
 and by the definition of amplitudes, $\sup\limits_{(x, y)\in K}|K_A(\cdot, \cdot, x, y)|\in L^p_{\mathcal{F}}(0, T)$ for any
 compact set $K\subseteq (G\times G)\setminus\{x=y\}$.

Similar to the above procedure, we can also get the desired result
(2), if $a\in S^{-\infty}_p(G\times G\times \dbR^n)$.
 \endpf
\begin{remark}\label{33}
 For a SPDO $A$, we call  $A$ a smoothing operator of order $p$
 if $K_A$ satisfies the following conditions:\\[1mm]
\noindent $(1)$ for a.e. $(t, \omega)\in (0, T)\times \Omega$,
$K_A(t,
\omega, \cdot, \cdot)\in C^\infty(G\times G)$;\\[1mm]
\noindent $(2)$ for any $(x, y)\in G\times G$, $K_A(\cdot, \cdot, x,
y)$ is $\{\mathcal{F}_t\}_{t\geq 0}$-adapted;\\[1mm]
\noindent $(3)$ for any compact set $K\subseteq G\times G$,
$\sup\limits_{(x, y)\in K}|K_A(\cdot, \cdot, x, y)|\in
L^p_{\mathcal{F}}(0, T)$.

It is easy to check that if a function $\widetilde{K}$ satisfies the
above conditions $(1)$-$(3)$, then an operator $A$ defined as
follows: $(Au)(t, \omega, x)=\displaystyle\int_{G} \widetilde{K}(t,
\omega, x, y)u(t, \omega, y)dy$ is a SPDO and its amplitude $a\in
S^{-\infty}_p(G\times G\times\dbR^n)$. Combining Proposition
\ref{t31} with the above fact, we obtain that $a\in
S^{-\infty}_p(G\times G\times\dbR^n)$ if and only if $A$ is a
smoothing operator of order $p$.
\end{remark}
\medskip

Next, we recall the notion of singular supports. For a distribution
$u$, the singular support of $u$ is the complement of the open set
on which $u$ is smooth and we write it sing supp $u$.

In the following, we present the pseudo-local property for SPDOs.
\begin{theorem}\label{t32}
Suppose that $A$ is a SPDO. Then for a.e. $(t, \omega)\in (0,
T)\times \Omega$ and any $u\in \mathcal{E}'_{s}(G)$,
$$\mbox{sing supp }(Au)(t, \omega, \cdot)\subseteq\mbox{sing supp }u(t, \omega, \cdot).$$
\end{theorem}
The proof of Theorem \ref{t32} follows by the similar method as that
used in the deterministic case. Therefore, we omit it here.

\subsection{Uniformly properly supported stochastic pseudo-differential
operator}

In order to present the composition of two SPDOs, in this
subsection, we introduce uniformly properly supported SPDOs with
respect to $(t, \omega)$. First of all, we recall the notion of
proper sets. A set $E\subseteq G\times G$ is called a proper set if
$E$ has compact intersection with $K\times G$ and with $G\times K$,
for any compact set $K\subseteq G$. Also, we give some relevant
definitions.
\begin{definition}\label{d41}
A function $\widetilde{K}\in \mathcal{D}'_{s}(G\times G)$ is said to
be uniformly properly supported with respect to $(t, \omega)$ if
there exists a proper set $E$ such that $\supp\widetilde{K}(t,
\omega, \cdot, \cdot)\subseteq E$ for a.e. $(t, \omega)\in (0,
T)\times \Omega$.
\end{definition}
\begin{definition}\label{d42}
We call $A$ a uniformly properly supported SPDO with respect to $(t,
\omega)$ if $A$ is a SPDO and its kernel $K_A$ is uniformly properly
supported with respect to $(t, \omega)$.
\end{definition}
\begin{definition}\label{d43}
Suppose that $a\in S^\ell_p(G\times G\times \dbR^n)$. $a$ is said to
have uniformly proper support with respect to $(t, \omega, \xi)$ if
there exists a proper set $E\subseteq G\times G$ such that $\supp
a(t, \omega, \cdot, \cdot, \xi)\subseteq E$ for a.e. $(t, \omega)\in
(0, T)\times \Omega$ and any $\xi\in \dbR^n$.
\end{definition}

Next, we give a characterization of  amplitudes for uniformly
properly supported SPDOs with respect to $(t, \omega)$.
\begin{proposition}\label{t41}
For a SPDO $A$, if its amplitude $a$ has uniformly proper support
with respect to $(t, \omega, \xi)$, then $A$ is a uniformly properly
supported SPDO with respect to $(t, \omega)$; Conversely, if $A$ is
a uniformly properly supported SPDO with respect to $(t, \omega)$,
then its amplitude $a$ can be replaced by another one, which has
uniformly proper support with respect to $(t, \omega, \xi)$.
\end{proposition}
\noindent{\bf Sketch of the proof. } First, it is easy to show the
fact: if $E$ is a proper set, then there exists a function $\psi\in
C^\infty(G\times G)$ such that $\psi=1$ in a neighborhood of $E$ and
$\supp\psi$ is proper.

Next, if $A\in\mathcal{L}^\ell_p(G\times G)$ is a uniformly properly
supported SPDO with respect to $(t, \omega)$, by the above fact,
there exist a function $\psi_1\in C^\infty(G\times G)$ and a proper
set $E_1$, such that $\psi_1=1$ in a neighborhood of $E_1$,
$\supp\psi_1$ is proper and $\supp K_A(t, \omega, \cdot,
\cdot)\subseteq E_1$ for a.e. $(t, \omega)\in(0, T)\times\Omega$.
This implies that for any $u, v\in C^\infty_0(G)$,
\begin{eqnarray*}
&&\langle(Au)(t, \omega, \cdot), v\rangle_{\mathcal{D}'(G),
\mathcal{D}(G)}=\langle K_A(t, \omega, \cdot, \cdot),
uv\rangle_{\mathcal{D}'(G\times G),
\mathcal{D}(G\times G)}\\[2mm]
&&=\langle\psi_1 K_A(t, \omega, \cdot, \cdot),
uv\rangle_{\mathcal{D}'(G\times G), \mathcal{D}(G\times G)}+
\langle(1-\psi_1) K_A(t, \omega, \cdot, \cdot),
uv\rangle_{\mathcal{D}'(G\times G), \mathcal{D}(G\times G)}\\[2mm]
&&=(2\pi)^{-n}\displaystyle\int_{\dbR^n}\displaystyle\int_{\dbR^n}
\displaystyle\int_{\dbR^n}e^{i(x-y)\cdot\xi}a(t, \omega, x, y,
\xi)\psi_1(x, y)u(y)v(x)dxdyd\xi.
\end{eqnarray*}
Set $a_*(t, \omega, x, y, \xi)=a(t, \omega, x, y, \xi)\psi_1(x, y)$
and denote by $A_*$ its associated SPDO. Then $a_*\in
S^\ell_p(G\times G\times\dbR^n)$ and  $A=A_*$. Moreover, since
$\supp a_*(t, \omega, \cdot, \cdot, \xi)\subseteq\supp \psi_1$ for
a.e. $(t, \omega)\in(0, T)\times\Omega$ and any $\xi\in\dbR^n$, the
amplitude $a_*$ has uniformly proper support with respect to $(t,
\omega, \xi)$.

On the other hand, if $a$ has uniformly proper support  with respect
to $(t, \omega, \xi)$, then there exists a proper set $E_2$, such
that  $\supp a(t, \omega, \cdot, \cdot, \xi)\subseteq E_2$ for a.e.
$(t, \omega)\in(0, T)\times\Omega$ and any $\xi\in\dbR^n$. This
leads to that for any $\varphi\in C^\infty_0(\overline{E_2}^c)$,
$$
\langle K_A(t, \omega, \cdot, \cdot),
\varphi\rangle_{\mathcal{D}'(G\times G), \mathcal{D}(G\times
G)}=(2\pi)^{-n}\displaystyle\int_{\dbR^n}
\displaystyle\int_{\dbR^n}\displaystyle\int_{\dbR^n} e^{i(x-y)\cdot
\xi}a(t, \omega, x, y, \xi)\varphi(x, y)dxdyd\xi=0.
$$
Therefore, $\supp K_A(t, \omega, \cdot, \cdot)\subseteq
\overline{E_2}$ for a.e. $(t, \omega)\in(0, T)\times\Omega$, which
implies that $A$ is a uniformly properly supported SPDO with respect
to $(t, \omega)$.
\endpf

\medskip

Based on Proposition \ref{t41}, we get that a uniformly properly
supported SPDO with respect to $(t, \omega)$ actually has better
properties than the usual operators.
\begin{theorem}\label{t42}
Suppose that $A$ is a uniformly properly supported SPDO of order
$(\ell, p)$ with respect to $(t, \omega)$. Then

\noindent $(1)$ $A: L^q_{\mathcal{F}}(0, T;
\mathcal{D}_G)\rightarrow
 L^{q^*}_{\mathcal{F}}(0, T; \mathcal{D}_G)$ is continuous;

\noindent $(2)$ the domain of $A$ can  be extended to be
$L^q_{\mathcal{F}}(0, T; \mathcal{E}(G))$. Moreover, $A:
L^q_{\mathcal{F}}(0, T; \mathcal{E}(G))\rightarrow
L^{q^*}_{\mathcal{F}}(0, T; \mathcal{E}(G))$ is continuous.
\end{theorem}
\noindent {\bf Proof. } Denote by $a$ the amplitude of $A$. Since
$A$ is a uniformly properly supported SPDO with respect to $(t,
\omega)$, by Proposition \ref{t41}, without loss of generality, we
suppose that there is a proper set $E_3$, such that for a.e. $(t,
\omega)\in (0, T)\times\Omega$ and any $\xi\in\dbR^n$, $\supp a(t,
\omega, \cdot, \cdot, \xi)\subseteq E_3$. Also, for any $u\in
L^q_{\mathcal{F}}(0, T; \mathcal{D}_G)$, there exists a compact set
$K^0$ such that $\supp u(t, \omega, \cdot)\subseteq K^0$ for a.e.
$(t, \omega)\in(0, T)\times\Omega$. Write $K^1=\{\ x\in G\ |\
\mbox{there exists a }y \in K^0 \mbox{ such that }(x, y)\in E_3\}$.
Then $K^1$ is a compact set, and if $x$ does not belong to $K^1$,
$$(Au)(t, \omega, x)=(2\pi)^{-n}\displaystyle\int_{\dbR^n}\int_{G}
e^{i(x-y)\cdot\xi}a(t, \omega, x, y, \xi)u(t, \omega, y)dyd\xi=0.$$
This means that $\supp Au(t, \omega, \cdot)\subseteq K^1$ for a.e.
$(t, \omega)\in(0, T)\times \Omega$. By Theorem  \ref{t21}, $Au\in
L^{q^*}_{\mathcal{F}}(0, T; \mathcal{E}(G))$. Therefore, $Au\in
L^{q^*}_{\mathcal{F}}(0, T; \mathcal{D}_G)$.

Moreover, if $\lim\limits_{j\rightarrow \infty}u_j=0$ in
$L^q_{\mathcal{F}}(0, T; \mathcal{D}_G)$, then by Proposition
\ref{104}, there exists a compact set $K^2\subseteq G$, such that
$\displaystyle\bigcup_{j\in\dbN}\supp u_j(t, \omega, \cdot)\subseteq
K^2$, for a.e. $(t, \omega)\in(0, T)\times\Omega$. Write $K^3=\{x\in
G\ |\mbox{ there exists a }y\in K^2,\mbox{ such that }(x, y)\in
E_3\}$. Then $K^3$ is a compact set and $\supp (Au_j)(t, \omega,
\cdot)\subseteq K^3$, for a.e. $(t, \omega)\in(0, T)\times\Omega$
and any $j\in\dbN$. Since $\lim\limits_{j\rightarrow\infty}Au_j=0 $
in $L^{q^*}_{\mathcal{F}}(0, T; \mathcal{E}(G))$, then by
Proposition \ref{104} again, $\lim\limits_{j\rightarrow
\infty}Au_j=0$ in $L^{q^*}_{\mathcal{F}}(0, T; \mathcal{D}_G)$.

On the other hand, for any open set $U\subseteq G$ with
$\overline{U}$ being compact, we write $K^4=\{y\in G\ |\ \mbox{there
exists an }x\in\overline{U}, \mbox{ such that }(x, y)\in E_3\}$, and
then $K^4$ is a compact set. Choose a function $\psi_2$, such that
$\psi_2\in C^\infty_0(G)$ and $\psi_2=1$ in $K^4$. We extend the
domain of a uniformly properly supported SPDO with respect to $(t,
\omega)$ as follows: $(Au)(t, \omega, x)=(A(\psi_2 u))(t, \omega,
x)$, for any $u\in L^q_{\mathcal{F}}(0, T; \mathcal{E}(G))$, a.e.
$(t, \omega)\in(0, T)\times\Omega$ and any $x\in U.$ It is easy to
check that the above definition of $A$ is well posed, namely, the
definition of $A$ is independent of $\psi_2$ and $U$.

Furthermore, for any compact set $K$, there exists an open set
$U_1\subseteq G$ such that $K\subseteq U_1$ and $\overline{U_1}$ is
compact. Write $K^5=\{ y\in G\ |\ \mbox{there exists an }x\in
\overline{U_1}\mbox{ such that }(x, y)\in E_3\}$, take a function
$\psi_3$ such that $\psi_3\in C_0^\infty(G)$ and $\psi_3=1$ in
$K^5$, and then set $K^6=\supp\psi_3$. Then for any $u\in
L^q_{\mathcal{F}}(0, T; \mathcal{E}(G))$ and multi-index $\alpha$,
\begin{eqnarray*}
&&|(\partial^\alpha_x(Au))(t, \omega, x)|\\[1mm]
&&=\left|(2\pi)^{-n}\int_{\dbR^n}\int_G
\partial_x^\alpha[e^{i(x-y)\cdot\xi}a(t, \omega, x, y, \xi)
u(t, \omega, y)\psi_3(y)]dyd\xi\right|\\[1mm]
&&=\left|(2\pi)^{-n}\int_{\dbR^n}\int_G\sum\limits_{\alpha_1+\alpha_2=\alpha}
\frac{\alpha!}{\alpha_1!\alpha_2!}i^{|\alpha_1|}
\xi^{\alpha_1}e^{i(x-y)\cdot\xi}\partial_x^{\alpha_2}a(t, \omega, x,
y, \xi)u(t, \omega, y)\psi_3(y)dyd\xi\right|\\[1mm]
&&=\left|(2\pi)^{-n}\int_{\dbR^n}\sum\limits_{\alpha_1+\alpha_2=\alpha}
\frac{\alpha!}{\alpha_1!\alpha_2!}i^{|\alpha_1|}(1+|\xi|^2)^{-|\alpha|-\ell-n-1}
\xi^{\alpha_1}\int_{G}
(1-\Delta_y)^{|\alpha|+\ell+n+1}e^{i(x-y)\cdot\xi}\right.\\[1mm]
&&\quad\cdot\partial_x^{\alpha_2}a(t, \omega, x,
y, \xi)u(t, \omega, y)\psi_3(y)dyd\xi\big|\\[1mm]
&&=\left|(2\pi)^{-n}\int_{\dbR^n}\sum\limits_{\alpha_1+\alpha_2=\alpha}
\frac{\alpha!}{\alpha_1!\alpha_2!}i^{|\alpha_1|}(1+|\xi|^2)^{-|\alpha|-\ell-n-1}
\xi^{\alpha_1}\int_{G}
e^{i(x-y)\cdot\xi}(1-\Delta_y)^{|\alpha|+\ell+n+1}\right.\\[1mm]
&&\left.\quad[\partial_x^{\alpha_2}a(t, \omega, x,
y, \xi)u(t, \omega, y)\psi_3(y)]dyd\xi\right|\\[1mm]
&&\leq C(n, \ell, \alpha, K)M_{n, \ell, \alpha, K}(t, \omega)
\sup\limits_{y\in K^6, \atop |\beta|\leq
2(|\alpha|+\ell+n+1)}|\partial^\beta_y u(t, \omega, y)|,
\end{eqnarray*}
for a function $M_{n, \ell, \alpha, K}(\cdot, \cdot)\in
L^p_{\mathcal{F}}(0, T)$, a.e. $(t, \omega)\in(0, T)\times \Omega$
and any $x\in K$. It follows that
\begin{eqnarray*}
&&\left|\sup\limits_{x\in K}|(\partial^\alpha_x(Au))(\cdot, \cdot,
x)|\right|_{L^{q^*}_{\mathcal{F}}(0, T)}\\[2mm]
&&\leq C(n, \ell, \alpha, K)|M_{n, \ell, \alpha, K}(\cdot,
\cdot)|_{L^p_{\mathcal{F}}(0, T)}\left|\sup\limits_{y\in K^6, \atop
|\beta|\leq 2(|\alpha|+\ell+n+1)}\left|\partial^\beta_y u(\cdot,
\cdot, y)\right|\right|_{L^q_{\mathcal{F}}(0, T)}.
\end{eqnarray*}
This implies that $A: L^q_{\mathcal{F}}(0, T;
\mathcal{E}(G))\rightarrow L^{q^*}_{\mathcal{F}}(0, T;
\mathcal{E}(G))$ is continuous. The proof is completed.
\endpf

\medskip

By Theorem \ref{t42}, for a.e. $(t, \omega)\in(0, T)\times\Omega$, a
finite number of uniformly properly supported SPDOs with respect to
$(t, \omega)$ can be composed. Moreover, for $\ell_1,
\ell_2\in\dbR$, suppose that $A\in \mathcal{L}^{\ell_1}_p(G)$ and
$B\in \mathcal{L}^{\ell_2}_q(G)$. If $pq\bar{q}\geq
pq+\bar{q}(p+q)$, then $B\circ A: L^{\bar{q}}_{\mathcal{F}}(0, T;
\mathcal{D}_G)\rightarrow L^{\tilde{q}}_{\mathcal{F}}(0, T;
\mathcal{D}_G)$ and $B\circ A: L^{\bar{q}}_{\mathcal{F}}(0, T;
\mathcal{E}(G))\rightarrow L^{\tilde{q}}_{\mathcal{F}}(0, T;
\mathcal{E}(G))$ are continuous, where
$\tilde{q}=\frac{pq\bar{q}}{pq+\bar{q}(p+q)}$. In addition, since a
kernel is smooth off the diagonal in $G\times G$, we may write a
SPDO as the sum of a uniformly properly supported operator and a
smoothing operator.
\begin{theorem}\label{t43}
If $A$ is a SPDO of order $(\ell, p)$, then $A=A^0+A^1$, where $A^0$
is a uniformly properly supported SPDO with respect to $(t, \omega)$
and $A^1$ is a smoothing operator.
\end{theorem}
\noindent {\bf Sketch of the proof. } Suppose that $a$ is the
amplitude of $A$. Choose a function $\psi_4\in C^\infty(G\times G)$
such that $\psi_4=1$ in a neighborhood of the set $\{ (x, y)\in
G\times G\ |\ x=y \}$ and $\supp\psi_4$ is proper. If we write
$a^0(t, \omega, x, y, \xi)=a(t, \omega, x, y, \xi)\psi_4(x, y)$,
then for a.e. $(t, \omega)\in (0, T)\times\Omega$ and any
$\xi\in\dbR^n$, $\supp a^0(t, \omega, \cdot, \cdot,
\xi)\subseteq\supp\psi_4$. It follows that the SPDO $A^0$ determined
by $a^0$ is a uniformly properly supported SPDO with respect to $(t,
\omega)$. On the other hand, let $K_A$ be the kernel of $A$, then
$(1-\psi_4)K_A\in \mathcal{E}_s(G\times G)$.  If we denote by $A^1$
the operator, whose kernel is $(1-\psi_4)K_A$, then its amplitude
$a^1$ turns out to be $(1-\psi_4)a$. It is easy to check that
$(1-\psi_4)K_A$ satisfies the conditions (1)-(3) mentioned in Remark
\ref{33}. Therefore, $a^1\in S^{-\infty}_p(G\times G)$, $A^1$ is a
smoothing operator and $A=A^0+A^1$.
\endpf

\medskip

For a uniformly properly supported SPDO with respect to $(t,
\omega)$ determined by an amplitude, one can reduce it to the form
represented by a symbol.
\begin{theorem}\label{t44}
Suppose that $A$ is a uniformly properly supported SPDO with respect
to $(t, \omega)$. Then for any $u\in L^q_{\mathcal{F}}(0, T;
\mathcal{D}(G))$,
$$(Au)(t, \omega, x)=(2\pi)^{-n}\displaystyle\int_{\dbR^n} e^{i
x\cdot\xi}\sigma_A(t, \omega, x, \xi)\hat{u}(t, \omega, \xi)d\xi,$$
where $\sigma_A(t, \omega, x, \xi)=e^{-i x\cdot\xi}(A e^{i
x\cdot\xi})(t, \omega, x)$.
\end{theorem}
The proof of Theorem \ref{t44} is similar to that in the
deterministic case. Here we omit it.

\medskip

For any SPDO $A$, by Theorem \ref{t43}, $A$ can be rewritten as a
sum of a uniformly properly supported SPDO $A^0$ with respect to
$(t, \omega)$ and  a smoothing operator $A^1$. On the other hand, by
Theorem \ref{t44}, $A^0$ can be represented by a symbol
$\sigma_{A^0}$.
 We call $\sigma_{A^0}$ a principal symbol of $A$ and denote it by
 $\sigma_A$ simply.

\medskip

 In the remainder of this subsection, we give an
 equivalent characterization of SPDOs and omit the proof, since it is similar to that in the deterministic
 case.
 \begin{corollary}
$(1)$ $A\in \mathcal{L}^{\ell}_p(G\times G)$ if and only if for a.e.
$(t, \omega)\in (0, T)\times\Omega$, $A: C^\infty_0(G)\rightarrow
C^\infty(G)$ is a continuous linear operator, and for any
$\varphi\in C^\infty_0(G)$, $e^{-ix\cdot\xi}A(\varphi e^{i
x\cdot\xi})\in S^{\ell}_p(G\times \dbR^n)$;

\noindent $(2)$ $A\in \mathcal{L}^{\ell}_p(G\times G)$ is uniformly
properly supported with respect to $(t, \omega)$ if and only if for
a.e. $(t, \omega)\in (0, T)\times\Omega$, $A: C^\infty(G)\rightarrow
C^\infty(G)$ is a continuous linear operator, its distribution
kernel $K_A$ is uniformly properly supported with respect to $(t,
\omega)$, and $e^{-ix\cdot\xi}A(e^{i x\cdot\xi})\in
S^{\ell}_p(G\times \dbR^n)$.
\end{corollary}

\subsection{Asymptotic expansions of a symbol}

In this subsection, we give the notion of asymptotic expansions of a
symbol and present some usual results on asymptotic expansions.
First, we have the following definition.
\begin{definition}\label{d51}
Suppose that $\{\ell_j\}_{j\in\dbN}$ is a monotone decreasing
sequence satisfying $\ell_j\rightarrow -\infty$ $(j\rightarrow
\infty)$ and $a_j\in S^{\ell_j}_p(G\times \dbR^n)$. Then
$\sum\limits_{j=0}^\infty a_j$ is called asymptotic expansions of a
symbol  $a\in S^{\ell_0}_p(G\times\dbR^n)$ if
\begin{equation}\label{e52}a
-\sum\limits^{k-1}_{j=0}a_{j}\in S^{\ell_k}_p(G\times\dbR^n), \
\forall\ k\in \dbN,\end{equation} and we write
$a\sim\sum\limits_{j=0}^\infty a_j$.
\end{definition}

Next, we give a useful lemma.
\begin{lemma}\label{t51}
Suppose that $\{\ell_j\}_{j\in\dbN}$ is a monotone decreasing
sequence satisfying $\ell_j\rightarrow -\infty$ $(j\rightarrow
\infty)$ and $a_j\in S^{\ell_j}_p(G\times \dbR^n)$. Then there
exists a symbol $a\in S^{\ell_0}_p(G\times\dbR^n)$ such that
$a\sim\sum\limits_{j=0}^\infty a_j.$
\end{lemma}
\noindent{\bf Proof. } Pick a sequence of compact sets
$\{K_j\}_{j\in\dbN}$ such that $K_0\subseteq K_1\subseteq\cdots$ and
$\displaystyle\bigcup_{j\in\dbN}K_j=G$, and choose a function
$\psi_5\in C^\infty(\dbR^n)$ satisfying $\psi_5(\xi)=0$ for
$|\xi|\leq \displaystyle\frac{1}{2}$  and $\psi_5(\xi)=1$ for
$|\xi|\geq 1$. We construct a function $a$ of the form
\begin{equation}\label{e51}
a(t, \omega, x, \xi)=\sum\limits_{j=0}^\infty
\psi_5(\varepsilon_j\xi)a_j(t, \omega, x, \xi),
\end{equation}
where $\varepsilon_j$ are  small constants, which will be specified
later. It is easy to verify that $a\in
\mathcal{E}_s(G\times\dbR^n)$.

On the other hand, since for any
$k\in\dbN$,\begin{equation}\label{51}a(t, \omega, x,
\xi)-\displaystyle\sum_{j=0} ^{k-1}a_j(t, \omega, x,
\xi)=\displaystyle\sum_{j=0} ^{k-1}[\psi_5(\varepsilon_j
\xi)-1]a_j(t, \omega, x, \xi)+\displaystyle\sum_{j=k}
^{\infty}\psi_5(\varepsilon_j \xi)a_j(t, \omega, x, \xi),
\end{equation}
and $\displaystyle\sum_{j=0} ^{k-1}[\psi_5(\varepsilon_j
\xi)-1]a_j(t, \omega, x, \xi)\in S^{-\infty}_{p}(G\times\dbR^n)$, it
remains to prove that $\displaystyle\sum_{j=k}
^{\infty}\psi_5(\varepsilon_j \xi)a_j(t, \omega, x, \xi)\in
S^{\ell_k}_{p}(G\times\dbR^n)$. Indeed, for any compact set
$K\subseteq G$, there exists an $i^*\in\dbN$ such that $K\subseteq
K_{i^*}$. For any two multi-indices $\alpha$ and $\beta$, write
$k^*=\max\{k, |\alpha|+|\beta|+i^*\}$. Then
$$
\displaystyle\sum_{j=k} ^{\infty}\psi_5(\varepsilon_j \xi)a_j(t,
\omega, x, \xi)=\displaystyle\sum_{j=k} ^{k^*}\psi_5(\varepsilon_j
\xi)a_j(t, \omega, x, \xi)+\displaystyle\sum_{j=k^*+1}
^{\infty}\psi_5(\varepsilon_j \xi)a_j(t, \omega, x, \xi).
$$
For a.e. $(t, \omega)\in(0, T)\times\Omega$ and any $(x, \xi)\in
K\times\dbR^n$,
\begin{eqnarray}\label{52}
\begin{array}{rl}
&\left|\partial^\alpha_\xi\partial^\beta_x
\left[\displaystyle\sum_{j=k} ^{k^*}\psi_5(\varepsilon_j \xi)a_j(t,
\omega, x, \xi)\right]\right| \leq C(\alpha, \beta, k,
K)M^1_{\alpha, \beta, K}(t, \omega) (1+|\xi|)^{\ell_k-|\alpha|},
\end{array}
\end{eqnarray}
for a function $M^1_{\alpha, \beta, K}(\cdot, \cdot)\in
L^p_{\mathcal{F}}(0, T)$. Also, for any $j\geq k^*+1$,
\begin{eqnarray*}
&&\left|\partial^\alpha_\xi\partial^\beta_x \left[
\psi_5(\varepsilon_j \xi)a_j(t, \omega, x, \xi)\right]\right| \leq
\displaystyle\sum_{\alpha_1+\alpha_2=\alpha}
 \displaystyle\frac{\alpha!}{\alpha_1!\alpha_2!}\left|
 \partial^{\alpha_1}_{\xi}\psi_5(\varepsilon_j \xi)
 \partial_\xi^{\alpha_2}\partial^\beta_x
 a_j(t, \omega, x, \xi) \right|\\[2mm]
&&\leq C(j, \alpha)M^2_{j}(t,
\omega)\left(1+\frac{1}{2\varepsilon_j}\right)
^{\ell_j-\ell_{j-1}}(1+|\xi|)^{\ell_{j-1}-|\alpha|},
\end{eqnarray*}
for a function $M^2_{j}(\cdot, \cdot)\in L^p_{\mathcal{F}}(0, T)$.
We choose $\varepsilon_j$ $(j\geq k^*+1)$ sufficiently small such
that
$$
\left(1+\frac{1}{2\varepsilon_j}\right) ^{\ell_j-\ell_{j-1}}\leq
\frac{2^{-j}}{C(j, \alpha)|M^2_{j}(\cdot,
\cdot)|_{L^p_{\mathcal{F}}(0, T)}+1}.
$$
Then,
$$\left\{\sum\limits_{j=k^*+1}^{N} C(j,
\alpha)M^2_{j}(\cdot, \cdot)\left(1+\frac{1}{2\varepsilon_j}\right)
^{\ell_j-\ell_{j-1}}\right\}_{N\in\dbN}$$ is a Cauchy sequence in
$L^p_{\mathcal{F}}(0, T)$. Therefore, we can find a nonnegative
function $M^2(\cdot, \cdot)\in L^p_{\mathcal{F}}(0, T)$ such that
$$
\sum\limits_{j=k^*+1}^{\infty} C(j, \alpha)M_{j}(t,
\omega)\left(1+\frac{1}{2\varepsilon_j}\right)
^{\ell_j-\ell_{j-1}}=M^2(t, \omega)\quad\quad \mbox{ for a.e. }(t,
\omega)\in(0, T)\times\Omega.
$$
This implies that
\begin{equation}\label{53}\left|\partial^\alpha_\xi\partial^\beta_x
\left[\displaystyle\sum_{j=k^*+1} ^{\infty}\psi_5(\varepsilon_j
\xi)a_j(t, \omega, x, \xi)\right]\right|\leq M^2(t,
\omega)(1+|\xi|)^{\ell_k-|\alpha|}.\end{equation} Hence, by
(\ref{51})-(\ref{53}), $a\sim\sum\limits_{j=0}^\infty a_j.$ Since
$a-a_0\in S^{\ell_1}_p(G\times\dbR^n)$ and $a_0\in
S^{\ell_0}_p(G\times\dbR^n)$, we get that $a\in
S^{\ell_0}_p(G\times\dbR^n)$.
\endpf

\medskip

Lemma \ref{t51} leads to the following theorem, which shows that the
asymptotic relation $a\sim\sum\limits_{j=0}^\infty a_j$ is valid if
an apparently weaker condition than (\ref{e52}) is assumed to hold.
\begin{theorem}\label{t52}
Suppose that $\{\ell_j\}_{j\in\dbN}$ is a monotone decreasing
sequence satisfying $\ell_j\rightarrow -\infty$ $(j\rightarrow
\infty)$ and $a_j\in S^{\ell_j}_p(G\times \dbR^n)$. If $a\in
\mathcal{E}_{s}(G\times\dbR^n)$
satisfies the following conditions:\\[2mm]
\noindent $(1)$ for any compact set $K\subseteq G$, and
multi-indices $\alpha$ and $\beta$, there exist  a  function
$M_{\alpha, \beta, K}(\cdot)\in L^p_{\mathcal{F}}(0, T)$ and a
constant $\rho=\rho(\alpha, \beta, K)\in \dbR$ such that
\begin{equation}\label{e53}\left|\partial^\beta_x\partial^\alpha_\xi a(t, \omega, x,
\xi)\right|\leq M_{\alpha, \beta, K}(t,
\omega)(1+|\xi|)^\rho,\end{equation} for a.e. $(t, \omega)\in
(0, T)\times\Omega$ and any $(x, \xi)\in K\times \dbR^n$;\\[1mm]
\noindent $(2)$ there exist a sequence of real numbers
$\{\rho_j\}_{j\in\dbN}$ satisfying $\rho_j\rightarrow -\infty\
(j\rightarrow\infty)$ and a sequence of functions $\{d_j(\cdot,
\cdot)\}_{j\in\dbN}\subseteq L^p_{\mathcal{F}}(0, T)$, such that for
any $j\in\dbN$,
\begin{equation}\label{e54}\left|a(t, \omega, x,
\xi)-\sum\limits_{k=0}^{j-1} a_k(t, \omega, x, \xi)\right| \leq
d_j(t, \omega)(1+|\xi|)^{\rho_j},\end{equation} for a.e. $(t,
\omega)\in (0, T)\times\Omega$ and any $(x, \xi)\in K\times \dbR^n$,
then $a\in S^{\ell_0}_p(G\times\dbR^n)$ and
$a\sim\sum\limits_{j=0}^\infty a_j.$
\end{theorem}

\medskip

Based on Theorem \ref{t52}, we  can get some results on asymptotic
expansions. The first one is the following asymptotic expansions of
a  principal symbol.
\begin{proposition}\label{t53}
Suppose that $A$ is a SPDO  and $a$ is its amplitude. Then its
principal symbol $\sigma_A$ satisfies that
$$\sigma_A(t, \omega, x, \xi)\sim\sum\limits_{|\alpha|\geq 0}
\frac{1}{\alpha ! \
i^{|\alpha|}}\partial^\alpha_\xi\partial^\alpha_y a(t, \omega, x, y,
\xi)|_{y=x}.$$
\end{proposition}

\medskip

The next result is the asymptotic expansions for a principal symbol
of transpose operators.

\begin{proposition}\label{86} Suppose that $A$ is a
 SPDO.
Then a principal symbol $\sigma_{^tA}$ of its transpose operator has
the following asymptotic expansions:
$$
\sigma_{^tA}(t, \omega, x,
\xi)\sim\sum\limits_{|\alpha|=0}^{\infty}\frac{1}{\alpha!i^{|\alpha|}}\partial^\alpha_\xi
\partial^\alpha_x \sigma_A(t, \omega, x, -\xi).
$$
\end{proposition}

Finally, for two uniformly properly supported SPDOs with respect to
$(t, \omega)$, we present the asymptotic expansions for a  principal
symbol of their composition operator.
\begin{proposition}\label{85}
Let $\ell_1, \ell_2\in\dbR$.  Suppose that $A\in
\mathcal{L}^{\ell_1}_{p}(G\times G)$ and $B\in
\mathcal{L}^{\ell_2}_{q}(G\times G)$ are uniformly properly
supported SPDOs with respect to $(t, \omega)$. Then  the composition
operator $B\circ A$ is a SPDO of order $(\ell_1+\ell_2, q^*)$.
Moveover, its principal symbol $\sigma_{B\circ A}$ has the following
asymptotic expansions:
$$
\sigma_{B\circ A}(t, \omega, x,
\xi)\sim\sum\limits_{|\alpha|=0}^\infty
\frac{1}{\alpha!i^{|\alpha|}}\partial^\alpha_\xi \sigma_B(t, \omega,
x, \xi)\cdot \partial^\alpha_x \sigma_A(t, \omega, x, \xi).
$$
\end{proposition}
The proofs of  Theorem \ref{t52} and Propositions \ref{t53}-\ref{85}
are similar to those for the usual pseudo-differential operators in
the deterministic case. Here we omit them.

\subsection{Algebra and generalized module of stochastic pseudo-differential operators}

In this subsection, we establish an algebra and generalized module
 of SPDOs. For this purpose, first of all, we give the following
result on the  composition of two uniformly properly supported SPDOs
with respect to $(t, \omega)$.
\begin{theorem}\label{t61}
Suppose that $A$ and $B$ are uniformly properly supported SPDOs with
respect to $(t, \omega)$. Then $B\circ A$ is  a uniformly properly
supported SPDO with respect to $(t, \omega)$.
\end{theorem}
\noindent{\bf Proof. } Denote by $a$ and $b$  the amplitudes of $A$
and $B$, respectively. Then for any $u\in \mathcal{D}_{s}(G)$,
\begin{eqnarray*}
&&((B\circ A)u)(t, \omega,
x)=(2\pi)^{-n}\displaystyle\int_{\dbR^n}\int_{\dbR^n}
e^{i(x-y)\cdot\xi} b(t,
\omega, x, y, \xi)(Au)(t, \omega, y)dyd\xi\\[2mm]
&&=(2\pi)^{-2n}\displaystyle\int_{\dbR^n}\int_{\dbR^n}
e^{i(x-y)\cdot\xi} b(t, \omega, x, y,
\xi)\left[\displaystyle\int_{\dbR^n}\int_{\dbR^n}e^{i
(y-z)\cdot\varsigma}a(t, \omega, y, z,
\varsigma)u(t, \omega, z)dzd\varsigma\right]dyd\xi\\[2mm]
&&=(2\pi)^{-2n}\displaystyle\int_{\dbR^n}\int_{\dbR^n}
e^{i(x-z)\cdot\xi}\left[\displaystyle\int_{\dbR^n}\int_{\dbR^n}
e^{i(z-y)(\xi-\varsigma)} b(t, \omega, x, y, \xi)a(t, \omega, y, z,
\varsigma)dyd\varsigma\right]u(t, \omega, z)dzd\xi.
\end{eqnarray*}
Set \begin{eqnarray*} &&c(t, \omega, x, z,
\xi)=(2\pi)^{-n}\displaystyle\int_{\dbR^n}\int_{\dbR^n}
e^{i(z-y)(\xi-\varsigma)} b(t, \omega, x, y, \xi)a(t, \omega, y, z,
\varsigma)dyd\varsigma\\[2mm]
&&=(2\pi)^{-n}\displaystyle\int_{\dbR^n}\int_{\dbR^n} e^{i(z-y)\eta}
b(t, \omega, x, y, \xi)a(t, \omega, y, z, \xi-\eta)dyd\eta.
\end{eqnarray*}
Since $A$ and $B$ are two uniformly properly supported SPDOs with
respect to $(t, \omega)$, there exist two proper sets $G_1$ and
$G_2$, such that supp $a(t, \omega, \cdot, \cdot, \xi)\subseteq G_1$
and supp $b(t, \omega, \cdot, \cdot, \xi)\subseteq G_2$ for a.e.
$(t, \omega)\in(0, T)\times\Omega$ and any $\xi\in \dbR^n$. If we
write
$$G^*=\{(x, z)\in G\times G\ |
\ \mbox{there exists a }y\in G, \mbox{ such that }(x, y)\in G_2, \
(y, z)\in G_1  \},$$ then it is easy to see that $G^*$ is proper and
supp $c(t, \omega, \cdot, \cdot, \xi)\subseteq G^*$ for a.e. $(t,
\omega)\in(0, T)\times\Omega$ and any $\xi\in \dbR^n$. This finishes
the proof.
\endpf

\medskip

In the following, we establish an algebra and generalized module of
SPDOs. Notice that a generalized module means a usual module, which
does not satisfy the associative law. For any given $p\in [1,
\infty]$,
$\displaystyle\bigcup_{\ell\in\dbR}\mathcal{L}^\ell_p(G\times G)$ is
a linear space. Also, for any $A$, $B\in
\displaystyle\bigcup_{\ell\in\dbR}\mathcal{L}^\ell_p(G\times G)$,
denote by $a$ and $b$  their amplitudes respectively and define an
equivalent relation $\sim:$ $$A\sim B\Leftrightarrow a-b\in
S^{-\infty}_p(G\times G\times \dbR^n)\Leftrightarrow A-B\in
\displaystyle\bigcap_{\ell\in\dbR}\mathcal{L}^\ell_p(G\times G).$$
In every equivalent class $[\cdot]$, there exists a uniformly
properly supported SPDO with respect to $(t, \omega)$. Therefore,
for any two equivalent classes $\mathcal{A}$ and $\mathcal{B}$, take
$A\in \mathcal{A}$ and $B\in \mathcal{B}$ such that $A$ and $B$ are
two uniformly properly supported SPDOs with respect to $(t,
\omega)$. We define   the composition of $\mathcal{A}$ and
$\mathcal{B}$ as follows:
$$\mathcal{B}\circ\mathcal{A}=[B\circ A].$$
It is easy to check that the definition is well posed. Moreover,
$\displaystyle\bigcup_{\ell\in\dbR}\mathcal{L}^\ell_\infty(G\times
G)$ constructs an algebra in the above sense, and
$\displaystyle\bigcup_{\ell\in\dbR}\mathcal{L}^\ell_p(G\times G)$
$(p\geq 1)$ is a generalized module over
$\displaystyle\bigcup_{\ell\in\dbR}\mathcal{L}^\ell_\infty(G\times
G)$.

\begin{remark}
For any $\mathcal{A}$, $\mathcal{B}\in
\displaystyle\bigcup_{\ell\in\dbR}\mathcal{L}^\ell_\infty(G\times
G)$, define $[\mathcal{A},
\mathcal{B}]=\mathcal{A}\circ\mathcal{B}-\mathcal{B}\circ\mathcal{A}$.
Then we conjecture that $[\cdot, \cdot]$ is a Poisson bracket
defined on
$\displaystyle\bigcup_{\ell\in\dbR}\mathcal{L}^\ell_\infty(G\times
G)$. However, we do not give a precise proof at this moment.
\end{remark}

\section{Boundedness of stochastic pseudo-differential operators}

Many papers have been devoted to a study of the continuity of
pseudo-differential operators in $L^p$ spaces. In this section, we
shall establish the $L^p$-estimates of SPDOs. The main idea borrows
from that used in \cite{h1}   and \cite{m}. However, different from
the deterministic case, results here involve integrability with
respect to the variables of time and sample point, and we have to
deal with the problem, generated by global estimates. For simplicity
of notation, for any bounded linear operator $L$, we denote by $|L|$
its operator norm.

\subsection{$L^2$-boundedness}

In this subsection, we first present the  $L^2$-estimates for SPDOs.
To begin with, we prove the following basic estimate.
\begin{lemma}\label{76}
Suppose that $\widetilde{K}$ is a function defined on $(0, T)\times
\Omega\times\dbR^n\times\dbR^n$ such that
\begin{equation}\label{77}
\sup\limits_{y\in\dbR^n}\int_{\dbR^n} |\widetilde{K}(t, \omega, x,
y)|dx\leq M(t, \omega),\quad\quad
\sup\limits_{x\in\dbR^n}\int_{\dbR^n} |\widetilde{K}(t, \omega, x,
y)|dy\leq M(t, \omega),
\end{equation}
for a nonnegative function $M(\cdot, \cdot)$ defined on $(0,
T)\times \Omega$. Then, for the linear operator $L$ defined as
follows: $(Lu)(t, \omega,
x)=\displaystyle\int_{\dbR^n} \widetilde{K}(t, \omega, x, y)u(t, \omega, y)dy$,
the following conclusions hold:\\[2mm]
\noindent $(1)$ for a.e. $(t, \omega)\in (0, T)\times \Omega$, $L:
L^p(\dbR^n)\rightarrow L^p(\dbR^n)$  is bounded and $ |L|\leq M(t, \omega)$;\\[2mm]
\noindent $(2)$ if $M(\cdot, \cdot)\in L^\infty_{\mathcal{F}}(0,
T)$, then $L: L^{q}_{\mathcal{F}}(0, T; L^p(\dbR^n))\rightarrow
L^{q}_{\mathcal{F}}(0, T; L^p(\dbR^n))$ is bounded and $|L|\leq
|M|_{L^\infty_{\mathcal{F}}(0, T)}$; \\[2mm]
\noindent $(3)$ if $M(\cdot, \cdot)\in
L^{\widetilde{q}}_{\mathcal{F}}(0, T)$ $(\widetilde{q}\in[1,
\infty))$, then $L: L^{q}_{\mathcal{F}}(0, T;
L^p(\dbR^n))\rightarrow L^{\widehat{q}}_{\mathcal{F}}(0, T;
L^p(\dbR^n))$ is bounded and $|L|\leq
|M|_{L^{\widetilde{q}}_{\mathcal{F}}(0, T)}$, where
$\widehat{q}=\frac{q\widetilde{q}}{q+\widetilde{q}}$ for
$\widetilde{q},\ q\geq 1$, $q\widetilde{q}\geq q+\widetilde{q}$;
 and $\widehat{q}=\widetilde{q}$ for $\widetilde{q}\geq 1$, $q=\infty$.
\end{lemma}

\noindent {\bf Proof. } For any $p>1$, a.e. $(t, \omega)\in(0,
T)\times \Omega$ and any $u\in L^p(\dbR^n)$, by H\"{o}lder's
inequality, we have that
\begin{eqnarray*}
&&|(Lu)(t, \omega, x)|^p\leq \int_{\dbR^n}|\widetilde{K}(t, \omega,
x, y)||u(y)|^p dy \cdot \left(\int_{\dbR^n}|\widetilde{K}(t, \omega,
x,
y)|dy\right)^{\frac{p}{p'}}\\[2mm]
&&\leq M^{\frac{p}{p'}}(t, \omega)\int_{\dbR^n}|\widetilde{K}(t,
\omega, x, y)||u(y)|^p dy,
\end{eqnarray*}
here and hereafter $p'$ denotes a constant satisfying
$\frac{1}{p}+\frac{1}{p'}=1$. Integrating on $\dbR^n$ with respect
to the variable $x$,  by (\ref{77}), we get that
\begin{eqnarray}\label{79}
&&\int_{\dbR^n}|(Lu)(t, \omega, x)|^pdx\leq M^{1+\frac{p}{p'}}(t,
\omega)\int_{\dbR^n}|u(y)|^p dy.
\end{eqnarray}
This means that for a.e. $(t, \omega)\in (0, T)\times \Omega$, $L$
is a bounded operator from $L^p(\dbR^n)$ to $L^p(\dbR^n)$ and
$|L|\leq M(t, \omega)$.
If $M\in L^{\infty}_{\mathcal{F}}(0, T)$,
then by (\ref{79}),
 for any $u\in L^{q}_{\mathcal{F}}(0, T;
L^p(\dbR^n))$ and a.e. $(t, \omega)\in (0, T)\times \Omega$,
\begin{eqnarray*}
&&|(Lu)(t, \omega,
\cdot)|_{L^p(\dbR^n)}\leq|M|_{L^{\infty}_{\mathcal{F}}(0, T)}|u(t,
\omega, \cdot)|_{L^p(\dbR^n)}.
\end{eqnarray*}
It follows that
\begin{eqnarray*}
&&\dbE\int^T_0|(Lu)(t, \omega, \cdot)|^{q}_{L^p(\dbR^n)}dt\leq
|M|^{q}_{L^{\infty}_{\mathcal{F}}(0, T)}\dbE\int^T_0 |u(t, \omega,
\cdot)|^{q}_{L^p(\dbR^n)}dt,\mbox{ for }q\geq 1;\\[1mm]
&&|Lu|_{L^\infty_{\mathcal{F}}(0, T;
L^p(\dbR^n))}\leq|M|_{L^{\infty}_{\mathcal{F}}(0,
T)}|u|_{L^\infty_{\mathcal{F}}(0, T; L^p(\dbR^n))},\mbox{ for
}q=\infty,
\end{eqnarray*}
which implies that $L: L^{q}_{\mathcal{F}}(0, T;
L^p(\dbR^n))\rightarrow L^{q}_{\mathcal{F}}(0, T; L^p(\dbR^n))$ is
bounded and $|L|\leq |M|_{L^\infty_{\mathcal{F}}(0, T)}$.
Furthermore, if $M\in L^{\widetilde{q}}_{\mathcal{F}}(0, T)$
$(\widetilde{q}\geq 1)$, by (\ref{79}), for any $u\in
L^{q}_{\mathcal{F}}(0, T; L^p(\dbR^n))$, we obtain that
\begin{eqnarray*}
|(Lu)(t, \omega, \cdot)|^{\widehat{q}}_{L^p(\dbR^n)}\leq
M^{\widehat{q}}(t, \omega)|u(t, \omega,
\cdot)|^{\widehat{q}}_{L^p(\dbR^n)}.
\end{eqnarray*}
Hence,
\begin{eqnarray*}
&&\dbE\int^T_0|(Lu)(t, \omega,
\cdot)|^{\widehat{q}}_{L^p(\dbR^n)}dt\leq
\left[\dbE\int^T_0M^{\widetilde{q}}(t,
\omega)dt\right]^{\frac{\widehat{q}}{\widetilde{q}}}\cdot\left[\dbE\int^T_0|u(t,
\omega,
\cdot)|^{q}_{L^p(\dbR^n)}dt\right]^{\frac{\widehat{q}}{q}},\mbox{
for }q\geq
1;\\[1mm]
&&\dbE\int^T_0|(Lu)(t, \omega,
\cdot)|^{\widehat{q}}_{L^p(\dbR^n)}dt\leq
\dbE\int^T_0M^{\widetilde{q}}(t, \omega)dt\cdot
|u|^{\widehat{q}}_{L^\infty_{\mathcal{F}}(0, T; L^p(\dbR^n))},\mbox{
for }q=\infty.
\end{eqnarray*}
This implies that $L: L^{q}_{\mathcal{F}}(0, T;
L^p(\dbR^n))\rightarrow L^{\widehat{q}}_{\mathcal{F}}(0, T;
L^p(\dbR^n))$ is bounded and $|L|\leq
|M|_{L^{\widetilde{q}}_{\mathcal{F}}(0, T)}$. Results for the cases
of $p=1$ and $p=\infty$ can be derived in the same way.
\endpf

\medskip

Next, we give the notion of adjoint operators.
\begin{definition}\label{83}
Suppose that $A$ is a SPDO of order $(\ell, p)$ and $a$ is its
amplitude. A linear operator $A^*$ is called the adjoint operator of
$A$ if for any $u\in L^q_{\mathcal{F}}(0, T; \mathcal{D}(G))$,
$$(A^*u)(t, \omega, x)=(2\pi)^{-n}\displaystyle\int_{\dbR^n}\displaystyle\int_{\dbR^n}
e^{i (x-y)\cdot\xi}\overline{a(t, \omega, y, x, \xi)}u(t, \omega,
y)dyd\xi,$$where $\overline{z}$ denotes the conjugate of a complex
number $z$.
\end{definition}
\begin{remark}\label{84}
It is easy to show that for any SPDO $A\in
\mathcal{L}^{\ell}_p(G\times G)$, $A^*\in
\mathcal{L}^{\ell}_{p}(G\times G)$. Moreover,  for a.e. $(t,
\omega)\in(0, T)\times\Omega$ and any $u,\ v\in C^\infty_0(G)$,
$((A^*u)(t, \omega, \cdot), v(\cdot))_{L^2(G)}=(u(\cdot), (Av)(t,
\omega, \cdot))_{L^2(G)}$.
\end{remark}
\begin{remark}\label{112}
Similar to  Proposition \ref{86}, for any SPDO $A$, it is easy to
show that a principal symbol of its adjoint operator $A^*$ has the
following asymptotic expansions:
$$
\sigma_{A^*}(t, \omega, x,
\xi)\sim\sum\limits_{|\alpha|=0}^{\infty}\frac{1}{\alpha!i^{|\alpha|}}\partial^\alpha_\xi
\partial^\alpha_x \overline{\sigma_A(t, \omega, x, \xi)}.
$$
\end{remark}

Now, we give the $L^2$-estimates for a class of the SPDOs of order
$(0, \infty)$.
\begin{theorem}\label{80}
Suppose that $A$ is a SPDO and $a$ is its symbol. If $a\in
S_{\infty}^{0}$, $A: L^{q}_{\mathcal{F}}(0, T; L^2(\dbR^n))
\rightarrow L^{q}_{\mathcal{F}}(0, T; L^2(\dbR^n))$ is bounded.
\end{theorem}

\noindent {\bf Proof. } {\bf Step 1. }First, we prove that $A:
L^q_{\mathcal{F}}(0, T; L^2(\dbR^n))\rightarrow
L^{q^*}_{\mathcal{F}}(0, T; L^2(\dbR^n))$ is bounded,
 if $a\in S^{-n-1}_p$. By the definition of a
kernel $K_A$, we have that
\begin{eqnarray}\label{81}
\begin{array}{rl}
&|K_A(t, \omega, x, y)|=(2\pi)^{-n}\left|\displaystyle\int_{\dbR^n}
e^{i(x-y)\cdot\xi}a(t, \omega, x, \xi)d\xi\right|\leq
(2\pi)^{-n}\displaystyle\int_{\dbR^n}|a(t, \omega, x,
\xi)|d\xi\\[2mm]
&\leq C(n)M_0(t,
\omega)\displaystyle\int_{\dbR^n}(1+|\xi|)^{-1-n}d\xi\leq C(n)M_0(t,
\omega),
\end{array}
\end{eqnarray}
where $M_0(\cdot, \cdot)\in L^p_{\mathcal{F}}(0, T)$ is a
nonnegative function. On the other hand, for any multi-index
$\alpha$,
\begin{eqnarray}\label{82}
\begin{array}{rl}
&\left|(x-y)^\alpha K_A(t, \omega, x,
y)\right|=(2\pi)^{-n}\left|\displaystyle\int_{\dbR^n}
\partial^{\alpha}_{\xi}e^{i(x-y)\cdot
\xi}a(t,
\omega, x, \xi)d\xi\right|\\[2mm]
&=(2\pi)^{-n}\left|\displaystyle\int_{\dbR^n}e^{i(x-y)\cdot\xi}\partial^{\alpha}_{\xi}a(t,
\omega, x, \xi)d\xi\right|\\[2mm]
&\leq C(n)M_\alpha(t, \omega)\displaystyle\int_{\dbR^n}
(1+|\xi|)^{-1-n-|\alpha|}d\xi\leq C(n)M_\alpha(t, \omega),
\end{array}
\end{eqnarray}
where $M_\alpha(\cdot, \cdot)\in L^p_{\mathcal{F}}(0, T)$ is a
nonnegative function. Then, by (\ref{81}) and (\ref{82}), it follows
that for a nonnegative function $M_1(\cdot, \cdot)\in
L^p_{\mathcal{F}}(0, T)$,
$$
|K_A(t, \omega, x, y)|\leq \displaystyle\frac{C(n)M_1(t,
\omega)}{(1+|x-y|)^{1+n}},
$$
which implies that
$$\sup\limits_{y\in\dbR^n}\displaystyle\int_{\dbR^n} \left|K_A(t,
\omega, x, y)\right|dx\leq C(n)M_1(t, \omega),\quad\quad
\sup\limits_{x\in\dbR^n}\displaystyle\int_{\dbR^n} \left|K_A(t,
\omega, x, y)\right|dy\leq C(n)M_1(t, \omega).
$$
By (2)-(3) in Lemma \ref{76}, we get the desired result.

\noindent {\bf Step 2. } Suppose that
 $a\in S^{\ell}_\infty\ (\ell<0)$. By Remark \ref{84}, for any SPDO
$A\in \mathcal{L}^{\frac{-n-1}{2}}_{\infty}$, $A^*\circ A\in
\mathcal{L}^{-n-1}_{\infty}$. Then, by the proof in Step 1 and
(\ref{79}), for a.e. $(t, \omega)\in (0, T)\times \Omega$,
\begin{eqnarray*}
&&|(Au)(t, \omega, \cdot)|^2_{L^2(\dbR^n)}=((Au)(t, \omega, \cdot),
(Au)(t, \omega, \cdot))_{L^2(\dbR^n)}\\[2mm]
&&=(((A^*\circ A)u)(t, \omega, \cdot), u(t, \omega,
\cdot))_{L^2(\dbR^n)}\leq M_2^2(t, \omega)|u(t, \omega,
\cdot)|^2_{L^2(\dbR^n)},
\end{eqnarray*}
for a nonnegative function $M_2(\cdot,\cdot)\in
L^{\infty}_{\mathcal{F}}(0, T)$. Therefore,
\begin{eqnarray*}
&&\left| |(Au)(\cdot, \cdot,
\cdot)|_{L^2(\dbR^n)}\right|_{L^{q}_{\mathcal{F}}(0, T)}\leq
|M_2(\cdot, \cdot)|_{L^{\infty}_{\mathcal{F}}(0, T)}\left| |u(\cdot,
\cdot, \cdot)|_{L^2(\dbR^n)}\right|_{L^{q}_{\mathcal{F}}(0, T)}.
\end{eqnarray*}
This implies that for any $k\in\dbN$, $A: L^q_{\mathcal{F}}(0, T;
L^2(\dbR^n))\rightarrow L^{q}_{\mathcal{F}}(0, T; L^2(\dbR^n))$ is
bounded, if $a\in S^{\frac{-n-1}{2^k}}_{\infty}$. For any $\ell<0$,
we take a positive integer $k_2$, such that
$\ell<\frac{-n-1}{2^{k_2}}$. Since $S^{\ell}_{\infty}\subseteq
S^{\frac{-n-1}{2^{k_2}}}_{\infty}$, we get the desired result of
Theorem \ref{80}.

\noindent{\bf Step 3. }If $a\in S^0_\infty$, there exists a constant
$C_0>0$ such that $|a(t, \omega, x, \xi)|\leq C_0$ for a.e. $(t,
\omega)\in (0, T)\times \Omega$ and any $(x, \xi)\in \dbR^{2n}.$
Take $C_1>2C_0^2+1$ and define $b(t, \omega, x, \xi)=(C_1-|a(t,
\omega, x, \xi)|^2)^{1/2}$. Then $b\in S^0_\infty$ and
$b\overline{b}+a\overline{a}=C_1$. If we denote by $B$ the SPDO
determined by $b$, then $ \sigma_{B^*\circ B}+\sigma_{A^*\circ
A}-C_1\in S^{-1}_{\infty}. $ Hence, we can find a SPDO $R\in
\mathcal{L}^{-1}_\infty$ such that for a.e. $(t, \omega)\in (0,
T)\times \Omega$,
$$
|(Au)(t, \omega, \cdot)|^2_{L^2(\dbR^n)}+|(Bu)(t, \omega,
\cdot)|^2_{L^2(\dbR^n)}=C_1|u(t, \omega,
\cdot)|^2_{L^2(\dbR^n)}+((Ru)(t, \omega, \cdot), u(t, \omega,
\cdot))_{L^2(\dbR^n)}.
$$
Combining the above equality with the proof in Step 2, we see that
$$
|(Au)(t, \omega, \cdot)|^2_{L^2(\dbR^n)}\leq  C(C_1, R)|u(t, \omega,
\cdot)|^2_{L^2(\dbR^n)}.
$$
This implies that
\begin{eqnarray*}
&&\left| |(Au)(\cdot, \cdot,
\cdot)|_{L^2(\dbR^n)}\right|_{L^q_{\mathcal{F}}(0, T)}\leq C(C_1,
a)\left| |u(\cdot, \cdot,
\cdot)|_{L^2(\dbR^n)}\right|_{L^q_{\mathcal{F}}(0, T)}.
\end{eqnarray*}
Therefore, the  proof is completed.
\endpf

\medskip

In the following, we give two corollaries. The first one generalizes
the result of Theorem \ref{80} to the space $L^q_{\mathcal{F}}(0, T;
H^{\delta}(\dbR^n))$ for $\delta\in \dbR$.
\begin{corollary}\label{88}
Suppose that $A$ is a SPDO and $a$ is its symbol. If $a\in
S_{\infty}^{\ell}$, then for any $\delta\in \dbR$, $A:
L^{q}_{\mathcal{F}}(0, T; H^{\delta}(\dbR^n))\rightarrow
L^{q}_{\mathcal{F}}(0, T; H^{\delta-\ell}(\dbR^n))$ is bounded.
\end{corollary}

\noindent{\bf Sketch of the proof. } First, we recall that for the
pseudo-differential operator $\Lambda^{\delta}$, whose symbol is
$(1+|\xi|^2)^{\frac{\delta}{2}}$,
$|\upsilon|_{H^\delta(\dbR^n)}=|\Lambda^{\delta}\upsilon|_{L^2(\dbR^n)}$
for any $\upsilon\in H^{\delta}(\dbR^n)$. Therefore, for a.e. $(t,
\omega)\in (0, T)\times \Omega$ and any $u\in L^{q}_{\mathcal{F}}(0,
T; H^{\delta}(\dbR^n))$, if we write $v=\Lambda^{\delta}u$, then
$$
|(Au)(t, \omega,
\cdot)|_{H^{\delta-\ell}(\dbR^n)}=|(\Lambda^{\delta-\ell}(Au))(t,
\omega,
\cdot)|_{L^2(\dbR^n)}=|\Lambda^{\delta-\ell}A(\Lambda^{-\delta}v)(t,
\omega, \cdot)|_{L^2(\dbR^n)}.
$$
Noticing that $\Lambda^{\delta-\ell}A\Lambda^{-\delta}\in
S^0_{\infty}$, by the result of Theorem \ref{80}, we have that
\begin{eqnarray*}
&&\dbE\int^T_0 |(Au)(t, \omega,
\cdot)|_{H^{\delta-\ell}(\dbR^n)}^{q}dt =\dbE\int^T_0
|\Lambda^{\delta-\ell}A(\Lambda^{-\delta}v)(t, \omega,
\cdot)|_{L^2(\dbR^n)}^{q}dt\\[2mm]
&&\leq C(a, q, \delta, \ell)\dbE\int^T_0 |v(t, \omega,
\cdot)|_{L^2(\dbR^n)}^{q}dt=C(a, q, \delta, \ell)\dbE\int^T_0 |u(t,
\omega, \cdot)|_{H^{\delta}(\dbR^n)}^{q}dt.
\end{eqnarray*}
This finishes the proof.
\endpf

\medskip

The second corollary involves the $L^2$-estimates of the SPDOs
defined on a local domain of $\dbR^n$. For this purpose, we
introduce some locally convex topological vector spaces:
\begin{eqnarray*}
&&L^q_{\mathcal{F}}(0, T; H^{\delta}_{loc}(G))=\left\{u\left|
\mbox{for any }\psi\in C^\infty_0(G)\mbox{ and a.e. }(t,
\omega)\in(0, T)\times\Omega,\ \right.\right.\\
&&\left.\quad\quad\quad\quad\quad\quad\quad\quad\quad\quad\ (\psi
u)(t, \omega, \cdot)\in H^{\delta}(\dbR^n);\mbox{ and }|u|_{q,
\delta, \psi}=|\psi u|_{L^q_{\mathcal{F}}(0, T; H^{\delta}(\dbR^n))}
<\infty\right\},\\
&&\mbox{which is generated by a family of semi-norms }\{|\cdot|_{q,
\delta, \psi}\}_{\psi\in C^\infty_0(G)};\\[2mm]
&&L^q_{\mathcal{F}}(0, T; H^{\delta}_{K})=\left\{u\left| \mbox{for
a.e. }(t, \omega)\in(0, T)\times\Omega,\ u(t, \omega, \cdot)\in
H^{\delta}(\dbR^n)\mbox{ and } \supp u(t, \omega, \cdot)\subseteq K\right.\right\},\\
&&\mbox{which is generated by the norm
}|\cdot|_{L^q_{\mathcal{F}}(0, T; H^{\delta}(\dbR^n))} \mbox{
 and }K\mbox{ is any given compact set};\\[2mm]
&&L^q_{\mathcal{F}}(0, T; H^{\delta}_{comp}(G))=\bigcup_{K\subseteq
G\mbox{ is compact}}L^q_{\mathcal{F}}(0, T; H^{\delta}_{K}),\\
&&\mbox{ which is endowed with the inductive topology.}
\end{eqnarray*}

Then, similar to Proposition \ref{104}, it is easy to show the
following conclusions:
\begin{eqnarray*}
&&(1)\ \mbox{the locally convex space }L^q_{\mathcal{F}}(0, T;
H^{\delta}_{comp}(G))\mbox{ satisfies the first countability
axiom};\\[1mm]
&&(2)\ \mbox{the locally convex space }L^q_{\mathcal{F}}(0, T;
H^{\delta}_{loc}(G))\mbox{ is a Fr\'{e}chet space};\\[1mm]
&&(3)\ \lim\limits_{j\rightarrow \infty}u_j=0 \mbox{ in
}L^q_{\mathcal{F}}(0, T; H^{\delta}_{comp}(G)),\mbox{ if and only if
there exists a compact set }K^*\mbox{ such that }\\[1mm]
&&\quad\quad\bigcup_{j\in\dbN}\supp u_j(t, \omega, \cdot)\subseteq
K^*\mbox{ for a.e. }(t, \omega)\in(0, T)\times\Omega\mbox{ and
}\lim\limits_{j\rightarrow \infty}|u_j|_{L^q_{\mathcal{F}}(0, T;
H^{\delta}(\dbR^n))}=0;\\[1mm]
&&(4)\ \lim\limits_{j\rightarrow \infty}u_j=0\mbox{ in
}L^q_{\mathcal{F}}(0, T; H^{\delta}_{loc}(G)),\mbox{ if and only if
for any }\psi\in
C^\infty_0(G),\\[1mm]
&&\quad\quad\lim\limits_{j\rightarrow\infty}|\psi
u_j|_{L^q_{\mathcal{F}}(0, T; H^{\delta}(\dbR^n))}=0.
\end{eqnarray*}

\begin{corollary}\label{89}
Suppose that $A$ is  a SPDO and $a$ is its amplitude.  \\[2mm]
\noindent $(1)$ If $a\in S_{\infty}^{\ell}(G\times G\times\dbR^n)$,
then for any $\delta\in \dbR$, $A: L^{q}_{\mathcal{F}}(0, T;
H^{\delta}_{comp}(G))\rightarrow L^{q}_{\mathcal{F}}(0, T;
H^{\delta-\ell}_{loc}(G))$ is continuous;\\[1mm]
\noindent $(2)$ If $A$ is a uniformly properly supported SPDO with
respect to $(t, \omega)$ and $a\in S_{\infty}^{\ell}(G\times
G\times\dbR^n)$, then for any $\delta\in \dbR$, both $A:
L^{q}_{\mathcal{F}}(0, T; H^{\delta}_{comp}(G))\rightarrow
L^{q}_{\mathcal{F}}(0, T; H^{\delta-\ell}_{comp}(G))$ and $A:
L^{q}_{\mathcal{F}}(0, T; H^{\delta}_{loc}(G))\rightarrow
L^{q}_{\mathcal{F}}(0, T; H^{\delta-\ell}_{loc}(G))$ are continuous.
\end{corollary}
\noindent {\bf Proof. } For any $u\in L^{q}_{\mathcal{F}}(0, T;
H^{\delta}_{comp}(G))$, by the definition of the space $
L^{q}_{\mathcal{F}}(0, T; H^{\delta}_{comp}(G))$, there exists a
compact set $K^7$ such that $u\in L^{q}_{\mathcal{F}}(0, T;
H^{\delta}_{K^7})$. Take a function $\varphi_1\in C^\infty_0(G)$
such that $\varphi_1=1$ in $K^7$. Then for any function $\psi\in
C^\infty_0(G)$,
$$
(\psi(Au))(t, \omega, x)=(2\pi)^{-n}\displaystyle\int_{\dbR^n}\int_G
e^{i(x-y)\cdot\xi}a(t, \omega, x, y, \xi)\psi(x)\varphi_1(y)u(t,
\omega, y)dyd\xi.$$ Write $\check{a}(t, \omega, x, y, \xi)=a(t,
\omega, x, y, \xi)\psi(x)\varphi_1(y)$ and denote by $\check{A}$ the
SPDO determined by $\check{a}$. Then $\check{a}\in S^{\ell}_\infty$
and $\check{A}u=\psi (Au)$. By Corollary \ref{88}, it follows that
$\psi (Au)\in L^{q}_{\mathcal{F}}(0, T; H^{\delta-\ell}(\dbR^n))$.
This implies that $Au\in L^{q}_{\mathcal{F}}(0, T;
H^{\delta-\ell}_{loc}(G))$.

On the other hand, if $\lim\limits_{j\rightarrow \infty}u_j=0$ in
$L^{q}_{\mathcal{F}}(0, T; H^{\delta}_{comp}(G))$, then there exists
a compact set $K^{8}$ such that
$\displaystyle\bigcup_{j\in\dbN}\supp u_j(t, \omega, \cdot)\subseteq
K^{8}$ for a.e. $(t, \omega)\in(0, T)\times \Omega$ and
$\lim\limits_{j\rightarrow \infty}|u_j|_{L^{q}_{\mathcal{F}}(0, T;
H^{\delta}(\dbR^n))}=0$. Similarly, take a function $\varphi_2\in
C^\infty_0(G)$ such that $\varphi_2=1$ in $K^{8}$. Then for any
function $\psi\in C^\infty_0(G)$, $\psi (Au)=\widetilde{A}u$, where
the amplitude of $\widetilde{A}$ is $a\psi\varphi_2$ and
$a\psi\varphi_2\in S^{\ell}_{\infty}$. Therefore, by Corollary
\ref{88}, $\lim\limits_{j\rightarrow \infty}\psi (Au_j)=0$ in
$L^q_{\mathcal{F}}(0, T; H^{\delta-\ell}(\dbR^n))$. This implies
that $\lim\limits_{j\rightarrow \infty}Au_j=0$ in
$L^q_{\mathcal{F}}(0, T; H^{\delta-\ell}_{loc}(G))$. Hence, we get
the desired result (1).

Similarly, we can show that (2) also holds.
\endpf

\subsection{$L^p$-boundedness }

In this subsection, we present the $L^p$-estimates $(p>1, p\neq 2)$
of SPDOs. To begin with, we give the following known lemma, which
will be used later.

\begin{lemma}\label{90}{\rm (\cite[Page 272]{7})}
There exist two functions $\psi^*(\cdot), \varphi^*(\cdot)\in
C^\infty_0(\dbR^n)$ satisfying $0\leq \psi^*(\xi),\
\varphi^*(\xi)\leq 1$, and
\begin{eqnarray*}
&&(1)\ \mbox{supp }\psi^*(\cdot)\subseteq B_1,\ \mbox{supp
}\varphi^*(\cdot)\subseteq G_0;\quad\quad(2) \
\psi^*(\xi)+\sum\limits_{j=0}^{\infty} \varphi^*(2^{-j}\xi)=1,
\forall\ \xi\in \dbR^n,
\end{eqnarray*}
where $B_1=\{\xi\in\dbR^n; \ |\xi|=1\}$ and $G_0=\{\xi\in\dbR^n; \
k_*^{-1}<|\xi|<2k_*\}$ for a constant $k_*>1$.
\end{lemma}

Also, we give a lemma, which is a reformation of the known
Calder\'{o}n-Zygmund decomposition and is adapted to the stochastic
case.
\begin{lemma}\label{91}
Suppose that $u\in L^1(\dbR^n; L^p_{\mathcal{F}}(0, T))$. Then for
any $r>0$, there exist the functions $v(\cdot),\ w_k(\cdot)\in
L^1(\dbR^n; L^p_{\mathcal{F}}(0, T))$ $(k=1, 2, \cdots)$ such that
$u(t, \omega, x)=v(t, \omega, x)+\sum\limits_{k=1}^{\infty}w_k(t,
\omega, x)$, where $v(\cdot)$ and $w_k(\cdot)$ satisfy the following
conditions:
\begin{eqnarray*}
&&\supp w_k(t, \omega, \cdot)\subseteq \overline{I_k},\mbox{ for
a.e. }(t, \omega)\in(0, T)\times\Omega,\mbox{ where }I_k \mbox{ are
disjoint cubes in }\dbR^n; \\
&&r\sum\limits_{k=1}^{\infty}|I_k|\leq
|u(\cdot)|_{L^1(\dbR^n; L^p_{\mathcal{F}}(0, T))};\\
&&\int_{\dbR^n}w_k(t, \omega, x)dx=0,\quad |v(\cdot, \cdot,
x)|_{L^p_{\mathcal{F}}(0, T)}\leq
2^n r,\mbox{ for a.e. }(t, \omega, x)\in(0, T)\times\Omega\times\dbR^n;\\
&&|v(t, \omega,
\cdot)|_{L^1(\dbR^n)}+\sum\limits_{k=1}^{\infty}|w_k(t, \omega,
\cdot)|_{L^1(\dbR^n)}\leq 3|u(t, \omega, \cdot)|_{L^1(\dbR^n)},
\mbox{ for a.e. }(t, \omega)\in(0, T)\times\Omega;\\
&&|v(\cdot, \cdot, x)|_{L^q_{\mathcal{F}}(0, T)}=|u(\cdot, \cdot,
x)|_{L^q_{\mathcal{F}}(0, T)},\mbox{ for a.e.
}x\in(\displaystyle\bigcup_{k\in\dbN} I_k)^c;\\
&&\int_{I_k}|v(\cdot, \cdot, x)|_{L^q_{\mathcal{F}}(0, T)}dx\leq
\int_{I_k}|u(\cdot, \cdot, x)|_{L^q_{\mathcal{F}}(0, T)}dx,\mbox{
for
 any }k\in\dbN.
\end{eqnarray*}
\end{lemma}
\noindent{\bf Proof. } First, we divide $\dbR^n$ into the cubes,
whose volumes are  greater than $r^{-1}|u(\cdot)|_{L^1(\dbR^n;
L^p_{\mathcal{F}}(0, T))}$. Then for every such cube $M$,
\begin{equation}\label{113}
\displaystyle\frac{|u(\cdot)|_{L^1(M; L^p_{\mathcal{F}}(0,
T))}}{|M|}\leq \displaystyle\frac{|u(\cdot)|_{L^1(\dbR^n;
L^p_{\mathcal{F}}(0, T))}}{|M|}<r.
\end{equation}
Again, we divide every cube $M$ into $2^n$ cubes equally and denote
by $I_{1k}$ $(k=1, 2, \cdots)$ those small cubes satisfying $
\displaystyle\frac{|u(\cdot)|_{L^1(I_{1 k}; L^p_{\mathcal{F}}(0,
T))}}{|I_{1k}|}\geq r. $ By (\ref{113}), we obtain that $
|u(\cdot)|_{L^1(I_{1 k}; L^p_{\mathcal{F}}(0, T))}< 2^n r|I_{1k}|$.
Next, we divide the cubes, which are not $I_{1k}$,  into $2^n$ small
cubes equally. Denote by $I_{2k}$ $(k=1, 2, \cdots)$ those small
cubes satisfying $\displaystyle\frac{|u(\cdot)|_{L^1(I_{2 k};
L^p_{\mathcal{F}}(0, T))}}{|I_{2k}|}\geq r$. Similarly, we get that
$ |u(\cdot)|_{L^1(I_{2 k}; L^p_{\mathcal{F}}(0, T))}< 2^n
r|I_{2k}|$.

 We proceed as
the above steps, and then get a sequence of cubes (denoted by
$\{I_k\}_{k\in\dbN}$), such that the average $
\displaystyle\frac{|u(\cdot)|_{L^1(I_k; L^p_{\mathcal{F}}(0,
T))}}{|I_{k}|}$ of $|u|$ on $I_k$ is greater than $r$ for all $k\in
\dbN$. Moreover,
\begin{equation}\label{100}
|u(\cdot)|_{L^1(I_k; L^p_{\mathcal{F}}(0, T))}< 2^n r|I_{k}|.
\end{equation}
 Define
\begin{eqnarray*}
v(t, \omega, x)=\left\{
\begin{array}{lll}
&\frac{\int_{I_{k}} u(t,
\omega, x)dx}{|I_{k}|} &x\in I_k,\\
&u(t, \omega, x) &x\in (\displaystyle\bigcup_{k\in\dbN} I_k)^c;
\end{array}\right.
\end{eqnarray*}
and
\begin{eqnarray*}
w_k(t, \omega, x)=\left\{
\begin{array}{lll}
&u(t, \omega, x)-v(t, \omega, x) &x\in I_k,\\
&0 &x\in (I_k)^c.
\end{array}\right.
\end{eqnarray*}
Then, if $x\in I_k$, by (\ref{100}), $|v(\cdot, \cdot,
x)|_{L^p_{\mathcal{F}}(0,
T)}=\displaystyle\frac{\left|\int_{I_k}u(\cdot, \cdot,
x)dx\right|_{L^p_{\mathcal{F}}(0, T)}}{|I_k|}\leq
\frac{\int_{I_k}|u(\cdot, \cdot, x)|_{L^p_{\mathcal{F}}(0,
T)}dx}{|I_k|}< 2^n r,$ for $k\in\dbN$. If $x\in
(\displaystyle\bigcup_{k\in\dbN} I_k)^c$, there exists a sequence of
cubes $\{I^*_k\}_{k\in\dbN}$ satisfying $x\in
\displaystyle\bigcap_{k\in\dbN} I^*_k$ and $|I^*_k|$ tends to 0
$(k\rightarrow \infty)$, such that
$\displaystyle\frac{|u(\cdot)|_{L^1(I_k^*; L^p_{\mathcal{F}}(0,
T))}}{|I_k^*|}< r$. Therefore, we get that $|v(\cdot, \cdot,
x)|_{L^p_{\mathcal{F}}(0, T)}\leq r\mbox{ for a.e.
}x\in(\displaystyle\bigcup_{k\in\dbN} I_k)^c.$ This implies that
$|v(\cdot, \cdot, x)|_{L^p_{\mathcal{F}}(0, T)}\leq 2^n r\mbox{ for
a.e. }x\in\dbR^n.$ Furthermore,  it is easy to check other
conclusions in Lemma \ref{91}. The proof is completed.
\endpf

\begin{remark}\label{122}
Similar to Lemma \ref{91}, we can get the following result: if $u\in
L^{\infty}_{\mathcal{F}}(0, T; L^1(\dbR^n))$, then for any $r>0$,
there exist the functions $v(\cdot),\ w_k(\cdot)\in
L^{\infty}_{\mathcal{F}}(0, T; L^1(\dbR^n))$ $(k=1, 2, \cdots)$ such
that $u(t, \omega, x)=v(t, \omega,
x)+\sum\limits_{k=1}^{\infty}w_k(t, \omega, x)$, where $v(\cdot)$
and $w_k(\cdot)$ satisfy the following conditions:
\begin{eqnarray*}
&&\supp w_k(t, \omega, \cdot)\subseteq \overline{I_k},\mbox{ for
a.e. }(t, \omega)\in(0, T)\times\Omega,\mbox{ where }I_k \mbox{ are
disjoint cubes in }\dbR^n; \\
&&r\sum\limits_{k=1}^{\infty}|I_k|\leq
|u(\cdot)|_{L^{\infty}_{\mathcal{F}}(0, T; L^1(\dbR^n))};\\
&&\int_{\dbR^n}w_k(t, \omega, x)dx=0,\quad |v(t, \omega, x)|\leq
2^n r,\mbox{ for a.e. }(t, \omega, x)\in(0, T)\times\Omega\times\dbR^n;\\
&&|v(t, \omega,
\cdot)|_{L^1(\dbR^n)}+\sum\limits_{k=1}^{\infty}|w_k(t, \omega,
\cdot)|_{L^1(\dbR^n)}\leq 3|u(t, \omega, \cdot)|_{L^1(\dbR^n)}.
\end{eqnarray*}
\end{remark}

Based on the above lemmas, first of all, we give a boundedness
result for a class of SPDOs, whose symbol $a=a(t, \omega, \xi)$ is
independent of the variable $x$. For the functions $\psi^*$ and
$\varphi^*$ given in Lemma \ref{90}, write $a_{-1, *}(t, \omega,
\xi)=\psi^*(\xi)a(t, \omega, \xi)$ and $a_{j, *}(t, \omega,
\xi)=\varphi^*(2^{-j}\xi)a(t, \omega, \xi)$ $(j\in\dbN)$. Then we
have the following result.
\begin{lemma}\label{92}
Suppose that $a=a(t, \omega, \xi)\in S^0_{p}$. Then there exists a
nonnegative function $M_n(\cdot, \cdot)\in L^{p}_{\mathcal{F}}(0,
T)$ such that for any $j\in\dbN\cup\{-1\}$,
$$
\displaystyle\int_{\dbR^n}|\breve{a}_{j, *}(t, \omega, x)|dx\leq
M_n(t, \omega), \mbox{ a.e. }(t, \omega)\in (0, T)\times\Omega,
$$
 where $\breve{a}_{j, *}(t, \omega, x)$ denotes the
Fourier inversion transform of $a_{j, *}(t, \omega, \xi)$ with
respect to the variable $\xi$. Moreover, for the associated SPDO
$A_{j,
*}$ determined by the symbol $a_{j, *}$ and any $p_*\in[1, \infty],$
$A_{j, *}: L^{q}_{\mathcal{F}}(0, T; L^{p_*}(\dbR^n))\rightarrow
L^{q^*}_{\mathcal{F}}(0, T; L^{p_*}(\dbR^n))$ is bounded.
\end{lemma}
The proof of  Lemma \ref{92} is similar to that in the deterministic
case. Therefore, we omit it.

\medskip

Next, we present another useful lemma. Notice that by Lemma
\ref{91}, for any $u\in L^1(\dbR^n; L^p_{\mathcal{F}}(0, T))$, there
exist the functions $v$ and $w_k$ $(k\in\dbN)$ satisfying all the
conditions mentioned in Lemma \ref{91}. Actually, for any $u\in
L^1(\dbR^n; L^p_{\mathcal{F}}(0, T))$ and the associated $v$, the
following conclusion also holds.
\begin{lemma}\label{94}
Suppose that $a=a(t, \omega, \xi)\in S^0_\infty$, $u\in L^1(\dbR^n;
L^p_{\mathcal{F}}(0, T))$ and $u(t, \omega, \cdot)\in
L^1(\dbR^n)\cap L^2(\dbR^n)$, for a.e. $(t, \omega)\in (0,
T)\times\Omega$. Then for any $r>0$ and a.e. $(t, \omega)\in (0,
T)\times \Omega$,
\begin{eqnarray*}
&&r\cdot \left| \{\ x\in\dbR^n; \ |(Au)(t, \omega,
x)|>r\}\right|\\[2mm]
&&\leq C(n, a)\left(|u|_{L^1(\dbR^n; L^p_{\mathcal{F}}(0, T))}+|u(t,
\omega, \cdot)|_{L^1(\dbR^n)}+r^{-1}|v(t, \omega,
\cdot)|^2_{L^2(\dbR^n)}\right),
\end{eqnarray*}
where $v$ is the function associated to $u$ in Lemma \ref{91}.
\end{lemma}
\noindent{\bf Proof. } For any $u\in L^1(\dbR^n;
L^p_{\mathcal{F}}(0, T))$ satisfying $u(t, \omega, \cdot)\in
L^1(\dbR^n)\cap L^2(\dbR^n)$, for a.e. $(t, \omega)\in (0,
T)\times\Omega$, suppose that $v$ and $w_k$ $(k\in\dbN)$ are the
functions mentioned in  Lemma \ref{91}. Then, for any $r>0$, it
follows that
\begin{eqnarray}\label{95}
\begin{array}{rl}
&|\{x\in \dbR^n; |(Au)(t, \omega, x)|>r\}|\\[2mm]
&\leq \left|\{x\in\dbR^n; |(Av)(t, \omega,
x)|>\frac{r}{2}\}\right|+\left|\{x\in\dbR^n;
|(A\sum\limits_{k=1}^{\infty}w_k)(t, \omega,
x)|>\frac{r}{2}\}\right|.
\end{array}
\end{eqnarray}
In the following, we estimate two terms in the right side of
(\ref{95}), respectively. By Theorem \ref{80},
\begin{eqnarray}\label{99}
\begin{array}{rl}
&\displaystyle\frac{r^2}{4}|\{x\in\dbR^n; |(Av)(t, \omega,
x)|>\frac{r}{2}\}|\leq \displaystyle\int_{\dbR^n} |(Av)(t, \omega,
x)|^2dx\leq C(a)\displaystyle\int_{\dbR^n} |v(t, \omega, x)|^2dx.
\end{array}
\end{eqnarray}
On the other hand, denote by $I_{k, *}$ the cube, with the same
center as $I_k$ and the length of  side twice than $I_k$. Let
$I^*=\displaystyle\bigcup_{k\in\dbN} I_{k, *}$. Then, for any
$k\in\dbN$, if we can prove that
\begin{equation}\label{96}
\int_{(I_{k, *})^c} |(Aw_k)(t, \omega, x)|dx\leq
C\int_{\dbR^n}|w_k(t, \omega, x)|dx,
\end{equation}
it follows that
\begin{eqnarray}\label{97}
\begin{array}{rl}
&\displaystyle\frac{r}{2}\left|\{x \in (I^*)^c;
|\sum\limits_{k=1}^{\infty} (A w_k)(t, \omega, x)|>\frac{r}{2}
\}\right|\leq \displaystyle\int_{(I^*)^c}
|\sum\limits_{k=1}^{\infty}(Aw_k)(t, \omega,
x)|dx\\[2mm]
&\leq \sum\limits_{k=1}^{\infty}\displaystyle\int_{(I_{k, *})^c}
|(Aw_k)(t, \omega, x)|dx\leq
C\sum\limits_{k=1}^{\infty}\displaystyle\int_{\dbR^n} |w_k(t,
\omega, x)|dx\leq C|u(t, \omega, \cdot)|_{L^1(\dbR^n)};\\[4mm]
&r|I^*|\leq r\sum\limits_{k=1}^{\infty}|I_{k,
*}|=2^nr\sum\limits_{k=1}^{\infty}|I_k| \leq 2^n
|u(\cdot)|_{L^1(\dbR^n; L^p_{\mathcal{F}}(0, T))}.
\end{array}
\end{eqnarray}
By (\ref{95}), (\ref{99}) and (\ref{97}), we obtain the desired
result. The remainder is devoted to the proof of (\ref{96}). We
denote by $\ell_*$ the length of side of $I_k$ and write $h_j(t,
\omega, x)=\breve{a}_{j, *}(t, \omega, \xi)$, where $a_{j, *}$ are
the functions mentioned in Lemma \ref{92}. Then we see that
\begin{eqnarray}\label{60}
\begin{array}{rl}
&\displaystyle\int_{(I_{k, *})^c}\left|(Aw_k)(t, \omega,
x)\right|dx\leq \sum\limits_{j=-1}^{\infty}\displaystyle\int_{(I_{k,
*})^c}
\left|\int_{\dbR^n} h_j(t, \omega, x-y)w_k(t, \omega, y)dy\right|dx\\
&\leq\sum\limits_{2^j\ell_*\geq 1}\displaystyle\int_{(I_{k,
*})^c}\left| \int_{\dbR^n} h_j(t, \omega, x-y)w_k(t, \omega,
y)dy\right|dx\\
&\quad+\sum\limits_{2^j\ell_*<1}\displaystyle\int_{(I_{k, *})^c}
\left| \int_{\dbR^n} [h_j(t, \omega, x-y)-h_j(t, \omega, x)]w_k(t,
\omega, y)dy\right|dx.
\end{array}
\end{eqnarray}
Similar to  the deterministic case, we estimate two terms in the
right side of (\ref{60}), respectively, and then get (\ref{96}).
This finishes the proof.\endpf

\begin{remark}\label{123}
Similar to Lemma \ref{94}, we can prove the following result:
suppose that $a=a(t, \omega, \xi)\in S^0_\infty$ and $u\in
L^\infty_{\mathcal{F}}(0, T; L^1(\dbR^n))$. Then for any $r>0$ and
a.e. $(t, \omega)\in (0, T)\times \Omega$,
\begin{eqnarray*}
&&r\cdot \left| \{\ x\in\dbR^n; \ |(Au)(t, \omega,
x)|>r\}\right|\leq C(n, a)|u|_{L^\infty_{\mathcal{F}}(0, T;
L^1(\dbR^n))}.
\end{eqnarray*}
\end{remark}

\medskip

Now, we give the $L^p$-estimates $(p>1,\ p\neq2)$ of a class of the
SPDOs of order $(0, \infty)$.
\begin{theorem}\label{t62}
Suppose that $a=a(t, \omega, \xi)\in S^0_\infty$ and $p>1$. Then for
the associated SPDO $A$, the following conclusions hold:\\[2mm]
\noindent $(1)$ if $1<p<2$, $A: L^p(\dbR^n; L^{p'}_{\mathcal{F}}(0,
T))\rightarrow L^p(\dbR^n; L^{p}_{\mathcal{F}}(0, T))$ is
bounded;\\[2mm]
\noindent $(2)$ if $p>2$, $A: L^p(\dbR^n; L^{p}_{\mathcal{F}}(0,
T))\rightarrow L^p(\dbR^n; L^{p'}_{\mathcal{F}}(0, T))$ is bounded.
\end{theorem}

\noindent {\bf Proof. } For any $r>0$, define
$$
u_r(t, \omega, x)=\left\{
\begin{array}{rll}
&u(t, \omega, x) &|u(\cdot, \cdot, x)|_{L^p_{\mathcal{F}}(0, T)}\geq r,\\
&0 &|u(\cdot, \cdot, x)|_{L^p_{\mathcal{F}}(0, T)}<r,
\end{array}
\right.
$$$$
 u^r(t, \omega, x)=\left\{
\begin{array}{rll}
&0 &|u(\cdot, \cdot, x)|_{L^p_{\mathcal{F}}(0, T)}\geq r,\\
&u(t, \omega, x) &|u(\cdot, \cdot, x)|_{L^p_{\mathcal{F}}(0, T)}<r.
\end{array}
\right.
$$
Then, for $1<p<2$, by Lemma \ref{94}, we have that
\begin{eqnarray}\label{l1}
\begin{array}{rl}
&|(Au_r)(t, \omega,
\cdot)|^p_{L^p(\dbR^n)}=\displaystyle\int^\infty_0 p
r^{p-1}|\{x\in\dbR^n; |(Au_r)(t, \omega, x)|>r\}|dr\\[2mm]
&\leq C(n, a, p)\displaystyle\int^\infty_0
r^{p-2}\left[|u_r|_{L^1(\dbR^n; L^p_{\mathcal{F}}(0, T))}+|u_r(t,
\omega, \cdot)|_{L^1(\dbR^n)}+r^{-1}|v_r(t, \omega,
\cdot)|^2_{L^2(\dbR^n)}\right] dr,
\end{array}
\end{eqnarray}
where $v_r$ is the function associated to $u_r$ in Lemma \ref{94}.
Next, we estimate every term in the right side of (\ref{l1})
respectively.
\begin{eqnarray}\label{115}
\begin{array}{rl}
&\dbE\displaystyle\int^T_0\int^\infty_0 r^{p-2}|u_r|_{L^1(\dbR^n;
L^p_{\mathcal{F}}(0, T))}drdt=\dbE\displaystyle\int^T_0\int^\infty_0
\int_{\dbR^n}r^{p-2}|u_r(\cdot, \cdot, x)|_{L^p_{\mathcal{F}}(0,
T)}dxdrdt\\[2mm]
&\leq
C(T)\displaystyle\int_{\dbR^n}\int_0^{|u|_{L^p_{\mathcal{F}}(0, T)}}
r^{p-2}|u(\cdot, \cdot, x)|_{L^p_{\mathcal{F}}(0, T)}drdx\leq C(T,
p)\displaystyle\int_{\dbR^n}|u(\cdot, \cdot,
x)|_{L^p_{\mathcal{F}}(0, T)}^pdx;
\end{array}
\end{eqnarray}
\begin{eqnarray}\label{116}
\begin{array}{rl}
&\dbE\displaystyle\int^T_0\int^\infty_0 r^{p-2}|u_r(t, \omega,
\cdot)|_{L^1(\dbR^n)}drdt=\dbE\displaystyle\int^T_0\int^\infty_0
\int_{\dbR^n}r^{p-2}|u_r(t, \omega, x)|dxdrdt\\[2mm]
&\leq C(T,
p)\displaystyle\int^\infty_0\int_{\dbR^n}r^{p-2}|u_r(\cdot, \cdot,
x)|_{L^p_{\mathcal{F}}(0, T)}dxdr\\[2mm]
&=C(T, p)\displaystyle\int_{\dbR^n}\int_0^{|u|_{L^p_{\mathcal{F}}(0,
T)}}r^{p-2}|u(\cdot, \cdot, x)|_{L^p_{\mathcal{F}}(0, T)}drdx=C(T,
p)\int_{\dbR^n}|u(\cdot, \cdot, x)|_{L^p_{\mathcal{F}}(0, T)}^pdx;
\end{array}
\end{eqnarray}
\begin{eqnarray}\label{117}
\begin{array}{rl}
&\dbE\displaystyle\int^T_0\int^\infty_0 r^{p-3}|v_r(t, \omega,
\cdot)|_{L^2(\dbR^n)}^2drdt=\dbE\displaystyle\int^T_0\int^\infty_0\int_{\dbR^n}
r^{p-3}|v_r(t, \omega, x)|^2dxdrdt\\[2mm]
&\leq\displaystyle\int^\infty_0\int_{\dbR^n}r^{p-3}|v_r(\cdot,
\cdot, x)|_{L^p_{\mathcal{F}}(0, T)}|v_r(\cdot, \cdot,
x)|_{L^{p'}_{\mathcal{F}}(0, T)}dxdr\\[2mm]
&\leq C(n)\displaystyle\int^\infty_0\int_{\dbR^n}r^{p-2}|v_r(\cdot,
\cdot, x)|_{L^{p'}_{\mathcal{F}}(0, T)}dxdr\\[2mm]
&= C(n)\displaystyle\int^\infty_0 [ \sum\limits_{k\in\dbN}\int_{I_k}
r^{p-2}|v_r(\cdot, \cdot, x)|_{L^{p'}_{\mathcal{F}}(0, T)}dx+
\int_{(\displaystyle\bigcup_{k\in\dbN} I_k)^c}r^{p-2}|v_r(\cdot,
\cdot, x)|_{L^{p'}_{\mathcal{F}}(0, T)}dx ] dr\\[2mm]
&\leq C(n)\displaystyle\int^\infty_0 [
\sum\limits_{k\in\dbN}\int_{I_k} r^{p-2}|u_r(\cdot, \cdot,
x)|_{L^{p'}_{\mathcal{F}}(0, T)}dx+
\int_{(\displaystyle\bigcup_{k\in\dbN} I_k)^c}r^{p-2}|u_r(\cdot,
\cdot, x)|_{L^{p'}_{\mathcal{F}}(0, T)}dx ] dr\\[2mm]
&=C(n)\displaystyle\int^\infty_0\int_{\dbR^n}r^{p-2}|u_r(\cdot,
\cdot, x)|_{L^{p'}_{\mathcal{F}}(0, T)}dxdr\\[2mm]
&=C(n)\displaystyle\int_{\dbR^n}\int_0^{|u|_{L^p_{\mathcal{F}}(0,
T)}}r^{p-2}|u(\cdot, \cdot, x)|_{L^{p'}_{\mathcal{F}}(0,
T)}drdx\\[2mm]
&=C(n, p)\displaystyle\int_{\dbR^n}|u(\cdot, \cdot,
x)|_{L^{p}_{\mathcal{F}}(0, T)}^{p-1}|u(\cdot, \cdot,
x)|_{L^{p'}_{\mathcal{F}}(0, T)}dx\leq C(n, p,
T)\displaystyle\int_{\dbR^n}|u(\cdot, \cdot,
x)|_{L^{p'}_{\mathcal{F}}(0, T)}^{p}dx.
\end{array}
\end{eqnarray}
By (\ref{l1})-(\ref{117}), we see that
\begin{equation}\label{118}
\dbE\displaystyle\int^T_0|(Au_r)(t, \omega,
\cdot)|^p_{L^p(\dbR^n)}dt\leq C(n, a, p,
T)\displaystyle\int_{\dbR^n}|u(\cdot, \cdot,
x)|_{L^{p'}_{\mathcal{F}}(0, T)}^{p}dx.
\end{equation}

On the other hand,
\begin{eqnarray*}\label{l2}
\begin{array}{rl}
&|(Au^r)(t, \omega,
\cdot)|^p_{L^p(\dbR^n)}=\displaystyle\int^\infty_0 p
r^{p-1}|\{x\in\dbR^n; |(Au^r)(t, \omega, x)|>r\}|dr\\[2mm]
&\leq\displaystyle\int^\infty_0 p r^{p-3}\int_{\dbR^n}|(Au^r)(t,
\omega, x)|^2dx dr\leq C(a)\displaystyle\int_{0}^{\infty}
\int_{\dbR^n}pr^{p-3}|u^r(t, \omega,
x)|^2dxdr\\[2mm]
&=C(a)\displaystyle\int_{\dbR^n}
\displaystyle\int^\infty_{|u|_{L^p_{\mathcal{F}}(0,
T)}}pr^{p-3}dr|u(t, \omega, x)|^2dx\leq C(a,
p)\int_{\dbR^n}|u(\cdot, \cdot, x)|^{p-2}_{L^p_{\mathcal{F}}(0,
T)}|u(t, \omega, x)|^2dx.
\end{array}
\end{eqnarray*}
Therefore, we obtain that
\begin{eqnarray}\label{121}
\begin{array}{rl}
&\dbE\displaystyle\int^T_0|(Au^r)(t, \omega,
\cdot)|^p_{L^p(\dbR^n)}dt \leq C(a, p)\dbE\displaystyle\int^T_0
\int_{\dbR^n}|u(\cdot, \cdot, x)|^{p-2}_{L^p_{\mathcal{F}}(0,
T)}|u(t, \omega, x)|^2dxdt\\[2mm]
&=C(a, p)\displaystyle\int_{\dbR^n}|u(\cdot, \cdot,
x)|^{p-2}_{L^p_{\mathcal{F}}(0,
T)}\left(\dbE\displaystyle\int^T_0|u(t, \omega, x)|^2dt\right)dx\\[2mm]
&\leq C(a, p)\displaystyle\int_{\dbR^n}|u(\cdot, \cdot,
x)|^{p-1}_{L^p_{\mathcal{F}}(0, T)}|u(\cdot, \cdot,
x)|_{L^{p'}_{\mathcal{F}}(0, T)}dx\leq C(a, p,
T)\displaystyle\int_{\dbR^n}|u(\cdot, \cdot,
x)|^{p}_{L^{p'}_{\mathcal{F}}(0, T)}dx.
\end{array}
\end{eqnarray}
(\ref{118}) and (\ref{121}) imply the desired result (1). By Lemma
2.1 in \cite{lyz} and a duality argument, we can get the result (2)
for $p>2$.
\endpf

\medskip

\begin{remark}
Similar to the proof of Theorem \ref{t62}, by Remark \ref{122} and
Remark \ref{123}, we can prove the following result: if $a=a(t,
\omega, \xi)\in S^0_\infty$, then for $1<p<2$,  the associated SPDO
$A: L^p(\dbR^n; L^{\infty}_{\mathcal{F}}(0, T))\rightarrow
L^{\infty}_{\mathcal{F}}(0, T; L^p(\dbR^n))$ is bounded. However,
for $p>2$, we do not establish the corresponding boundedness result,
because we fail to give a characterization for the dual spaces of
$L^p(\dbR^n; L^{\infty}_{\mathcal{F}}(0, T))$ and
$L^{\infty}_{\mathcal{F}}(0, T; L^p(\dbR^n))$ at this moment.
\end{remark}

\medskip

In the remainder of this subsection, we establish the
$L^p$-estimates of the SPDOs defined on a local domain. First of
all, we give a preliminary.
\begin{corollary}\label{t63}
Suppose that $a\in S^0_\infty(\dbR^n\times \dbR^n\times\dbR^n)$ and
$A$ is the associated SPDO. If there exists a bounded set
$B_0\subseteq\dbR^{2n}$, such that $\supp a(t, \omega, \cdot, \cdot,
\xi)\subseteq B_0$ for a.e. $(t, \omega)\in(0, T)\times\Omega$ and
any $\xi\in\dbR^n$, then the conclusions $(1)$-$(2)$ mentioned in
Theorem \ref{t62} hold.
\end{corollary}
\noindent{\bf Sketch of the proof. } The method of the proof is
similar to that used in the deterministic case. First, we construct
a family of amplitudes $\{a_{\eta, \varsigma}\}_{(\eta,
\varsigma)\in\dbR^{2n}}$, whose associated SPDOs $\{A_{\eta,
\varsigma}\}_{(\eta, \varsigma)\in\dbR^{2n}}$ have good estimates
with respect to parameters $(\eta, \varsigma)$. Indeed, write
$$
\hat{a}(t, \omega, \eta, \varsigma,
\xi)=\displaystyle\int_{\dbR^n}\int_{\dbR^n}e^{-i(x\cdot\eta+y\cdot\varsigma)}
a(t, \omega, x, y, \xi)dxdy.
$$
Then for any multi-indices $\alpha_1$, $\alpha_2$ and $\beta$, it
follows that
\begin{eqnarray*}
&&\eta^{\alpha_1}\varsigma^{\alpha_2}\partial^{\beta}_{\xi}\hat{a}(t,
\omega, \eta, \varsigma,
\xi)=\eta^{\alpha_1}\varsigma^{\alpha_2}\partial^{\beta}_{\xi}
\int_{\dbR^n}\int_{\dbR^n}e^{-i(x\cdot\eta+y\cdot\varsigma)} a(t,
\omega, x, y, \xi)dxdy\\[2mm]
&&=\int_{\dbR^n}\int_{\dbR^n}\frac{1}{(-i)^{|\alpha_1|+|\alpha_2|}}
\partial^{\alpha_1}_x\partial^{\alpha_2}_y
e^{-i(x\cdot\eta+y\cdot\varsigma)} \partial^\beta_\xi a(t, \omega,
x, y, \xi)dxdy\\[2mm]
&&=\int_{\dbR^n}\int_{\dbR^n}\frac{1}{i^{|\alpha_1|+|\alpha_2|}}
e^{-i(x\cdot\eta+y\cdot\varsigma)}
\partial^{\alpha_1}_x\partial^{\alpha_2}_y\partial^\beta_\xi a(t,
\omega, x, y, \xi)dxdy.
\end{eqnarray*}
This implies that for any $k\in\dbN$, a.e. $(t, \omega)\in(0,
T)\times\Omega$ and any $(\eta, \varsigma, \xi)\in\dbR^{3n}$,
$$
|\partial^\beta_\xi\hat{a}(t, \omega, \eta, \varsigma, \xi)|\leq
C(k, \beta, B_0)(1+|\xi|)^{-|\beta|}(1+|\eta|+|\varsigma|)^{-k}.
$$

Next, denote by $A_{\eta, \varsigma}$ the SPDO, whose amplitude is
$\hat{a}$. Based on the result of Theorem \ref{t62}, we obtain that
\begin{eqnarray*}
&&|A_{\eta, \varsigma}|_{\mathcal{L}(L^p(\dbR^n;
L^{p'}_{\mathcal{F}}(0, T)), L^p(\dbR^n; L^p_{\mathcal{F}}(0,
T)))}\leq C(n, k, T, a, B_0,
p)(1+|\eta|+|\varsigma|)^{-k},\mbox{ for }1<p<2; \\[2mm]
&&|A_{\eta, \varsigma}|_{\mathcal{L}(L^p(\dbR^n;
L^{p}_{\mathcal{F}}(0, T)), L^p(\dbR^n; L^{p'}_{\mathcal{F}}(0,
T)))}\leq C(n, k, T, a, B_0, p)(1+|\eta|+|\varsigma|)^{-k},\mbox{
for }p>2.
\end{eqnarray*}
Moreover, notice that
\begin{eqnarray*}
&&(Au)(t, \omega, x)=(2\pi)^{-n}\int_{\dbR^n}\int_{\dbR^n}
e^{i(x-y)\cdot\xi}a(t, \omega, x, y, \xi)u(t, \omega, y)dxdy\\[2mm]
&&=(2\pi)^{-3n}\int_{\dbR^n}\int_{\dbR^n}e^{i(x-y)\cdot\xi}
\left[\int_{\dbR^n}\int_{\dbR^n}e^{i(x\cdot\eta+y\cdot\varsigma)}
\hat{a}(t, \omega, \eta, \varsigma, \xi)d\eta d\varsigma\right]u(t,
\omega, y)dxdy\\[2mm]
&&=(2\pi)^{-2n}\int_{\dbR^n}\int_{\dbR^n}e^{i(x\cdot\eta+y\cdot\varsigma)}
(A_{\eta, \varsigma}u)(t, \omega, x)d\eta d \varsigma.
\end{eqnarray*}
Therefore, for $1<p<2$,
\begin{eqnarray*}
&&|Au|_{L^p(\dbR^n; L^p_{\mathcal{F}}(0, T))}\leq C(n)
\int_{\dbR^n}\int_{\dbR^n} |A_{\eta, \varsigma}u|_{L^p(\dbR^n;
L^p_{\mathcal{F}}(0,
T))}d\eta d\varsigma\\[2mm]
&&\leq C(n, T, a, B_0, p)\int_{\dbR^n}\int_{\dbR^n}
(1+|\eta|+|\varsigma|)^{-2n-2}|u|_{L^p(\dbR^n;
L^{p'}_{\mathcal{F}}(0,
T))}d\eta d\varsigma\\[2mm]
&&\leq C(n, T, a, B_0, p)|u|_{L^p(\dbR^n; L^{p'}_{\mathcal{F}}(0,
T))}.
\end{eqnarray*}
We can get the  result in the case of $p>2$ in the same way.
\endpf

\medskip

Next, we introduce some locally convex topological spaces. For any
$p, q>1$ and any compact set $K\subseteq\dbR^n$,
\begin{eqnarray*}
&&L^p(K; L^q_{\mathcal{F}}(0, T))=\left\{\ u\in L^p(G;
L^q_{\mathcal{F}}(0, T))\ |\ \supp u(t, \omega, \cdot) \subseteq
K,\mbox{ for a.e. }(t, \omega)\in(0,
T)\times\Omega\right\},\\
&&\mbox{which is generated by the norm }|\cdot|_{L^p(G;
L^q_{\mathcal{F}}(0, T))};\\[2mm]
&&L^p_{comp}(G; L^q_{\mathcal{F}}(0, T))=\bigcup_{K\subseteq G\mbox{
is compact}}L^p(K; L^q_{\mathcal{F}}(0, T)),\\
&&\mbox{ which is endowed
with the inductive topology};\\[2mm]
&&L^p_{loc}(G; L^q_{\mathcal{F}}(0, T))=\left\{\ u\ |\ \ |u|_{p, q,
\psi}=|\psi u|_{L^p(G; L^q_{\mathcal{F}}(0, T))}<\infty,\mbox{ for
any }\psi\in
C^\infty_0(G) \right\},\\
&&\mbox{which is generated by a family of semi-norms }\{|\cdot|_{p,
q, \psi}\}_{\psi\in C^\infty_0(G)}.
\end{eqnarray*}

\medskip

Then, based on Corollary \ref{t63},  by the method similar to that
used in the proof of Corollary \ref{89}, we get the following
boundedness result.
\begin{corollary}
Suppose that $A$ is a SPDO and $a$ is its amplitude.\\[2mm]
\noindent $(1)$ If $a\in S^0_\infty(G\times G\times \dbR^n)$, then
for $1<p<2$, $A: L^p_{comp}(G; L^{p'}_{\mathcal{F}}(0,
T))\rightarrow L^p_{loc}(G; L^{p}_{\mathcal{F}}(0, T))$ is
continuous;\\[2mm]
\noindent $(2)$ If $a\in S^0_\infty(G\times G\times \dbR^n)$, then
for $p>2$, $A: L^p_{comp}(G; L^{p}_{\mathcal{F}}(0, T))\rightarrow
L^p_{loc}(G; L^{p'}_{\mathcal{F}}(0, T))$ is continuous;\\[2mm]
\noindent $(3)$ If $A$ is a uniformly properly supported SPDO with
respect to $(t, \omega)$ and $a\in S^0_\infty(G\times G\times
\dbR^n)$, than for $1<p<2$, both $A: L^p_{comp}(G;
L^{p'}_{\mathcal{F}}(0, T))\rightarrow L^p_{comp}(G;
L^{p}_{\mathcal{F}}(0, T))$ and $A: L^p_{loc}(G;
L^{p'}_{\mathcal{F}}(0, T))\rightarrow L^p_{loc}(G;
L^{p}_{\mathcal{F}}(0,
T))$ are continuous;\\[2mm]
\noindent $(4)$ If $A$ is a uniformly properly supported SPDO with
respect to $(t, \omega)$ and $a\in S^0_\infty(G\times G\times
\dbR^n)$, than for $p>2$, both $A: L^p_{comp}(G;
L^{p}_{\mathcal{F}}(0, T))\rightarrow L^p_{comp}(G;
L^{p'}_{\mathcal{F}}(0, T))$ and $A: L^p_{loc}(G;
L^{p}_{\mathcal{F}}(0, T))\rightarrow L^p_{loc}(G;
L^{p'}_{\mathcal{F}}(0, T))$ are continuous.
\end{corollary}

\section{Elliptic stochastic pseudo-differential operators}

In this subsection, we  introduce the notion of elliptic SPDOs and
point out their invertibility. Also, we give the G{\aa}rding
inequality.

First of all, we  define an elliptic operator and its parametrix.
\begin{definition}\label{30}
Suppose that $A$ is a SPDO of order $(\ell, p)$ and $a$ is its
symbol. $A$ is called an elliptic SPDO if  for any compact set $K$,
there exist two positive constants $C_K$ and $R_K$, such that
$$
|a(t, \omega, x, \xi)|\geq C_K(1+|\xi|)^{\ell}, \ \mbox{for a.e. }
(t, \omega)\in (0, T)\times \Omega\mbox{ and any }(x, \xi)\in
K\times\{\ \xi\in \dbR^n;\ |\xi|\geq R_K\}.
$$
\end{definition}
\begin{definition}\label{31}
A left parametrix $Q_1$ (or a right parametrix $Q_2$) for a SPDO $A$
of order $(\ell, p)$ is a SPDO, which is a single-sided inverse for
$A$ modulo smoothing operators satisfying
$$Q_1A-I\in \mathcal{L}^{-\infty}_\infty\quad\quad (\mbox{ or }
AQ_2-I\in \mathcal{L}^{-\infty}_\infty).$$
\end{definition}

Next, we present invertibility of elliptic SPDOs.
\begin{theorem}\label{32}
Suppose that $A$ is an elliptic SPDO of order $(\ell, \infty)$. Then
there exist two SPDOs $Q_1,\ Q_2\in
\mathcal{L}^{-\ell}_{\infty}(G)$, which are the left parametrix and
right parametrix for $A$, respectively.
\end{theorem}

\medskip

Also, we give a lemma, which is a preliminary for the proof of the
G{\aa}rding inequality.
\begin{lemma}\label{61}
Suppose that $A$ is a SPDO of order $(0, \infty)$ and $a$ is its
symbol. If for a.e.  $(t, \omega)\in (0, T)\times\Omega$ and any
$(x, \xi)\in \dbR^{2n}$, Re $a(t, \omega, x, \xi)\geq C^*$ for a
positive constant $C^*$, then there exists a SPDO $B\in
\mathcal{L}^0_\infty$, such that
$$
\displaystyle\frac{A+A^*}{2}-B^*B\in \mathcal{L}^{-\infty}_\infty.
$$
\end{lemma}
The proofs of Theorem \ref{32} and Lemma \ref{61} are similar to
those in the deterministic case. Here we omit them.

\medskip

Now, we prove the following  G{\aa}rding inequality.
\begin{theorem}\label{101}
Suppose that $A$ is a SPDO of order $(\ell, \infty)$ and $a$ is its
symbol. If there exist two positive constants $\delta^*$ and $R^*$,
such that Re $a(t, \omega, x, \xi)\geq \delta^*|\xi|^{\ell}$ for
a.e. $(t, \omega)\in(0, T)\times\Omega$ and any $(x,
\xi)\in\dbR^{n}\times\{\ \xi\in\dbR^n;\ |\xi|\geq R^*\}$, then for
any $\varepsilon>0$, $r\in\dbR$ and $u\in L^2_\mathcal{F}(0, T;
H^r(\dbR^n)\cap H^{\frac{\ell}{2}}(\dbR^n))$,
\begin{eqnarray*}
\begin{array}{rl}
&\dbE\displaystyle\int^T_0\mbox{Re }\left(Au(t, \omega, \cdot), u(t,
\omega, \cdot)\right)_{L^2(\dbR^n)}dt\\[2mm]
&\geq (\delta^*-\varepsilon)\dbE\displaystyle\int^T_0 |u(t, \omega,
\cdot)|^2_{H^{\frac{\ell}{2}}(\dbR^n)}dt-C\dbE\displaystyle\int^T_0|u(t,
\omega, \cdot)|^2_{H^r(\dbR^n)}dt.
\end{array}
\end{eqnarray*}
\end{theorem}

\noindent{\bf Sketch of the proof. } If $\ell=0$, without of
generality, we  assume that Re $a(t, \omega, x, \xi)\geq
\delta^*|\xi|^{\ell}$ for a.e.
 $(t, \omega)\in(0, T)\times\Omega$ and any $(x, \xi)\in\dbR^{2n}$.
 Write
$a_\varepsilon=a-(\delta^*-\varepsilon)$ and let $A_\varepsilon$ be
the associated SPDO. Then, we see that Re $a_\varepsilon\geq
\varepsilon$, and by Lemma \ref{61},  there exists a SPDO $B\in
\mathcal{L}^0_\infty$ such that
$\displaystyle\frac{A_\varepsilon+A_\varepsilon^*}{2}-B^*B=R_1\in
\mathcal{L}^{-\infty}_{\infty}$. Therefore, $
\displaystyle\frac{A+A^*}{2}-(\delta^*-\varepsilon)I-B^*B=R_1. $ By
the $L^2$-estimates, this leads to that
\begin{eqnarray*}
&&\dbE\int^T_0 \mbox{Re }((Au)(t, \omega, \cdot), u(t, \omega,
\cdot))_{L^2(\dbR^n)}dt\\[2mm]
&&=\dbE\int^T_0 \mbox{Re }((\delta^*-\varepsilon)u(t, \omega,
\cdot), u(t, \omega, \cdot))_{L^2(\dbR^n)}dt+\dbE\int^T_0 \mbox{Re
}((Bu)(t, \omega, \cdot), (Bu)(t, \omega,
\cdot))_{L^2(\dbR^n)}dt\\[2mm]
&&\quad+\dbE\int^T_0 \mbox{Re }((R_1u)(t,
\omega, \cdot), u(t, \omega, \cdot))_{L^2(\dbR^n)}dt\\[2mm]
&&\geq (\delta^*-2\varepsilon)\dbE\int^T_0 |u(t, \omega,
\cdot)|_{L^2(\dbR^n)}^{2}dt-C(\varepsilon)\dbE\displaystyle\int^T_0|u(t,
\omega, \cdot)|^2_{H^r(\dbR^n)}dt.
\end{eqnarray*}
We can also get the desired results for $a\in S^{\ell}_\infty$
$(\ell\neq0)$ by the technique used in Corollary \ref{88}. This
finishes the proof.
\endpf

\section{Calder\'{o}n-type theorem of stochastic pseudo-differential operators}

In this section, as an application of the theory of SPDOs, we
establish a Calder\'on-type uniqueness theorem on the Cauchy problem
of SPDEs. In order to present the main idea, we focus on the case of
at most double characteristics. The key point of the proof is to
establish a new Carleman-type estimate for the SPDOs of order $(1,
\infty)$.

\subsection{Statement of the main result}

Let $U$ be a neighborhood of the origin in $\dbR^n$. Set
$\mathcal{X}_\infty=\displaystyle\bigcap_{j\in\dbN}
L_{\mathcal{F}}^\infty(\O;C^1([0,T]; C^j(\overline{U})))$ and
$\mathcal{X}_m=\displaystyle\bigcap_{k_1+k_2=m, \atop
k_1<m}L_{\mathcal{F}}^2(\O;C^{k_1}([0,T];H^{k_2}(U)))$. Consider the
Cauchy problem for the following linear SPDE of order $m$:
 \begin{eqnarray}\label{x2}
\left\{
\begin{array}{ll}
\frac{1}{i}d D^{m-1}_tu=\displaystyle\sum^{m-1}_{k=0}
\displaystyle\sum_{|\alpha|=m-k}a_\alpha(t,\o, x
)D^\alpha_xD^k_tu dt\\[3mm]
\qq\qq\qq+ \sum\limits_{|\b|<m}\left[b_\b(t,\o, x )D^\b_{t,x}u
dt+c_\b(t,\o, x )D^\b_{t,x}u dw(t)\right]  \quad\quad\mbox{in
}(0,T)\times\Omega\times
U,\\[3mm]
u(0)=D_tu(0)=\cdots=D_t^{m-1}u(0)=0
\qq\qq\qq\qq\qq\quad\quad\quad\quad\mbox{in }\Omega\times U,
\end{array}\right.
\end{eqnarray}
where $D_t=\frac{1}{i}\frac{\partial}{\partial t}$,
$D_{x_k}=\frac{1}{i}\frac{\partial}{\partial {x_k}}$,  $\alpha$ and
$\beta$ denote two multi-indices, and  $a_\a$, $b_\b$, $c_\b$ $\in
\mathcal{X}_\infty$.

Write $\displaystyle p_m(t,\o, x, \lambda,
\xi)=\lambda^m-\sum^{m-1}_{k=0} \sum_{|\alpha|=m-k}a_\alpha(t,\o,
x)\xi^\alpha\lambda^k$ and denote by $\{\lambda_k(t,\o, x, \xi);
k=1, \cdots, m\}$ the characteristic roots of $p_m(t,\o, x, \lambda,
\xi)$ for a.e. $(t, \o)\in (0,T)\times\Omega$ and any $(x, \xi)\in
U\times \dbR^n$, i.e., $p_m(t,\o, x, \lambda_k(t,\o, x, \xi),
\xi)=0$.  Also, for a.e. $(t,\o)\in (0,T)\times\Omega$ and any $(x,
\xi)\in U\times\dbR^n$ satisfying $|\xi|=1$, we introduce the
following hypotheses:

\vspace{1mm}

\begin{enumerate}
\item[{\bf (H1)}] all real roots are simple and the multiplicity of
all complex roots is at most two;\vspace{1mm}
\item[{\bf (H$1'$)}] all roots $\lambda_k(t,\o, x,  \xi)$ $(k=1, \cdots, m)$ are
simple;\vspace{1mm}
\item[{\bf (H2)}] there exists a positive constant $\varepsilon$, (which is independent of $t,
\omega$, $x$ and $\xi$) such that for any complex root
$\lambda_k(t,\o, x, \xi)$, $|\mbox{Im} \lambda_k(t,\o, x, \xi)|\geq
\varepsilon$;\vspace{1mm}
\item[{\bf (H3)}] If a root $\lambda_k(t,\o, x, \xi)$ is real
(complex) at a point, it remains real (complex) at every
point;\vspace{1mm}
\item[{\bf (H4)}] The algebraic multiplicity of all complex roots is
constant with respect to every variable, and the geometric
multiplicity of all complex roots is constant with respect to
$\omega$.
\end{enumerate}

\vspace{1mm}

Then, the main results in this section are stated as follows. The
first one is a uniqueness result for equation (\ref{x2}) in the case
of single characteristics.

\begin{theorem}\label{x1}
Suppose that the hypotheses $(H1'), (H2)$ and $(H3)$ hold. If $\ds
u\in \mathcal{X}_m$ is a strong solution of equation $(\ref{x2})$.
Then there exist a neighborhood $V(\subset U)$ of the origin in
$\dbR^n$ and a sufficiently small $T'>0$ such that  $u$ vanishes in
$(0, T')\times\Omega\times V$.
\end{theorem}

\vspace{1mm}

The other one generalizes the result of Theorem \ref{x1} to the case
of at most double characteristics.

\begin{theorem}\label{x3}
Suppose that the hypotheses $(H1), (H2)$, $(H3)$ and $(H4)$ hold. If
$\ds u\in \mathcal{X}_m$ is a strong solution of equation
$(\ref{x2})$. Then there exist a neighborhood $V(\subset U)$ of the
origin in $\dbR^n$ and a sufficiently small $T'>0$ such that  $u$
vanishes in $(0, T')\times\Omega\times V$.
\end{theorem}

\medskip

We shall give the proofs of Theorem \ref{x1} and Theorem \ref{x3} in
the next two subsections respectively.

\subsection{Proof of Theorem \ref{x1}}

In order to get a Calder\'{o}n-type uniqueness result, we first
point out that it suffices to establish a suitable estimate for a
strong solution of a SPDE of order $m$.
 Next,
by introducing a pseudo-differential operator,   we reduce the
desired estimate to a new Carleman-type estimate for a SPDO of order
$(1, \infty)$. Finally, we give the proof of the Carleman estimate
(see Lemma \ref{6}), based on the theory of SPDOs mentioned in
Sections 2-4.

To begin with, we introduce a smooth real function $\zeta$
satisfying that $\zeta(t)=1$ in $[0, 2T/3]$ and $\zeta(t)=0$ in $[T,
\infty)$. Write $B_r=\{x\in \dbR^n; |x|<r\}$. Then, we have the
following result.
\begin{proposition}\label{3}
Under the assumptions $(H1'),\ (H2)$ and $(H3)$, if
$u\in\mathcal{X}_m$ is a strong solution of equation $(\ref{x1})$
satisfying that $\supp u(t, \omega, \cdot)\subseteq B_r$ for a.e.
$(t, \omega)\in (0, T)\times\Omega$, then there exists a constant
$C$, (which is independent of $u$ and $\mu$) such that for
sufficiently small positive constants $r$, $T$ and $\mu^{-1}$,
\begin{eqnarray}\label{7}
\begin{array}{rl}
&\dbE\displaystyle\int^T_0
e^{\mu(t-T)^2}\displaystyle\sum_{|\alpha|<m}|D^\alpha_{t, x} (\zeta
u)(t)|^2_{L^2(B_r)}dt \\[2mm]
&\leq
C(T+\mu^{-1})\dbE\displaystyle\int^T_0e^{\mu(t-T)^2}|\sum\limits_{|\b|<m}f_\b(t,\o,
\cdot )D^\b_{t,x}u(t) |^2_{L^2(B_r)}dt,
\end{array}
\end{eqnarray}
where $f_\beta\in \mathcal{X}_\infty$ $(0\leq|\beta|<m)$ depends
only on $\{a_\alpha\}_{1\leq|\alpha|\leq m},
\{b_\beta\}_{|\beta|<m}$ and $\{c_\beta\}_{|\beta|<m}$, and
$f_\beta(t, \omega, x)=0$ for $(t, w, x)\in[0,
2T/3]\times\Omega\times U$.
\end{proposition}

\begin{remark}\label{20}
Let us show how to deduce Theorem $\ref{x1}$ from Proposition
$\ref{3}$. In fact, let $u$ be a strong solution of equation
$(\ref{x2})$. Without loss of generality, we suppose that $\supp
u\subseteq \dbR^+\times\Omega\times\dbR^n$ by the given initial
condition. In order to construct a function with compact support
with respect to the variable $x$, we make the following Holmgren
transformation:
$$(t, \omega, x)\rightarrow (t', \omega', x'), \ x'=x, \ \omega'=\omega, \ t'=t+\delta'|x|^2$$
where $\delta'$ is a sufficiently small positive constant. Then
after the coordinate transformation, all conditions in Theorem
$\ref{x1}$ still hold, and for sufficiently small positive constants
$T$ and $r$, $\supp u(t, \omega, \cdot)\subseteq B_r$ for a.e. $(t,
\omega)\in (0, T)\times \Omega$. Therefore, by $(\ref{7})$, it
follows that
\begin{eqnarray*}
\begin{array}{rl}
&e^{\frac{\mu T^2}{4}}\dbE\displaystyle\int^\frac{T}{2}_0
|u(t)|_{L^2(B_r)}^2dt
\leq\dbE\displaystyle\int^\frac{T}{2}_0e^{\mu(t-T)^2}|u(t)|_{L^2(B_r)}^2dt\\
&\leq C(T, u, r, f_\beta,
m)\mu^{-1}\displaystyle\int^T_\frac{2T}{3}e^{\mu(t-T)^2}dt\leq C(T,
u, r, f_\beta, m)\mu^{-1}e^{\frac{\mu}{9}T^2}.
\end{array}
\end{eqnarray*}
Letting $\mu\rightarrow \infty$ in the above inequality, we obtain
that $u=0$ in $(0, T/2)\times \Omega\times B_r$, which implies
Theorem $\ref{x1}$.
\end{remark}

In the following, we divide the proof of Proposition \ref{3} into
four parts.

\medskip

\noindent{\bf Step 1 }We transform a SPDE of order $m$ to a
stochastic pseudo-differential system of order 1. Write $A_k(t,\o,
x, D)=\displaystyle\sum_{|\alpha|=m-k}a_\alpha(t,\o, x)D^\alpha_x$,
and then the symbol $a_k(t,\o, x, \xi)$ of $A_k$ is
$\displaystyle\sum_{|\alpha|=m-k}a_\alpha(t,\o, x)\xi^\alpha$.
 Therefore, if we denote by $u$ a strong solution of
equation (\ref{x2}) and let
$$Y=(\Lambda^{m-1}(\zeta u), D_t\Lambda^{m-2}(\zeta u), \cdots,
D_t^{m-1}(\zeta u))^\top,$$  it is easy to see that
 $$
  \frac{1}{i}dY=AYdt +fdt+Fdw(t),
 $$
where
$$
 \ba{ll}
 A=\left(
\begin{array}{lccccr}
&0 &\Lambda &0 &\cdots &0\\
&0 &0 &\Lambda &\cdots &0\\
&\vdots &\vdots &\vdots &\vdots &\vdots\\
&0 &0  &0 &\cdots &\Lambda\\
&A_0\Lambda^{1-m} &A_{1}\Lambda^{2-m} &A_{2}\Lambda^{3-m} &\cdots
&A_{m-1}
\end{array}
\right), \\[4mm]\ns\ds f=(0, 0, \cdots, \sum_{|\b|<m}f^1_\b(t,\o, x
)D^\b_{t,x}u+\sum_{|\b|<m}b_\b(t,\o, x )D^\b_{t,x}(\zeta
u))^\top,\\[3mm]\ns\ds
 F=(0, 0, \cdots,
\sum_{|\b|<m}f^2_\b(t,\o, x )D^\b_{t,x}u+\sum_{|\b|<m} c_\b(t,\o, x
)D^\b_{t,x}(\zeta u) )^\top,
 \ea
$$
and $f^1_\beta,\ f^2_\beta \in\mathcal{X}_\infty$ $(|\beta|<m)$
depend only on $m$, $\{a_\alpha\}_{1\leq |\alpha|\leq m},\
\{b_\beta\}_{|\beta|< m},\ \{c_\beta\}_{|\beta|< m}$ and $\zeta$.
Moreover, $f^j_\beta(t, \omega, x)=0$ for $(t, w, x)\in[0,
2T/3]\times\Omega\times U$ $(j=1, 2)$. Furthermore, if we denote by
$A_0$ a SPDO, with the symbol
$$
\sigma(A_0)=\left(
\begin{array}{lccccr}
&0 &|\xi| &0 &\cdots &0\\
&0 &0 &|\xi| &\cdots &0\\
&\vdots &\vdots &\vdots &\vdots &\vdots\\
&0 &0  &0 &\cdots &|\xi|\\
&a_0|\xi|^{1-m} &a_{1}|\xi|^{2-m} &a_{2}|\xi|^{3-m} &\cdots &a_{m-1}
\end{array}
\right),
$$
then $A=A_0+B$, where $B\in \mathcal{L}^{0}_\infty(U )$.

\medskip

\noindent {\bf Step 2} We make the diagonalization for  operator
$A_0$. By the assumption $(H1')$, for a.e. $(t,\o)\in(0, T)\times
\Omega$ and any $(x, \xi)\in U\times\dbR^n$ satisfying $|\xi|=1$,
there exists an invertible matrix $r^*(t,\o, x, \xi)$ such that
$j^*=r^*\cdot\sigma(A_0)\cdot {r^{*}}^{-1}$ is a diagonal matrix.
Since $\sigma(A_0)$ is homogeneous of order 1 with respect to $\xi$,
$j^*$ is still a diagonal matrix after $r^*$ is extended
homogeneously of order 0  to $\dbR^n\backslash\{0\}$ with respect to
$\xi$. Denote by $R$, $S$ and $J$ the SPDOs, whose symbols are
$r^*$, ${r^*}^{-1}$ and $j^*$, respectively. Then
$J\in\mathcal{L}^1_{\infty}(U)$ is diagonal. Let $Z=RY$. By the
assumption (H2), if an element of the diagonal of $J$ is $A_1+i
B_1$, then either $B_1=0$ or $B_1$ is an elliptic SPDO of order $(1,
\infty)$ with a real symbol.

\medskip

\noindent{\bf Step 3} We give a new Carleman-type estimate for a
SPDO of order $(1, \infty)$.
\begin{lemma}\label{6}
Suppose that $A_1$ and $B_1$ are two SPDOs of order $(1, \infty)$
and  their symbols are real. If $B_1=0$ or $B_1$ is elliptic, and
$z\in L^2_{\mathcal{F}}(\Omega; C([0, T]; H^1(\dbR^n)))$ is an
$H^1(\dbR^n)$-valued semimartingale satisfying $z(0)=z(T)=0$ a.s.,
then for sufficiently small $\mu^{-1}$ and $T$, it holds
\begin{eqnarray}\label{10}
\begin{array}{rl}
&\displaystyle\dbE\int^T_0e^{\mu(t-T)^2}|z|_{L^2(\dbR^n)}^2dt
+\displaystyle\frac{1}{\mu} \dbE\int^T_0e^{\mu(t-T)^2}|\mu(t-T)z-
B_1(t)z|^2_{L^2(\dbR^{n})}dt\\[3mm]
&\leq\displaystyle\frac{4}{\mu} \Re
\dbE\int^T_0\int_{\dbR^n}e^{\mu(t-T)^2}\left[\frac{1}{i}dz-A_1(t)zdt-i
B_1(t)zdt\right]\cdot\overline{[i\mu(t-T)z(t)-i
B_1(t)z]}dx \\[3mm]
&\quad-\displaystyle\frac{2}{\mu}\Im
\dbE\int^T_0\int_{\dbR^n}e^{\mu(t-T)^2}\left[\frac{1}{i}dz-A_1(t)zdt-i
B_1(t)zdt\right]\cdot\overline{(B_1(t)-B_1^*(t))z} dx\\[3mm]
&\ds\q-2\dbE\int^T_0\int_{\dbR^n}(t-T)e^{\mu(t-T)^2}|dz|^2dx-\frac{2}{\mu}\Re
\dbE\int^T_0e^{\mu(t-T)^2}(dz, B_1(t)(dz))_{L^2(\dbR^n)},
\end{array}
\end{eqnarray}
where $B_1^*$ denotes the conjugate operator of $B_1$.
\end{lemma}

\noindent {\bf Proof. } Set $\theta=e^{\frac{\mu(t-T)^2}{2}}$ and
$\varphi=\theta z$. Then it is easy to show that
\begin{equation}\label{11}
\theta \left[\frac{1}{i}dz-A_1(t)zdt-i
B_1(t)zdt\right]=\frac{1}{i}d\varphi-A_1(t)\varphi
dt+i\mu(t-T)\varphi dt-iB_1(t)\varphi dt.
\end{equation}
Multiplying both sides of (\ref{11}) by
$\overline{i\mu(t-T)\varphi-iB_1(t)\varphi}$, and integrating on
$(0, T)\times \Omega \times\dbR^n$, we obtain that
\begin{eqnarray}\label{12}
\begin{array}{rl}
&\mbox{Re }\displaystyle \dbE\int^T_0\int_{\dbR^n}\theta
\left[\frac{1}{i}dz-A_1(t)z dt-i B_1(t)z
dt\right]\overline{[i\mu(t-T)\varphi-i
B_1(t)\varphi]}dx\\[3mm]
&=\mbox{Re }\dbE\displaystyle\int^T_0\int_{\dbR^n}\left[|
\mu(t-T)\varphi-
B_1(t)\varphi|^2dt-\mu(t-T)\overline{\varphi}d\varphi\right.\\[2mm]
&\quad\quad\quad\quad\quad\quad\quad\left.+i \mu(t-T)A_1(t)\varphi
\cdot\overline{\varphi}dt+\overline{ B_1(t)\varphi}d\varphi-i
A_1(t)\varphi\cdot\overline{B_1(t)\varphi}dt\right]dx .
\end{array}
\end{eqnarray}

We estimate one by one the last four terms in (\ref{12}). First, by
It\^o's formula, it follows that
\begin{equation}\label{13}
-\mu\mbox{Re }\dbE\displaystyle\int^T_0\int_{\dbR^n}(t-T)
\overline{\varphi}d\varphi dx
=\frac{\mu}{2}\dbE\left[\int^T_0|\varphi|^2_{L^2(\dbR^n)}dt+\int^T_0\int_{\dbR^n}(t-T)\th^2|dz|^2dx\right].
\end{equation}
Also, by the $L^2$-estimates of SPDOs, we see that
\begin{eqnarray}\label{14}
\begin{array}{rl}
&\mbox{Re }\displaystyle\dbE\int^T_0\int_{\dbR^n}i
\mu(t-T)A_1(t)\varphi\cdot\overline{\varphi}dx
dt=\dbE\int^T_0\frac{\mu(T-t)}{2 i}((A_1(t)-A_1^*(t))\varphi,
\varphi) dt\\[3mm]
&\geq-C\mu T\displaystyle\dbE\int^T_0|\varphi|^2_{L^2(\dbR^n)}dt,
\end{array}
\end{eqnarray}
here and hereafter, we denote by $(\cdot, \cdot)$ the inner product
in $L^2(\dbR^n)$. Next, we notice that
\begin{eqnarray*}
\begin{array}{rl}
&\mbox{Re }\displaystyle\dbE\int^T_0\int_{\dbR^n}\overline{
B_1(t)\varphi}d\varphi dx\\[3mm]
&=\displaystyle\frac{1}{2}\mbox{Re }\dbE\int^T_0[(d\varphi,
B_1(t)\varphi)-(\varphi, B_{1, t}(t)\varphi)dt-(\varphi,
B_1(t)(d\varphi))-(d\varphi,
B_1(t)(d\varphi))]\\[3mm]
&=-\displaystyle\frac{1}{2}\mbox{Re }\dbE\int^T_0(\varphi, B_{1,
t}(t)\varphi)dt\\[3mm]
&\quad+\displaystyle\frac{1}{2}\mbox{Re }\dbE\int^T_0(d\varphi,
(B_1(t)-B_1^*(t))\varphi)-\displaystyle\frac{1}{2}\mbox{Re
}\dbE\int^T_0(d\varphi, B_1(t)(d\varphi)).
\end{array}
\end{eqnarray*}
Therefore, for sufficiently small $\varepsilon>0$, by H\"{o}lder's
inequality and the $L^2$-boundedness of SPDOs, it follows that
\begin{eqnarray}\label{15}
\begin{array}{rl}
&\mbox{Re }\displaystyle\dbE\int^T_0\int_{\dbR^n}\overline{
B_1(t)\varphi}d\varphi dx\\[3mm]
&\geq -C\dbE\displaystyle\int^T_0 |\varphi|_{L^2(\dbR^n)}|B_{1,
t}(t)\varphi|_{L^2(\dbR^n)}dt\\[3mm]
&\quad-\displaystyle\frac{1}{2}\mbox{Im
}\dbE\displaystyle\int^T_0\left(\frac{1}{i}d\varphi-A_1(t)\varphi
dt,
(B_1(t)-B_1^*(t))\varphi\right)\\[3mm]
&\quad-\displaystyle\frac{1}{2}\mbox{Im }\dbE\int^T_0(A_1(t)\varphi,
(B_1(t)-B_1^*(t))\varphi)dt-\frac{1}{2}\mbox{Re
}\dbE\int^T_0(d\varphi,
B_1(t)(d\varphi))\\[3mm]
&=-C\displaystyle\dbE\int^T_0 |\varphi(t)|_{L^2(\dbR^n)}|B_{1,
t}(t)\varphi|_{L^2(\dbR^n)}dt\\[3mm]
&\quad+\displaystyle\frac{1}{2}\mbox{Im }\dbE\int^T_0\left(\theta
\left[\frac{1}{i}dz-A_1(t)zdt-i
B_1(t)zdt\right]-i\mu(t-T)\varphi dt\right.\\[3mm]
&\left.\quad+iB_1(t)\varphi dt,
(B_1(t)-B_1^*(t))\varphi(t)\right)\\[3mm]
&\quad-\displaystyle\frac{1}{2}\mbox{Im }\dbE\int^T_0(A_1(t)\varphi,
(B_1(t)-B_1^*(t))\varphi)dt-\frac{1}{2}\mbox{Re
}\dbE\int^T_0(d\varphi,
B_1(t)(d\varphi))\\[3mm]
&\geq -\varepsilon\displaystyle\dbE\int^T_0
|\varphi|_{H^1(\dbR^n)}^2dt-C(\varepsilon)\displaystyle\dbE\int^T_0
|\varphi|^2_{L^2(\dbR^n)}dt-\displaystyle\frac{1}{2}\dbE\int^T_0|\mu(t-T)\varphi
-B_1(t)\varphi|^2_{L^2(\dbR^n)}dt\\[3mm]
&\quad+\displaystyle\frac{1}{2}\mbox{Im
}\dbE\int^T_0\theta^2\left(\frac{1}{i}dz-A_1(t)z dt-i
B_1(t)zdt, (B_1(t)-B_1^*(t))z\right)\\[3mm]
&\quad-\displaystyle\frac{1}{2}\mbox{Re }\dbE\int^T_0\th^2(dz,
B_1(t)(dz)).
\end{array}
\end{eqnarray}
Moreover, for sufficiently small $\varepsilon>0$, it holds
\begin{eqnarray}\label{16}
\begin{array}{rl} &-\mbox{Re
}\displaystyle\dbE\int^T_0\int_{\dbR^n}
iA_1(t)\varphi\cdot\overline{B_1(t)\varphi}dx
dt=\mbox{Im }\dbE\int^T_0(\varphi,
A_1^*(t)B_1(t)\varphi)dt\\[3mm]
&=\displaystyle\dbE\int^T_0\frac{(\varphi,
A_1^*(t)B_1(t)\varphi)-(A_1^*(t)B_1(t)\varphi,
\varphi)}{2i}dt\\[3mm]
&=\dbE\displaystyle\int^T_0\frac{(\varphi,
(A_1^*(t)B_1(t)-(A_1^*(t)B_1(t))^*)\varphi)}{2i}dt\\[3mm]
&\geq
-C(\varepsilon)\dbE\displaystyle\int^T_0|\varphi|_{L^2(\dbR^n)}^2dt-
\varepsilon\dbE \displaystyle\int^T_0|\varphi|_{H^1(\dbR^n)}^2dt.
\end{array}
\end{eqnarray}

By (\ref{12})--(\ref{16}), we end up with
\begin{eqnarray}\label{17}
\begin{array}{rl}
&\mbox{Re }\displaystyle \dbE\int^T_0\int_{\dbR^n}\theta
\left[\frac{1}{i}dz-A_1(t)zdt-i
B_1(t)zdt\right]\overline{[i\mu(t-T)\varphi-i
B_1(t)\varphi]}dx\\[3mm]
&\geq \displaystyle\frac{1}{2}\dbE\int^T_0\int_{\dbR^n}|\mu(t-T)
\varphi-B_1\varphi|^2dt+\frac{\mu}{2}\dbE\left[\int^T_0|\varphi|^2_{L^2(\dbR^n)}dt
+\int^T_0\int_{\dbR^n}(t-T)\th^2|dz|^2dx\right]
\\[3mm]
&\quad\ds -C\mu
T\dbE\int^T_0|\varphi|^2_{L^2(\dbR^n)}dt-\varepsilon\displaystyle\dbE\int^T_0
|\varphi|^2_{H^1(\dbR^n)}dt-C(\varepsilon)\dbE\int^T_0
|\varphi|^2_{L^2(\dbR^n)}dt\\[3mm]
&\quad+\displaystyle\frac{1}{2}\mbox{Im
}\dbE\int^T_0\theta^2\left(\frac{1}{i}dz-A_1(t)zdt-i B_1(t)zdt,
(B_1(t)-B_1^*(t))z\right)\\[3mm]
&\quad-\displaystyle\frac{1}{2}\mbox{Re }\dbE\int^T_0\th^2(dz,
B_1(dz)).
\end{array}
\end{eqnarray}
If $B_1$ is an elliptic SPDO, by Theorem \ref{32}, it follows that
\begin{eqnarray}\label{18}
\begin{array}{rl}
&\dbE\displaystyle\int^T_0|\varphi|^2_{H^1(\dbR^n)}dt\leq
C\dbE\displaystyle\int^T_0
|B_1(t)\varphi|^2_{L^2(\dbR^n)}dt+C\dbE\displaystyle\int^T_0
|\varphi|^2_{L^2(\dbR^n)}dt\\[3mm]
&\leq
C\dbE\displaystyle\int^T_0|\mu(t-T)\varphi-B_1(t)\varphi|^2_{L^2(\dbR^n)}dt
+C(1+T\mu)\dbE\displaystyle\int^T_0|\varphi|^2_{L^2(\dbR^n)}dt.
\end{array}
\end{eqnarray}
By (\ref{17}) and (\ref{18}), if we take $\mu^{-1}$ and $T$
sufficiently small, then
\begin{eqnarray}\label{19}
\begin{array}{rl}
&\mbox{Re }\displaystyle \dbE\int^T_0\int_{\dbR^n}\theta
\left[\frac{1}{i}dz-A_1(t)zdt-i
B_1(t)zdt\right]\overline{[i\mu(t-T)\varphi-i
B_1(t)\varphi]}dx\\[3mm]
&\geq
\displaystyle\frac{\mu}{4}\dbE\int^T_0|\varphi|^2_{L^2(\dbR^n)}dt
+\displaystyle\frac{1}{4}\dbE\int^T_0\int_{\dbR^n}|\mu(t-T)\varphi
-B_1(t)\varphi|^2dx
dt\\[3mm]
&\quad+\displaystyle\frac{1}{2}\mbox{Im
}\dbE\int^T_0\theta^2\left(\frac{1}{i}dz-A_1(t)zdt-i B_1(t)zdt,
(B_1(t)-B_1^*(t))z\right)\\[3mm]
&\quad\ds+\frac{\mu}{2}\dbE\int^T_0\int_{\dbR^n}(t-T)\th^2|dz|^2dx-\frac{1}{2}\mbox{Re
}\dbE\int^T_0\th^2(dz, B_1(t)(dz)),
\end{array}
\end{eqnarray}
which implies the result of Lemma \ref{6}.

On the other hand, the case of $B_1=0$  can be treated in the same
way.\endpf

\medskip

By Lemma \ref{6}, (\ref{10}) holds for every component of $Z$.
Therefore, we obtain the following result.
\begin{proposition}\label{5}
For sufficiently small $\mu^{-1}$ and $T$, the following inequality
holds:
\begin{eqnarray}\label{9}
\begin{array}{rl}
&\displaystyle\dbE\int^T_0e^{\mu(t-T)^2}|Z(t)|_{L^2(\dbR^n)}^2dt
+\displaystyle\frac{1}{\mu} \dbE\int^T_0e^{\mu(t-T)^2}|\mu(t-T)Z(t)-
B_1(t)Z(t)|^2_{L^2(\dbR^{n})}dt\\[3mm]
&\leq\displaystyle\frac{4}{\mu} \Re
\dbE\int^T_0\int_{\dbR^n}e^{\mu(t-T)^2}\left(\frac{1}{i}dZ
-J(t)Z(t)dt\right)\cdot\overline{[i\mu(t-T)Z(t)-i
B_1(t)Z(t)]}dx \\[3mm]
&\quad-\displaystyle\frac{2}{\mu}\Im
\dbE\int^T_0\int_{\dbR^n}e^{\mu(t-T)^2}\left(\frac{1}{i}dZ
-J(t)Z(t)dt\right)\cdot\overline{(B_1(t)-B_1^*(t))Z(t)} dx\\[3mm]
&\ds\q-2\dbE\int^T_0\int_{\dbR^n}(t-T)e^{\mu(t-T)^2}|dZ|^2dx-\frac{2}{\mu}\Re
\dbE\int^T_0e^{\mu(t-T)^2}(dZ, B_1(t)(dZ)).
\end{array}
\end{eqnarray}
\end{proposition}

\medskip

\noindent {\bf Step 4 }By the result of Proposition \ref{5}, we give
an estimate for the vector $Y$. First of all, we notice that there
exist a constant $C>0$ and a SPDO
$T_0^*\in\mathcal{L}^{0}_\infty(U)$, such that
\begin{eqnarray}\label{24}
\left\{
\begin{array}{rl} &\dbE\displaystyle\int^T_0\theta^2|Y(t)|^2_{L^2(\dbR^n)}dt
\leq C\dbE\displaystyle\int^T_0\theta^2(|Z(t)|^2_{L^2(\dbR^n)}
+|Y(t)|^2_{H^{-1}(\dbR^n)})dt;\\[3mm]
&\ds R\left(\frac{1}{i}dY(t)-AY(t)dt\right)=\frac{1}{i}dZ(t)
-JZ(t)dt+T_0^*Y(t)dt.
\end{array}
\right.
\end{eqnarray}
In fact, by $I_{m\times m}=r^*\cdot {r^*}^{-1}$, there exists a SPDO
$T_{-1}^*\in \mathcal{L}^{-1}_{\infty}(U)$ such that
$I=SR+T^*_{-1}$. Hence,
$$Y=SRY+T^*_{-1}Y=SZ+T^*_{-1}Y.
$$
By the $L^2$-boundedness of SPDOs, the first inequality of
(\ref{24}) holds. On the other hand,
 $$
 \ba{ll}\ds
 R\left(\frac{1}{i}dY-AYdt\right)
 =\frac{1}{i}dZ+\frac{1}{i}R_tYdt-RA(SZ+T^*_{-1}Y)dt\\[3mm]
 \ds=\frac{1}{i}dZ-JZdt+JZdt+\frac{1}{i}R_tYdt-RA(SZ+T^*_{-1}Y)dt\\[3mm]
 \ds=\frac{1}{i}dZ-JZdt+\left(D_tR-RAT^*_{-1}+JR-RASR\right)Ydt,
 \ea
 $$
where $T_0^*=D_tR-RAT^*_{-1}+JR-RASR$ is a SPDO of order $ (0,
\infty)$.

\medskip

Next, by (\ref{24}),  H\"{o}lder's inequality and Proposition
\ref{5}, we get that
\begin{eqnarray}\label{25}
\begin{array}{rl}
&\ds\dbE\int^T_0e^{\mu(t-T)^2}|Y|_{L^2(\dbR^n)}^2dt
+\displaystyle\frac{1}{\mu} \dbE\int^T_0e^{\mu(t-T)^2}|\mu(t-T)RY-
B_1(t)RY|^2_{L^2(\dbR^{n})}dt\\[3mm]
&\leq\displaystyle\frac{4}{\mu} \Re
\dbE\int^T_0\int_{\dbR^n}e^{\mu(t-T)^2}
\left[R\left(\frac{1}{i}dY-AYdt\right)-T_0^*Ydt\right]\cdot\overline{[i\mu(t-T)RY-i
B_1(t)RY]}dx \\[3mm]
&\quad-\displaystyle\frac{2}{\mu}\Im
\dbE\int^T_0\int_{\dbR^n}e^{\mu(t-T)^2}\left[R\left(\frac{1}{i}dY-AYdt\right)-T_0^*Ydt\right]\cdot\overline{(B_1(t)-B_1^*(t))RY} dx\\[3mm]
&\ds\q+C\dbE\int^T_0e^{\mu(t-T)^2}|Y|_{H^{-1}(\dbR^n)}^2dt\\[3mm]
&\ds\q-2\dbE\int^T_0\int_{\dbR^n}(t-T)\th^2|RdY|^2dx-\frac{2}{\mu}\Re
\dbE\int^T_0\th^2(RdY, B_1(t)(RdY))\\[3mm]
&\leq\displaystyle\frac{4}{\mu} \mbox{Re
}\displaystyle\dbE\int^T_0\int_{\dbR^n}
e^{\mu(t-T)^2}R\left(\frac{1}{i}dY-AYdt\right)
\cdot\overline{[i\mu(t-T)RY-i
B_1(t)RY]}dx\\[3mm]
&\quad-\displaystyle\frac{2}{\mu}\mbox{Im
}\displaystyle\dbE\int^T_0\int_{\dbR^n}
e^{\mu(t-T)^2}R\left(\frac{1}{i}dY-AYdt\right)
\cdot\overline{(B_1(t)-B_1^*(t))RY}dx
\\[3mm]
&\quad+\displaystyle\frac{1}{2\mu}\dbE\int^T_0
e^{\mu(t-T)^2}|\mu(t-T)RY-
B_1(t)RY|^2_{L^2(\dbR^n)}dt\\[3mm]
&\quad+\displaystyle\frac{C}{\mu}\dbE\int^T_0e^{\mu(t-T)^2}
|Y|^2_{L^2(\dbR^{n})}dt+C\dbE\int^T_0e^{\mu(t-T)^2}
|Y|^2_{H^{-1}(\dbR^n)}dt\\[3mm]
&\ds\q+CT\dbE\int^T_0\int_{\dbR^n}\th^2|RF|^2dxdt-\frac{2}{\mu}\Re
\dbE\int^T_0\th^2(RF, B_1(t)RF)dt
\end{array}
\end{eqnarray}
Since $B_1$ is an elliptic SPDO and its symbol is real, by Theorem
\ref{101}, we see that
$$
-\frac{2}{\mu}\Re \dbE\int^T_0\th^2(RF, B_1(t)RF)dt\leq
\frac{C}{\mu}\dbE\displaystyle\int^T_0
\theta^2|F|^2_{L^2(\dbR^n)}dt.
$$
It follows that
\begin{eqnarray*}
\begin{array}{rl}
&\ds\dbE\int^T_0e^{\mu(t-T)^2}|Y|_{L^2(\dbR^n)}^2dt
+\displaystyle\frac{1}{\mu} \dbE\int^T_0e^{\mu(t-T)^2}|\mu(t-T)RY-
B_1(t)RY|^2_{L^2(\dbR^{n})}dt\\[3mm]
&\leq \displaystyle\frac{4}{\mu} \mbox{Re
}\displaystyle\dbE\int^T_0\int_{\dbR^n}e^{\mu(t-T)^2}Rf\cdot
\overline{[i\mu(t-T)RY-i
B_1(t)RY]}dtdx\\[3mm]
&\quad-\displaystyle\frac{2}{\mu}\mbox{Im
}\displaystyle\dbE\int^T_0\int_{\dbR^n}e^{\mu(t-T)^2}Rf\cdot
\overline{(B_1(t)-B_1^*(t))RY}
dtdx\\[3mm]
&\quad+\displaystyle\frac{1}{2\mu}\dbE\int^T_0e^{\mu(t-T)^2}|\mu(t-T)RY-
B_1(t)RY|^2_{L^2(\dbR^n)}dt\\[3mm]
&\quad+\displaystyle\frac{C}{\mu}\dbE\int^T_0e^{\mu(t-T)^2}
|Y|^2_{L^2(\dbR^{n})}dt+C\dbE\int^T_0e^{\mu(t-T)^2}
|Y|^2_{H^{-1}(\dbR^n)}dt\\[3mm]
&\ds\q+C(T+\frac{1}{\mu})\dbE\int^T_0e^{\mu(t-T)^2}|F|^2_{L^2(\dbR^n)}dt.
\end{array}
\end{eqnarray*}
Therefore, by H\"{o}lder's inequality, we have that
\begin{eqnarray*}
\begin{array}{rl}
&\ds\dbE\int^T_0e^{\mu(t-T)^2}|Y|_{L^2(\dbR^n)}^2dt
+\displaystyle\frac{1}{\mu} \dbE\int^T_0e^{\mu(t-T)^2}|\mu(t-T)RY-
B_1(t)RY|^2_{L^2(\dbR^{n})}dt\\[3mm]
&\leq\displaystyle\frac{3}{4\mu}\dbE
\int^T_0e^{\mu(t-T)^2}|i\mu(t-T)RY-i
B_1(t)RY|^2_{L^2(\dbR^n)}dt+\displaystyle\frac{C}{\mu}
\dbE\int^T_0e^{\mu(t-T)^2}|f|^2_{L^2(\dbR^n)}dt\\[3mm]
&\quad+\displaystyle\frac{C}{\mu}\dbE\int^T_0e^{\mu(t-T)^2}
|Y|^2_{L^2(\dbR^{n})}dt+Cr^2\dbE\displaystyle\int^T_0e^{\mu(t-T)^2}
|Y|^2_{L^2(\dbR^n)}dt\\[3mm]
&\quad+C(T+\displaystyle\frac{1}{\mu})
\dbE\displaystyle\int^T_0e^{\mu(t-T)^2}|F|^2_{L^2(\dbR^n)}dt.
\end{array}
\end{eqnarray*}
Then, for sufficiently small $r$, $T$ and $\mu^{-1}$ , we obtain
that
\begin{eqnarray}\label{26}
\begin{array}{rl}
&\dbE\displaystyle\int^T_0e^{\mu(t-T)^2}|Y|_{L^2(\dbR^n)}^2dt
+\displaystyle\frac{1}{\mu} \dbE\int^T_0e^{\mu(t-T)^2}|\mu(t-T)RY-
B_1(t)RY|^2_{L^2(\dbR^{n})}dt\\[3mm]
&\leq \displaystyle\frac{C}{\mu}
\dbE\int^T_0e^{\mu(t-T)^2}|f|^2_{L^2(\dbR^{n})}dt+C(T+\displaystyle\frac{1}{\mu})
\dbE\displaystyle\int^T_0e^{\mu(t-T)^2}|F|^2_{L^2(\dbR^n)}dt.
\end{array}
\end{eqnarray}
By the definition of $Y$,
$\displaystyle\sum_{|\alpha|<m}|(D^\alpha_{t, x} (\zeta
u))(t)|^2_{L^2(\dbR^n)}\leq C|Y(t)|^2_{L^2(\dbR^n)}$. Notice that
(\ref{26}) implies (\ref{7}).

\subsection{Proof of Theorem 5.2}

The idea of the proof of Theorem 5.2 is similar to that used in
Theorem \ref{x1}. We only need to prove Proposition \ref{3} for a
strong solution $u$ of equation (\ref{x2}). Different from Theorem
\ref{x1}, since the multiplicity of a complex root may be two, the
matrix which makes the principal symbol $\sigma(A_0)$ become a
Jordan canonical form is only defined locally. Therefore, by a
localization technique, we get the desired estimate. In the
following, we give a sketch of the proof of Theorem 5.2.

\medskip

\noindent{\bf Step 1 }First, under the conditions (H1), (H3) and
(H4), there exist a finite covering $\{G_k\}_{k=1}^{N_1}$
$(N_1\in\dbN)$ of $\{\ \xi\in\dbR^n; \ |\xi|=1\}$, and sufficiently
small $r$ and $T$ such that
$$r^*_k\cdot \sigma(A_0)\cdot {r^*_k}^{-1}=j^*_k, \ \ \ \mbox{for
a.e. }(t, \omega)\in(0, T)\times\Omega\mbox{ and any }(x, \xi)\in
B_r\times G_k,$$ where $r^*_k$, ${r^*_k}^{-1}\in
\bigcap\limits_{j\in\dbN}L^\infty_{\mathcal{F}}(0, T; C^j(B_r\times
G_k))$ and $j^*_k$ is a Jordan canonical form of $\sigma(A_0)$.
Denote by $\{\widetilde{\varphi}_k\}_{k=1}^{N_1}$ a finite number of
smooth functions such that $\{\widetilde{\varphi}_k^2\}_{k=1}^{N_1}$
is the partition of unity for $\{G_k\}_{k=1}^{N_1}$ and
$\sum\limits_{k=1}^{N_1}\widetilde{\varphi}_k^2=1$ on $\{\
\xi\in\dbR^n; \ |\xi|=1\}$. After $\widetilde{\varphi}_k$ is
extended homogeneously of order 0 to $\dbR^n\backslash\{0\}$ with
respect to $\xi$, we denote by $\Phi_k$ the associated
pseudo-differential operator determined by symbol $\varphi_k$ and
write $Y_k=\Phi_k Y$. Also, we choose smooth functions
$\widetilde{\psi}_k: \{\ \xi\in\dbR^n; \ |\xi|=1\}\rightarrow G_k$
satisfying that $\widetilde{\psi}_k(\xi)=\xi$ on
$\supp\widetilde{\varphi}_k$ and set
$$
a^0_k(t, \omega, x, \xi)=\sigma(A_0)(t, \omega, x,
\widetilde{\psi}_k(\xi)),\ r_k(t, \omega, x, \xi)=r^*_k(t, \omega,
x, \widetilde{\psi}_k(\xi)).
$$
Then, after we extend $a^0_k$ ({\it resp.} $r_k$ and $r_k^{-1}$)
homogeneously of order 1 ({\it resp.} order 0) with respect to
$\xi$, we get that $r_k\cdot a^0_k \cdot r_k^{-1}=j_k$ for a.e. $(t,
\omega)\in (0, T)\times\Omega$ and any $(x, \xi)\in B_r\times
\dbR^n\backslash\{0\}$, where $j_k$ is a Jordan canonical form of
$a^0_k$. Denote by $R_k$, $A_k$, $S_k$ and $J_k$ the SPDOs
determined by $r_k, a_k, r_k^{-1}$ and $j_k$, respectively.

\medskip

\noindent{\bf Step 2 }Set $Z_k=R_k Y_k$.  If the geometric
multiplicity and the algebraic multiplicity for every root equal,
then $j_k$ is a diagonal matrix and we can derive the desired result
in the same way as that in the last subsection. Otherwise, by Lemma
\ref{6} and invertibility of elliptic SPDOs, it is easy to see that
the following result holds.
\begin{lemma}\label{124}
Suppose that $A_1$ and $B_1$ are two SPDOs of order $(1, \infty)$
and  their symbols are real. If $B_1=0$ or $B_1$ is elliptic, and
$z_1, z_2\in L^2_{\mathcal{F}}(\Omega; C([0, T]; H^1(\dbR^n)))$ is
an $H^1(\dbR^n)$-valued semimartingale satisfying $z_j(0)=z_j(T)=0$
$(j=1, 2)$ a.s., then there exists a constant $C=C(B_1, n)$ such
that for sufficiently small $\mu^{-1}$ and $T$, it holds
\begin{eqnarray*}
\begin{array}{rl}
&\displaystyle\dbE\int^T_0e^{\mu(t-T)^2}|z_1|_{L^2(\dbR^n)}^2dt
+\displaystyle\frac{1}{\mu} \dbE\int^T_0e^{\mu(t-T)^2}|\mu(t-T)z_1-
B_1(t)z_1|^2_{L^2(\dbR^{n})}dt\\[3mm]
&+\displaystyle\dbE\int^T_0e^{\mu(t-T)^2}|z_2|_{L^2(\dbR^n)}^2dt
+\displaystyle\frac{1}{\mu} \dbE\int^T_0e^{\mu(t-T)^2}|\mu(t-T)z_2-
B_1(t)z_2|^2_{L^2(\dbR^{n})}dt\\[3mm]
&\leq\displaystyle\frac{8}{\mu} \Re
\dbE\int^T_0\int_{\dbR^n}e^{\mu(t-T)^2}\left[\frac{1}{i}dz_1-A_1(t)z_1dt-i
B_1(t)z_1dt+\Lambda z_2dt\right]\cdot\overline{[i\mu(t-T)z_1-i
B_1(t)z_1]}dx \\[3mm]
&\quad-\displaystyle\frac{4}{\mu}\Im
\dbE\int^T_0\int_{\dbR^n}e^{\mu(t-T)^2}\left[\frac{1}{i}dz_1-A_1(t)z_1dt-i
B_1(t)z_1dt+\Lambda z_2dt\right]\cdot\overline{(B_1(t)-B_1^*(t))z_1} dx\\[3mm]
&\ds\q-4\dbE\int^T_0\int_{\dbR^n}(t-T)e^{\mu(t-T)^2}|dz_1|^2dx-\frac{4}{\mu}\Re
\dbE\int^T_0e^{\mu(t-T)^2}(dz_1, B_1(t)(dz_1))_{L^2(\dbR^n)}\\[3mm]
&\quad+\displaystyle\frac{8C}{\mu} \Re
\dbE\int^T_0\int_{\dbR^n}e^{\mu(t-T)^2}\left[\frac{1}{i}dz_2-A_1(t)z_2dt-i
B_1(t)z_2dt\right]\cdot\overline{[i\mu(t-T)z_2-i
B_1(t)z_2]}dx \\[3mm]
&\quad-\displaystyle\frac{4C}{\mu}\Im
\dbE\int^T_0\int_{\dbR^n}e^{\mu(t-T)^2}\left[\frac{1}{i}dz_2-A_1(t)z_2dt-i
B_1(t)z_2dt\right]\cdot\overline{(B_1(t)-B_1^*(t))z_2} dx\\[3mm]
&\ds\q-4C\dbE\int^T_0\int_{\dbR^n}(t-T)e^{\mu(t-T)^2}|dz_2|^2dx-\frac{4C}{\mu}\Re
\dbE\int^T_0e^{\mu(t-T)^2}(dz_2, B_1(t)(dz_2))_{L^2(\dbR^n)}.
\end{array}
\end{eqnarray*}
\end{lemma}

\noindent{\bf Sketch of the proof. }By Lemma \ref{6}, it remains to
estimate the following two terms:
\begin{eqnarray*}
&&\mathcal{I}=-\displaystyle\frac{4}{\mu} \Re
\dbE\int^T_0\int_{\dbR^n}e^{\mu(t-T)^2}\Lambda
z_2dt\cdot\overline{[i\mu(t-T)z_1-i
B_1(t)z_1]}dx\\[3mm]
&&\quad\quad+\displaystyle\frac{2}{\mu}\Im
\dbE\int^T_0\int_{\dbR^n}e^{\mu(t-T)^2}\Lambda
z_2dt\cdot\overline{(B_1(t)-B_1^*(t))z_1} dx.
\end{eqnarray*}
It is easy to check that
\begin{eqnarray*}
&&\mathcal{I}\leq \frac{C}{\mu}\dbE\int^T_0\theta^2
\left(|z_2|_{H^1(\dbR^n)}|\mu(t-T)z_1-B_1(t)z_1|_{L^2(\dbR^n)}
+|z_2|_{H^1(\dbR^n)}|z_1|_{L^2(\dbR^n)}\right)dt.
\end{eqnarray*}
Also, similar to (\ref{18}), we see that
\begin{eqnarray*}
&&\dbE\int^T_0\theta^2|z_2|_{H^1(\dbR^n)}^2dt\leq
C\dbE\int^T_0\theta^2\left[|\mu(t-T)z_2-B_1(t)z_2|_{L^2(\dbR^n)}^2
+(1+T\mu)|z_2|_{L^2(\dbR^n)}^2\right]dt.
\end{eqnarray*}
It follows that
\begin{eqnarray}\label{125}
\begin{array}{rl}
&\mathcal{I}\leq
\displaystyle\frac{1}{2\mu}\dbE\int^T_0\theta^2|\mu(t-T)z_1-
B_1(t)z_1|^2_{L^2(\dbR^{n})}dt+\frac{C}{\mu}\dbE\int^T_0\theta^2|\mu(t-T)z_2-
B_1(t)z_2|^2_{L^2(\dbR^{n})}dt\\[3mm]
&\quad+\displaystyle\frac{1}{2\mu}\dbE\int^T_0\theta^2|\mu(t-T)z_2-
B_1(t)z_2|^2_{L^2(\dbR^{n})}dt
+\frac{C}{\mu}(1+T\mu)\dbE\int^T_0\theta^2|z_2|^2_{L^2(\dbR^n)}dt\\[3mm]
&\quad+\displaystyle\frac{C}{\mu}\dbE\int^T_0\theta^2|z_1|^2_{L^2(\dbR^n)}dt+
\frac{C}{\mu}(1+T\mu)\dbE\int^T_0\theta^2
(|z_1|^2_{L^2(\dbR^n)}+|z_2|^2_{L^2(\dbR^n)})dt.
\end{array}
\end{eqnarray}
Applying the result of Lemma \ref{6} to $z_1$ and $z_2$,
respectively, by (\ref{125}), we obtain the desired result.\endpf

\medskip

Lemma \ref{124} implies the following result.
\begin{proposition}
There exists a constant $C>0$, such that for sufficiently small
$\mu^{-1}$ and $T$,
\begin{eqnarray}\label{9}
\begin{array}{rl}
&\displaystyle\dbE\int^T_0e^{\mu(t-T)^2}|Z_k(t)|_{L^2(\dbR^n)}^2dt
+\displaystyle\frac{1}{\mu}
\dbE\int^T_0e^{\mu(t-T)^2}|\mu(t-T)Z_k(t)-
B_1(t)Z_k(t)|^2_{L^2(\dbR^{n})}dt\\[3mm]
&\leq\displaystyle\frac{C}{\mu}
\left|\dbE\int^T_0\int_{\dbR^n}e^{\mu(t-T)^2}\left(\frac{1}{i}dZ_k
-J_k(t)Z_k(t)dt\right)\cdot\overline{[i\mu(t-T)Z_k(t)-i
B_1(t)Z_k(t)]}dx\right| \\[3mm]
&\quad+\displaystyle\frac{C}{\mu}
\left|\dbE\int^T_0\int_{\dbR^n}e^{\mu(t-T)^2}\left(\frac{1}{i}dZ_k
-J_k(t)Z_k(t)dt\right)\cdot\overline{(B_1(t)-B_1^*(t))Z_k(t)} dx\right|\\[3mm]
&\ds\q-C\dbE\int^T_0\int_{\dbR^n}(t-T)e^{\mu(t-T)^2}|dZ_k|^2dx-\frac{C}{\mu}\Re
\dbE\int^T_0e^{\mu(t-T)^2}(dZ_k, B_1(t)(dZ_k)).
\end{array}
\end{eqnarray}
\end{proposition}

\medskip

\noindent {\bf Step 3. }Finally, we give an estimate of $Y$, which
implies Proposition \ref{3}. Similar to (\ref{24}), we can show that
there exist a constant $C>0$ and a SPDO $T_{0
k}^*\in\mathcal{L}^{0}_\infty(U)$, such that
\begin{eqnarray*}
\left\{
\begin{array}{rl} &\dbE\displaystyle\int^T_0\theta^2|Y_k(t)|^2_{L^2(\dbR^n)}dt
\leq C\dbE\displaystyle\int^T_0\theta^2(|Z_k(t)|^2_{L^2(\dbR^n)}
+|Y_k(t)|^2_{H^{-1}(\dbR^n)})dt;\\[3mm]
&\ds R_k\left(\frac{1}{i}dY_k(t)-AY_k(t)dt\right)=\frac{1}{i}dZ_k(t)
-J_kZ_k(t)dt+T_{0 k}^*Y_k(t)dt.
\end{array}
\right.
\end{eqnarray*}
Therefore,  similar to (\ref{25}), we see that
\begin{eqnarray*}
\begin{array}{rl}
&\ds\dbE\int^T_0e^{\mu(t-T)^2}|Y_k|_{L^2(\dbR^n)}^2dt
+\displaystyle\frac{1}{2\mu}
\dbE\int^T_0e^{\mu(t-T)^2}|\mu(t-T)R_kY_k-
B_1(t)R_kY_k|^2_{L^2(\dbR^{n})}dt\\[3mm]
&\leq\displaystyle\frac{C}{\mu}
\left|\displaystyle\dbE\int^T_0\int_{\dbR^n}
e^{\mu(t-T)^2}R_k\left(\frac{1}{i}dY_k-AY_kdt\right)
\cdot\overline{[i\mu(t-T)R_kY_k-i
B_1(t)R_kY_k]}dx\right|\\[3mm]
&\quad+\displaystyle\frac{C}{\mu}\left|\displaystyle\dbE\int^T_0\int_{\dbR^n}
e^{\mu(t-T)^2}R_k\left(\frac{1}{i}dY_k-A_kY_kdt\right)
\cdot\overline{(B_1(t)-B_1^*(t))R_kY_k}dx\right|
\\[3mm]
&\quad+\displaystyle\frac{C}{\mu}\dbE\int^T_0e^{\mu(t-T)^2}
|Y_k|^2_{L^2(\dbR^{n})}dt+C\dbE\int^T_0e^{\mu(t-T)^2}
|Y_k|^2_{H^{-1}(\dbR^n)}dt\\[3mm]
&\ds\q-\frac{C}{\mu}\dbE\int^T_0\int_{\dbR^n}e^{\mu(t-T)^2}(t-T)|R_kd
Y_k|^2dx-\frac{C}{\mu}\Re \dbE\int^T_0\th^2(R_kdY_k, B_1(t)R_k dY_k)
\end{array}
\end{eqnarray*}

By the definition of $Y_k$, it is easy to show that
$|Y|_{H^s(\dbR^n)}=\sum\limits_{k=1}^{N_1}|Y_k|_{H^s(\dbR^n)}$ for
any $s\in\dbR$ and $A\Phi_k Y-\Phi_k A Y=T_{0 k}Y$, where $T_{0
k}\in \mathcal{L}^0_\infty(U)$. Hence, it follows that
\begin{eqnarray*}
\begin{array}{rl}
&\ds\dbE\int^T_0e^{\mu(t-T)^2}|Y|_{L^2(\dbR^n)}^2dt\\[3mm]
&\leq\displaystyle\frac{C}{\mu}
\dbE\int^T_0e^{\mu(t-T)^2}|f|^2_{L^2(\dbR^n)}dt+\displaystyle\frac{C}{\mu}\dbE\int^T_0e^{\mu(t-T)^2}
|Y|^2_{L^2(\dbR^{n})}dt+Cr^2\dbE\displaystyle\int^T_0e^{\mu(t-T)^2}
|Y|^2_{L^2(\dbR^n)}dt\\[3mm]
&\quad+C(T+\displaystyle\frac{1}{\mu})
\dbE\displaystyle\int^T_0e^{\mu(t-T)^2}|F|^2_{L^2(\dbR^n)}dt.
\end{array}
\end{eqnarray*}
If we take $r$ and $\displaystyle\frac{1}{\mu}$ sufficiently small,
then the above inequality implies Proposition \ref{3}.

\begin{remark}
In \cite{9}, for deterministic partial differential equations of
order $m$, principal symbols are required to have $C^\infty$
coefficients. However, by the proofs of Theorem 5.1 and Theorem 5.2,
we notice that the coefficients of equation (\ref{x2}) should only
belong to the space $\mathcal{X}_\infty$. This is because we
actually consider $t$ and $\omega$ as two parameters. In fact, we
find that in the deterministic case, all pseudo-differential
operators appeared in the proof of uniqueness theorem may  regard
the variable $t$ as a parameter. Therefore, in the classical
Calder\'{o}n uniqueness theorem, it is sufficient that the
coefficients of principal symbols are required to be $C^1$ with
respect to $t$.
\end{remark}

\end{document}